\documentclass[11pt]{elsarticle}
\usepackage[margin=1.2in]{geometry}
\usepackage[table]{xcolor}
\usepackage[%
pdftitle={The mixing time of switch Markov chains: a unified
approach},
pdfauthor={Péter L. Erdős; Catherine Greenhill; Tamás Róbert Mezei;
István Miklós; Dániel Soltész; Lajos Soukup},
pdfkeywords={degree sequences;
realizations; switch Markov chain; rapidly mixing; MCMC; Sinclair's
multi-commodity flow method; P-stability; strong stability; power-law
distribution},
pdfsubject={Mathematics},unicode=true%
]{hyperref}
\usepackage{microtype}
\usepackage{caption}
\makeatletter
\def\ps@pprintTitle{%
	\let\@oddhead\@empty
	\let\@evenhead\@empty
	\def\@oddfoot{\reset@font\hfil\thepage\hfil}
	\let\@evenfoot\@oddfoot
}

\renewcommand\subsection{\@startsection{subsection}{2}{\z@}%
           {12\p@ \@plus 6\p@ \@minus 3\p@}%
           {3\p@ \@plus 6\p@ \@minus 3\p@}%
           {\normalfont\normalsize\bfseries}}

\newcommand{\customlabel}[2]{%
   \protected@write \@auxout {}{\string \newlabel {#1}{{#2}{\thepage}{#2}{#1}{}} }%
   \hypertarget{#1}{}
}
\makeatother

\usepackage{amsmath}
\usepackage{amsthm}
\usepackage{amssymb}
\usepackage{enumerate}
\usepackage{pgfplots}
\pgfplotsset{compat=newest}
\usepackage{tikz}
\usetikzlibrary[decorations.pathmorphing,decorations.pathreplacing,matrix,calc]

\usepackage{float}
\usepackage{dsfont}
\usepackage{enumitem}
\usepackage{centernot}

\usepackage{mdframed}

\newcommand\DOI[1]{{\tt DOI:#1}}

\usepackage{verbatim}

\def\BrbG{balanced red-blue graph}
\numberwithin{equation}{section}
\newtheorem{theorem}{Theorem}[section]
\newtheorem*{theorem*}{Theorem}
\newtheorem{conjecture}[theorem]{Conjecture}
\newtheorem{lemma}[theorem]{Lemma}

\newtheorem{corollary}[theorem]{Corollary}
\theoremstyle{definition}
\newtheorem{definition}[theorem]{Definition}

\def\qed{{\hfill$\Box$}}
\def\G{\mathbb{G}}
\def\S{\mathbb{S}}
\def\C{\mathcal{C}}

\newcommand\old[1]{}
\def\B{\mathbb{B}}

\newcommand{\posp}[1]{(#1)^+}
\newcommand{\posm}[1]{(#1)^-}
\def\Pos{\mathit{Pos}}
\newcommand{\mbb}[1]{\mathbb{#1}}

\def\bd{\mathbf{d}}
\def\bD{\mathbf{D}}

\def\F{\mathcal{F}}

\def\m{\mu}

\def\UC{unconstrained degree}
\def\<{\left\langle}
\def\>{\right\rangle}

\DeclareMathOperator{\Prob}{Prob}
\DeclareMathOperator{\poly}{poly}

\usepackage[ruled,section]{algorithm}
\usepackage{algpseudocode}
\algnewcommand{\IIf}[1]{\State\algorithmicif\ #1\ \algorithmicthen}
\algnewcommand{\EndIIf}{\unskip\ \algorithmicend\ \algorithmicif}
\algblock[indent]{BeginIndent}{EndIndent}{}{}

\begin{document}

\begin{frontmatter}
	\title{The mixing time of switch Markov chains:\texorpdfstring{\\}{} a
	unified approach}

	\author[renyi]{Péter L. Erdős\texorpdfstring{\fnref{elp}}{}}
	\author[sydney]{Catherine Greenhill\texorpdfstring{\fnref{green}}{}}
	\author[renyi]{Tamás Róbert Mezei\texorpdfstring{\fnref{elp}}{}}
	\author[renyi]{István Miklós\texorpdfstring{\fnref{elp,miklos}}{}}
	\author[renyi]{Dániel Soltész\texorpdfstring{\fnref{elp,soltesz}}{}}
	\author[renyi]{Lajos Soukup\texorpdfstring{\fnref{soukup}}{}}
	\address[renyi]{Alfréd Rényi Institute of Mathematics, Reáltanoda u 13--15 Budapest, 1053 Hungary\\
	\texttt{$<$erdos.peter,mezei.tamas.robert,miklos.istvan,\\ soltesz.daniel,soukup.lajos$>$@renyi.hu}}
	\address[sydney]{School of Mathematics and Statistics, UNSW Sydney, NSW 2052, Australia \\ \texttt{c.greenhill@unsw.edu.au} }
	\fntext[elp]{The authors were supported in part by the National Research, Development and Innovation Office -- NKFIH grant K~116769, KH~126853, SNN~135643, K~132696}
	\fntext[miklos]{IM was supported in part by the National Research,
	Development and Innovation Office -- NKFIH grant SNN~116095 and SNN~135643.}
	\fntext[green]{Research supported by the Australian Research Council, Discovery Project DP190100977.}
	\fntext[soukup]{LS was supported in part by the National Research, Development and Innovation Office -- NKFIH grants  K~113047 and K~129211.}
	\fntext[soltesz]{DS was supported in part by the National Research, Development and Innovation Office -- NKFIH grans K~120706 and KH~130371.}

	\begin{abstract}
		Since 1997 a considerable effort has been spent to study  the
		mixing time of \textbf{switch} Markov chains on the realizations
		of graphic degree sequences of simple graphs. Several results
		were proved on rapidly mixing Markov chains on unconstrained,
		bipartite, and directed sequences, using different mechanisms.

		The aim of this paper is to unify these approaches. We will
		illustrate the strength of the unified method by showing that on
		any $P$-stable family of  unconstrained/bipartite/directed
		degree sequences the switch Markov chain is rapidly mixing. This
		is a common generalization of every known result that shows the
		rapid mixing nature of the switch Markov chain on a region of
		degree sequences. Among the applications of this general result
		is an almost uniform sampler for power-law and heavy-tailed
		degree sequences. Another application shows that the switch
		Markov chain on the degree sequence of an Erdős-Rényi random
		graph $G(n,p)$ is asymptotically almost surely rapidly mixing if
		$p$ is bounded away from 0 and 1 by at least $\frac{5\log
		n}{n-1}$.
	\end{abstract}
\end{frontmatter}

\section{Introduction}\label{sec:intro}

An important problem in network science is to algorithmically construct typical
instances of networks with predefined properties. In particular, special
attention has been devoted to sampling simple graphs with a given degree
sequence. In this paper only graphs without parallel edges and loops are
considered and we restrict our study to degree sequences which have at least one
realization (\textbf{graphic}). A realization of a degree sequence $\bd\in
\mathbb{N}^n$ is a graph on the vertex set $[n]=\{1,\ldots,n\}$ whose degree
sequence $d(G)$ is equal to $\bd$ (i.e., the graphs are labeled).  Sometimes we
use a different labeling of the vertex set $V=\{v_1,\ldots,v_n\}$, where $V\ni
v_i\longleftrightarrow i\in [n]$. We will use the more general notation $[a,b]$
and $[a,b)$ for the closed-closed and closed-open intervals of integers between
$a$ and $b$. We study the three most common degree sequence types: bipartite
degree sequences, directed degree sequences and the usual degree sequences which
we call \textbf{unconstrained} degree sequences (these are the degree sequences
of simple graphs).

In 1997 Kannan, Tetali, and Vempala (\cite{KTV97}) proposed the use of the
\textbf{switch Markov chain} (also known as the swap chain~\cite{MES}) for
uniformly sampling realizations of a degree sequence. For all three degree
sequence types, the switch Markov chain can be thought of as the Markov chain of
smallest possible modifications. To illustrate this, we give an informal
description of the switch Markov chain on \UC\ sequences. If $G_1,G_2$ are two
realizations of the same \UC\ sequence, it is easy to see that the minimum size
of the symmetric difference $E(G_1) \triangle E(G_2)$ is four. We say that $G_1$
and $G_2$ differ by a switch if this symmetric difference is exactly four. The
states of the switch Markov chain are the realizations of the degree sequence
and the probability of going from realization $G_1$ to $G_2$ is nonzero if and
only if they differ by a switch (and this nonzero quantity is independent of
$G_1$ and $G_2$). For the precise definition of this chain, and for the
definition of the chains for other degree sequence types, we refer the reader to
Section~\ref{sec:sinclair} (unconstrained and bipartite) and
Section~\ref{sec:directed} (directed).

\medskip

Because the symmetric difference of graphs is ubiquitous in this paper, we often
shorten the notation $E(X)\triangle E(Y)$ to $X\triangle Y$.

\medskip

The following conjecture has been named after Kannan, Tetali, and Vempala, in
recognition of their pioneering work.

\begin{conjecture}[the KTV conjecture]\label{conj:rapidmix}
	The switch Markov chain is rapidly mixing for any bipartite, directed,
	or unconstrained degree sequence.
\end{conjecture}

To give some context to the Conjecture, we say that a Markov chain is rapidly
mixing if the distribution on the state space is close in $\ell_1$ norm to the
unique stationary distribution after $\poly(n)$ steps, where $n$ is the size of
the input. In our case, $n$ is the number of vertices (or the number of entries
in the degree sequence). This property means that sampling from the state space
with the stationary distribution is a more or less tractable problem, even if
the state space has exponential size.

\medskip

It is not uncommon that uniformly randomly applied, small local modifications of
combinatorial objects result in rapid mixing. This is the case for solutions of
the 0--1 knapsack problem~\cite{MS04} (Morris and Sinclair, 2004), for the union
of perfect and almost perfect matchings of a graph~\cite{JS89} (Jerrum and
Sinclair, 1989), and two-rowed contingency tables~\cite{DG00} (Dyer and
Greenhill, 2000), for example. In all of these cases, applying the smallest
possible modifications of the respective combinatorial objects randomly, result
in rapid mixing of the corresponding Markov chain.

Although Conjecture~\ref{conj:rapidmix} is still open, there is a series of
results that prove the rapid mixing of the switch Markov chain on various
special degree sequence classes. We summarize these results in a very compact
way in Table~\ref{table:results} without presenting the sometimes lengthy
definitions of the special classes, which can be found in the references
provided. Some rapid mixing results on directed degree sequences work with
directed graphs~\cite{G11} (Greenhill, 2011), while some work in the bipartite
representation with a forbidden perfect matching~\cite{EKMS} (Erdős, Kiss,
Miklós and Soukup, 2015). Since we use the bipartite representation in the
present paper, we will not discuss directed degree sequences until
Section~\ref{sec:directed}.

\begin{table}[H]
	\begin{center}
		\renewcommand*{\arraystretch}{1.5}
		\begin{tabular}{| c | c | c | }
			\hline
			unconstrained deg.\ seq.
			& bipartite deg.\ seq.
			& directed deg.\ seq.
			\\ \hline \hline
			regular~\cite{CDG07}
			& {(half-)}regular~\cite{MES}
			& regular~\cite{G11}
			\\ \hline
			& \multicolumn{2}{|c|}{almost half regular~\cite{EKMS}}
			\\ \hline
			$ \Delta \leq \frac{1}{3}\sqrt{2m}$~\cite{G18}                            & $\Delta \leq \frac{1}{\sqrt{2}} \sqrt{m}$~\cite{linear}       & $ \Delta< \frac{1}{\sqrt{2}}\sqrt{m-4}$~\cite{linear} \\ \hline
			power-law density-                                                        &           \cellcolor{gray!50!white}  &   \cellcolor{gray!50!white}                                                       \\
			bound, $\gamma>2.5$~\cite{G18}                                            &               \cellcolor{gray!50!white}                                                &   \cellcolor{gray!50!white}                                                    \\ \hline
			${(\Delta-\delta+1)}^2\le$                                                & ${(\Delta-\delta)}^2\le$                                      & similar to bipartite case                             \\
			$\le 4\cdot \delta(n-\Delta-1)$                               & $\le \delta(\frac{n}{2}-\Delta)$~\cite{linear-pre, linear}    & \cite{linear-pre,linear}                              \\
			proof in \cite{AK19,AK20}		  & (or Corollary~18 in
			\cite{AK19,AK20})								  &													      \\ \hline
			\cellcolor{gray!50!white} & bipartite Erdős-Rényi \cite{linear-pre,linear} & similar to bipartite case                             \\ \cellcolor{gray!50!white}
						  & $p,\:1-p\ge 4\sqrt{\frac{2\log n}{n}}$ & \cite{linear-pre,linear}                              \\ \hline
			\multicolumn{2}{|c|}{strongly stable degree sequence
		classes~\cite{AK19,AK20}} & \cellcolor{gray!50!white}                                                                                                                       \\ \hline
		\end{tabular}
	\end{center}
	\caption{Some classes of degree sequences for which the rapid mixing of
	the switch Markov chain is was already known when the first
	draft of this paper was published on arXiv. Here $\Delta$ and
	$\delta$ denote the maximum and minimum degrees, respectively.
	Half of the sum of the degrees is $m$, and $n$ is the number of
	vertices. The notation is similar for bipartite and directed degree sequences.
	Some technical conditions are omitted.}\label{table:results}
\end{table}

Notice, that some, but not all of the results came in pairs for unconstrained
and bipartite (directed) degree sequences. The reason for this discrepancy is
the following: while both set of results are based on Sinclair's multicommodity
flow method, one of them has to deal with special circuits (one of the vertices
may be visited at most twice) instead of just cycles in the decomposition of
symmetric differences of two realizations of a degree sequence. The main goal of
this paper to remedy this discrepancy between the machineries used for the
bipartite and unconstrained degree sequences by decomposing into circuits where
each vertex is visited at most twice. Along the way we also give new, much more
transparent proofs for the main results in~\cite{MES} (Miklós, Erdős and Soukup, 2013).

Greenhill and Sfragara suggested exploring the connection between the mixing
rate of the switch Markov chain and stable degree sequences~\cite[Section
1.1]{G18} (Greenhill and Sfragara, 2018).  The first notion of stability for
unconstrained and bipartite degree sequences was $P$-\emph{stability},
introduced by Jerrum and Sinclair~\cite{JS90}.  An unconstrained degree sequence
on $n$ vertices, usually denoted by $\bd$, is an element of ${[n-1]}^n$. The set
of graphs with the degree sequence $\bd$ is denoted by $\G(\bd)$. To make the
following and later stability definitions more readable, let
$\G(\bd):=\emptyset$ for a non-graphic degree sequence $\bd$.
\begin{definition}\label{def:Pstable}
	Let $\mathcal{D}$ be an infinite set of unconstrained degree sequences. We say that $\mathcal{D}$ is \textbf{\boldmath$P$-stable}, if there exists a polynomial $p\in \mathbb{R}\left[x\right]$ such that for any $n\in \mathbb{N}$ and any  degree sequence $\bd \in \mathcal{D}$ on $n$ vertices we have
	\begin{equation*}
		\left|\G(\bd)\cup \left( \bigcup_{x,y\in [n],\ x\neq y}\G(\bd+\mathds{1}_x+\mathds{1}_y) \right)\right|\le p(n)\cdot |\G(\bd)|,
	\end{equation*}
	where $\mathds{1}_x$ is the $x^\text{th}$ unit vector.
\end{definition}
Without proof we state, that the notion of $P$-stability does not change even if we require
\begin{equation}
	\label{eq:statement-without-proof}
	\left|\left\{ G\ |\ G\in\G(\bd'),\ \bd'\in \mathbb{N}^n,\ \ell_1(\bd,\bd')\le 2\right\}\right|\le p(n)\cdot |\G(\bd)|.
\end{equation}
With a bit of care, Definition~\ref{def:Pstable} generalizes to bipartite and
directed degree sequences: we require that the realizations of perturbed degree
sequences are also bipartite or directed, respectively.

Informally, a class of \UC\ sequences is $P$-stable if making slight
perturbations to the degree sequence from a $P$-stable class cannot increase the
number of realizations too much. Thus $P$-stability is a property of the degree
sequences directly, in the sense that it only cares about the number of their
(perturbed) realizations.
\begin{theorem}[proved in Section~\ref{sec:appl}]\label{thm:unified}
	The switch Markov chain is rapidly mixing on $P$-stable unconstrained,
	bipartite, and directed degree sequence classes.
\end{theorem}

For the sake of transparency, we would like to mention that we were aware that a
quirky proof of Theorem~\ref{thm:unified} for the bipartite and directed models
can be read out already from~\cite{linear-pre,linear} (Erdős, Mezei,
Miklós and Soltész, 2018). However, that paper used a complicated and
opaque technique for describing the switch sequence
(Section~\ref{sec:switch-sequence}), which cannot be applied to unconstrained
degree sequences. Instead of publishing an immature proof for the bipartite and
directed cases, we opted to find a unified proof which naturally accommodates
unconstrained degree sequences too.

\medskip

\textbf{Comparison with previous results.}
Amanatidis and Kleer~\cite{AK19,AK20} introuced and studied the concept of
\emph{strongly stable} degree sequences.  In contrast with $P$-stability, strong
stability requires that for any realization of the perturbed degree sequence
there exists a realization of the original degree sequence such that the two
graphs are not too different. The precise definition of strong stability is
given at the beginning of Section~\ref{sec:regions}.  In~\cite{AK19,AK20}, it is
shown that strongly stable degree sequences are $P$-stable. It is not known
whether these concepts of stability are equivalent.

\begin{theorem*}[\cite{AK19,AK20}, Amanatidis and Kleer, 2020]
	The switch Markov chain is rapidly mixing on strongly stable
	unconstrained and bipartite degree sequence classes.
\end{theorem*}

Since~\cite{AK19} appeared, it
has become a commonly held belief that strong stability and $P$-stability are
equivalent. It is mentioned in~\cite{AK19} that the bipartite case does not
immediately imply the directed one, because their proof still needs some work
to accommodate \emph{forbidden edge sets}. Whether the conjecture on the
equivalence of the two notions of stability holds or not,
Theorem~\ref{thm:unified} is a strictly stronger theorem than the one above.

\medskip

The main reason we could unify the proofs for the unconstrained, bipartite, and
directed models is because of a technical novelty introduced in this paper:
extending the $T$-operator of~\cite{MES} (Miklós, Erdős and Soukup, 2013) from
bipartite to unconstrained graphs. The atoms of a decomposition generated by the
$T$-operator are so-called \emph{primitive circuits} (see
Definition~\ref{def:primitive}), along which we recursively construct the
multi-commodity flow required by the Sinclair method. The processing of the
primitive circuits via Algorithm~\ref{alg:sweep} is in turn an extension of the
\textsc{Sweep} algorithm of~\cite{MES}. These refinements require a fairly
extensive and detailed analysis of the structures of the realizations under
study.  In exchange, the unified framework in which we prove
Theorem~\ref{thm:unified} allows us to treat the three models with minimal
branching in the proof.

\medskip

There are two interesting direct consequences of Theorem~\ref{thm:unified}
concerning popular unconstrained random graph models. It turns out that
asymptotically almost surely, the degree sequence of an Erdős-Rényi random graph
$G(n,p)$ is $P$-stable if $p$ is bounded away from 0 and 1 by $\frac{5{\log
n}}{n-1}$, see Corollary~\ref{cor:ER}. Gao and Greenhill~\cite{GG20} show that
power-law distribution-bounded degree sequences with parameter $\gamma > 2$ are
$P$-stable. Consequently, Theorem~\ref{thm:unified} implies that the switch
Markov chain is rapidly mixing on these degree sequences: see
Section~\ref{sec:power}. This gives the first formal verification of the
validity of generating random power-law (distribution-bounded) graphs via the
switch Markov chain~\cite{BA1} (Jia, and Barabási, 2013) and~\cite{BA2} (Yan,
Vértes, Towlson, Chew, Walker, Schafer and Barabási, 2017).

\medskip

The proof of Theorem~\ref{thm:unified} relies on Sinclair's multicommodity flow
method (Section~\ref{sec:sinclair}), which can be described informally in the
case of the switch Markov chain as follows. Suppose that the chain has $N$
states, and let $\G$ be the graph underlying the Markov chain: the vertex set of
$\G$ is the state space of the chain (in our applications, the states are
realizations of a given degree sequence), and two states are adjacent in $\G$ if
the transition probability between them is non-zero in the chain. Sinclair's
multicommodity flow method ensures the rapid mixing of the chain if we can
design a multicommodity flow on $\G$ which transfers a unit amount of
commodities between each pair of vertices (different commodities for different
pairs), and no more than $N\cdot\poly(\log(N))$ amount of commodities go through
any vertex (no vertex is overloaded). Hence most of the present paper is devoted
to the design of a flow and the proof that it does not overload any vertex.

\medskip

\subsection{Some related Markov chain approaches}
In the literature there is at least one other Markov chain application  where the
applied operation is the switch. This studies binary contingency tables with
fixed marginals. Such tables can be considered as bipartite graphs with fixed
degree sequences.

\medskip

\textbf{Switch Markov chain on perfect matchings \emph{(or Diaconis chain)}.}
Diaconis, Graham and Holmes in 2001 (\cite{DGH01}) considered applications of
$(0,1)$-permanents to problems in statistics. The permanent is equal to the
number of the perfect matchings in the corresponding bipartite graph. They study
the switch Markov chain approach on the perfect matchings: the state space is
the set of perfect matchings of the fixed bipartite graphs and the transitions
are generated by switches. It is easy to see that this Markov chain is not
irreducible in general. For example, the hexagon has two distinct perfect
matchings (every second edge from the cycle), but no switch can be applied to it.
In~\cite{DGH01} a structural constraint was imposed on
the bipartite graph making the chain irreducible, and the authors conjectured
that the \emph{switch Markov chain for perfect matchings} is rapidly mixing on
that graph class.

Dyer, Jerrum and Müller in 2017 (\cite{DJM17}) showed that the switch Markov
chain for perfect matchings is irreducible if and only if the bipartite graph is
\emph{chordal}. However, as it turned out, the chain is not rapidly mixing on
all chordal bipartite graphs. The largest hereditary subclass of the chordal
balanced bipartite graphs (so called \emph{monotone balanced bipartite graphs})
for which the switch Markov chain is rapidly mixing was also determined
in~\cite{DJM17}. The proof is based on Sinclair's multicommodity flow
method.

\medskip\textbf{Applications of the switch Markov chain in practice.}
The switch Markov chain is often used in everyday practice. In social sciences
or in ecology it is standard practice to generate an initial example graph with
some required properties (for example, the degree sequence shall obey a
power-law, see Section~\ref{sec:power}) and perturb the initial graph with a
(not too long) series of switch operations. While there is no theoretical
guarantee that the resulting graph is much ``closer'' to a random one than the
initial graph, in practice, the resulting graph is sufficiently random.

\medskip\textbf{Sampling binary contingency tables with simulated annealing.}
Returning to binary contingency tables with arbitrary but fixed marginals, the
first problem is to count and randomly sample them for statistical purposes. As
it can be known from Table~\ref{table:results}, we cannot do that with
confidence using the switch Markov chain. Jerrum, Sinclair and
Vigoda~\cite{JSV04} attacked the problem with simulated annealing in 2004. Their
approach was improved by Bezáková, Bhatnagar and Vigoda in 2006 (\cite{BBV06}),
who used a greedy starting point and subsequently perturbed it with a moderate
number of switches to obtain a faster running time for bipartite graphs.
However, Štefankovič, Vigoda and Wilmes~\cite{STF18} showed that there is no
weighting scheme for which the Jerrum-Sinclair-Vigoda chain mixes rapidly in
general. It is not known whether the switch chain rapidly mixes over the set of
all graphic degree sequences.

\subsection{Structure of the paper}
The remainder of the paper is structured as follows. In Section~\ref{sec:def},
we give a slightly more detailed description of the flow and we present the
necessary graph theoretic tools which we will use to construct paths that will
form the flow. In Section~\ref{sec:sinclair} we give the formal definition of
Sinclair's multicommodity flow method and its simplified version that is
tailored to our needs. In Section~\ref{sec:order} we provide a high level
description of the multicommodity flow.  In Section~\ref{sec:design} we
introduce the $T$-operator and complete the description of the multicommodity
flow.  In Section~\ref{sec:directed} we define the switch Markov chain for
directed degree sequences. Before turning into the home stretch, an auxiliary
structure that tracks the defined flow is analyzed in Section~\ref{sec:aux}.  In
Section~\ref{sec:appl} we finally prove Theorem~\ref{thm:unified} and we also
provide the necessary modifications to deal with bipartite and directed degree
sequences. Lastly, we describe the known $P$-stable regions of degree sequences
in Section~\ref{sec:regions} and present the connections between
Theorem~\ref{thm:unified} and the aforementioned popular graph models.

\section{Definitions and preliminaries, the structure of the sets of realizations}\label{sec:def}

Let us recall some well known notions and notations. Let
$\bd=(d(v_1),\ldots,d(v_n))$  denote a \UC\ sequence on vertex set $V=[n]$ and let
\begin{displaymath}
	\bD=\big(\bigl(d(u_1), \ldots, d(u_{n_1}) \bigr),
	\bigl(d(v_1),\ldots,d(v_{n_2})\bigr)\big)
\end{displaymath}
denote a bipartite degree sequence on the bipartition
$(U,V):=([n_1],[n_1+1,n_1+n_2])$. (For convenience
we assume that $n_1 \ge n_2$ and that $\bD$ can also be considered as an
$n_1+n_2$ long vector.) We will use the notations $\G(\bd)$ and $\G(\bD)$ for
the sets of all realizations of the corresponding degree sequences.

The \textbf{switch} \textbf{operation} exchanges two disjoint edges $ac$ and
$bd$ in the realization $G$ with $ad$ and $bc$ if the resulting configuration
$G':=G-ac-bd+ad+bc$ is again a simple graph (we denote the operation by $ac, bd
\Rightarrow ad,bc$).

For an $ac, bd \Rightarrow ad,bc$ switch operation to be valid, it is necessary
but not always sufficient that both $ac, bd \in E(G)$ and $ad, bc \not\in E(G)$
hold.  For each setting (graphs, bipartite graphs, directed graphs) we will
define a set of \textbf{non-chords}, which are pairs of vertices which are
forbidden from forming edges. See Definition~\ref{def:chordsUCbip} below.  A
pair of vertices which is not a non-chord will be called a \textbf{chord}: such
pairs are allowed to be edges, so they may be inserted or deleted.  We emphasize
that whether or not a pair of vertices forms a \textbf{chord} does not depend on
the current realization.

\medskip

The term chord/non-chord is motivated from Algorithm~\ref{alg:sweep}, which
constructs a switch sequence that exchanges edges and non-edges of a circuit. To
do so, it has to include some of the chords of the circuit in switches. Those
chords that may not appear in a switch (because that would violate the graph
model) are hence called non-chords.

\medskip

Next, we reformulate the definition of the switch operation to avoid inserting
non-chords: an $ac, bd \Rightarrow ad,bc$ switch operation can be applied if
$ac, bd \in E(G)$, $ad, bc \not\in E(G)$, and $ad, bc$ are both chords. We now
define the set of chords and non-chords in the case of unconstrained and
bipartite graph models.

\begin{definition}\label{def:chordsUCbip}
	For simple graphs, the non-chords are exactly the pairs of the form
	$(v,v)$, as loops are forbidden. Because no further constraints have to
	be set, we call their degree sequences \textbf{unconstrained}.  In
	\textbf{bipartite} graphs, $(u,v)$ is a chord if and only if $u$ and $v$
	are in different vertex classes.
\end{definition}
In the case of directed graphs (Section~\ref{sec:directed}), we further restrict
the set of \textbf{chords}.

\medskip

It is a well-known fact that the set of all possible realizations of a
\textbf{graphic} \UC\ sequence is connected under the switch operation. See for
example~\cite{H55} (Havel, 1957) or~\cite{H62} (Hakimi, 1962). It is interesting
to know, however, that the first known proof is from 1891~\cite{pet} (Petersen,
1891).  For bipartite graphs the equivalent results were proved in 1957
in~\cite{G57} (Gale, 1957) and~\cite{R57} (Ryser, 1957). The ``classical''
proofs work through so called ``canonical'' realizations. However, the paths
between different realizations, created in this way, are very far from efficient
for the purpose of applying them in Sinclair's multicommodity-flow method.
Therefore another way has to be designed to select the appropriate paths.

To that end, let us consider two realizations of the same (bipartite or
unconstrained) degree sequence. Those edges, that are present or missing in both
realizations, do not need to be changed. The remaining pairs of vertices, each
of which is an edge in exactly one of the realizations, form the symmetric
difference of the two realizations, usually denoted by $\nabla$. To any
alternating circuit decomposition of $\nabla$, we are going to assign a sequence
of switches that transform the first realization into the second (this is
described right after the proof of Lemma~\ref{th:switch-sequence}).

\bigskip

A graph $H$, with edges colored by either red or blue, will be called a
\textbf{red-blue graph}. For vertex $v$ let $d_r(v)$ and $d_b(v)$ be the degree
of vertex $v$ in red and blue edges, respectively. This red-blue graph is
\textbf{balanced} if for each $v\in V(H)$ equality $d_r(v)=d_b(v)$ holds.

\medskip

Let $X, Y$ both be realizations (on the same vertex set) of an unconstrained
degree sequence $\bd$ or a bipartite degree sequence $\bD$. Let the symmetric
difference of the edges be
\begin{equation*}
	\nabla:=E(X)\triangle E(Y).
\end{equation*}
Color the edges of $\nabla$ according to which graph they come from: the $E(X)$
edges are colored red and the $E(Y)$ edges are colored blue. Equipped with this
coloring, $\nabla$~is a \BrbG.

\medskip

A \textbf{circuit} in a graph $H$ is a closed trail (so any edge is traversed at
most once). As the graph is simple, a circuit is determined by the sequence of
the vertices $v_0,\ldots,v_t$, where $v_0=v_t$. Note that there can also be
other indices $i<j$ such that $v_i=v_j$. A circuit is called a \textbf{cycle},
if its simple, i.e., for any $i<j$, $\;v_i=v_j$ only if $i=0$ and $j=t$.

A circuit (or, in particular, a cycle) in a \BrbG\ is called
\textbf{alternating}, if the color of its edges alternates. In other words, the
color of the edge from $v_i$ to $v_{i+1}$ differs from the color of the edge
from $v_{i+1}$ to $v_{i+2}$, and also edges $v_0v_1$ and $v_{t-1}v_t$ have
different colors. Consequently, alternating circuits have even length.

The meaning of ``alternating'' slightly varies with context. For a subgraph of
$\nabla$, it is used in the red/blue alternating sense; we may explicitly say
$(X,Y)$-alternating in this case. For realizations, ``alternating'' refers to
the conventional edge/non-edge alternation. In practice, this should not be
confusing, since:
\begin{lemma}\label{lemma:alternation}\hfill
	\vspace{-6pt}
	\begin{itemize}
		\item A circuit $C$ which is alternating in the red-blue graph
			$\nabla$ is alternating between edges and non-edges of
			$X$ (and $Y$ as well).
		\item If $C$ is an alternating circuit in $X$ then $C$ is a
			red-blue alternating circuit with respect to $X$ and
			$Y:=X\triangle C$.
	\end{itemize}
\end{lemma}
\begin{proof}
	The blue edges are missing from $X$, while the red edges are contained
	in $X$.
\end{proof}

The following observations are easy to see.
\begin{lemma}[adapted from~\cite{EKM}, Erdős, Király and Miklós, 2013]\label{th:list}
	\hfill
	\begin{enumerate}[label=\rm(\roman*)]
		\item If $H$ is a \BrbG\ then the edge set can be decomposed
			into alternating circuits.\label{item:listDecomp}
		\item Let $C = v_0, v_1, \ldots, v_{2t}=v_0$ be an alternating
			circuit in a \BrbG\ $H$, in which for some $i<j<2t$,
			$j-i$ is even and $v_i = v_j$. Then the circuit can be
			decomposed into two shorter alternating circuits:
			\begin{equation*}
				C=(v_i,v_{i+1},\ldots,v_{j-1},v_j)\uplus
				(v_j,v_{j+1},\ldots,v_{2t-1},
				v_0,v_1,\ldots,v_{i-1},v_i).
			\end{equation*}
			\label{item:listShorter}
	\end{enumerate}
\end{lemma}
\vspace{-1em}
It is clearly possible that a vertex occurs twice in an alternating circuit
without the possibility to divide it into two, smaller alternating circuits. The
smallest example is a ``bow-tie'' circuit: $v_1,v_2, v_3,v_1,v_4, v_5,v_1$ with
an alternating edge coloring. (The very first and very last occurrences of $v_1$
shows the closing of the alternating circuit.) Recalling our earlier discussion,
these two copies of the vertex $v_1$ form a \textbf{non-chord}.
\begin{definition}\label{def:primitive}
	An alternating circuit is \textbf{primitive}, if it cannot be decomposed
	further in the way described in
	Lemma~\ref{th:list}\ref{item:listShorter}.
\end{definition}
From Lemma 2.1 it follows, that if a vertex appears twice in a primitive
circuit, then the distance of two copies of the same vertex must be odd.
Consequently, a vertex cannot be visited three times by a primitive circuit,
because the three pairwise distances cannot be simultaneously odd. Therefore, in
a bipartite graph, a primitive alternating circuit is an alternating cycle.

The definition of primitive alternating circuits is different from the
definition of elementary alternating circuits introduced in~\cite{EKM} (Erdős,
Király and Miklós, 2013). Without providing the definition here, we mention
that, for example, the red-blue graph obtained from a red $C_5$ and its blue
complement is a primitive alternating circuit, see Figure~\ref{fig:c10}. Still,
it can be decomposed into an alternating $C_4$ ($x_1,x_5,x_7,x_6$) and an
alternating bow-tie, so it is not elementary.

\begin{figure}[ht]
	\centering
	\begin{tikzpicture}[scale=0.30]
		\begin{scope}[shift={(20,10)}]
		\foreach \i in {1,...,10}
		{
			\pgfmathsetmacro\deg{\i*36+180}
			\node[very thick,circle,draw,label={\deg:$x_{\i}$}]
				(x\i) at (\deg:10) {};
		}
		\end{scope}

			\begin{scope}[shift={(-5,10)}]
				\foreach \i in {1,3,5,7,9}
				{
					\pgfmathtruncatemacro\j{mod(\i+2,10)+1};
				\pgfmathsetmacro\deg{\i*72-36}
				\node[very thick,circle,draw,label={\deg:$x_{\i},x_{\j}$}]
					(y\i) at (\deg:7) {};

				\draw [very thick,dashed,gray!50] (x\i) -- (x\j);
			}
		\end{scope}

		\foreach \i in {1,3,5,7,9}
		{
			\pgfmathtruncatemacro\j{\i+1}
			\draw [very thick,dotted,red] (x\i) -- (x\j);
			\pgfmathtruncatemacro\result{mod(\j+6,10)+1};
		\draw [very thick,dotted,red] (y\i) -- (y\result);

		\pgfmathtruncatemacro\k{mod(\j,10)+1};
	\pgfmathtruncatemacro\result{mod(\k+3,10)+1};
\draw [very thick] (x\j) -- (x\k);
\draw [very thick] (y\result) -- (y\k);
		}

		\begin{scope}[shift={(4,0)}]
			\draw [very thick] (-8,-3) -- (-6,-3);
			\node at (-4,-3) {edge};

			\draw [very thick,dotted,red] (-1.5,-3) -- (0.5,-3);
			\node at (5.5,-3) {chord, non-edge};

			\draw [very thick,dashed,gray!50] (11,-3) -- (13,-3);
			\node at (16,-3) {non-chord};

		\end{scope}
	\end{tikzpicture}

	\caption{The cycle of length 5 has an alternating circuit traversing all of its edges and non-edges}\label{fig:c10}
\end{figure}

\begin{lemma}\label{th:C6}
       Let $C$ be an alternating circuit of length 6 in $\nabla$. If the only
       non-chords are loops, then there is at most one vertex which is visited
       more than once by $C$. In other words, at most one of the three main
       diagonals of $C$ is a non-chord.
\end{lemma}
\begin{proof}
       If $v$ is a vertex that is visited at least twice by $C$, it has at least
       four other neighbors that are pairwise distinct from each other and $v$.
       Since $v$ is counted twice in the length of $C$, every vertex of $C$ is
       accounted for, and the claim holds (and $C$ is a bow-tie).
 \end{proof}

We will use Sinclair's multicommodity flow method (Theorem~\ref{th:sinc}) to
bound the mixing time of the switch Markov chain. The multicommodity flow is
given by a set of switch sequences between any two realizations of the degree
sequence. The main idea behind the definition of the flow can be described
roughly as follows:

For each pair of realizations $(X,Y)$ and every possible complete matching of
the $X$-edges with $Y$-edges in $\nabla=X\triangle Y$ at every vertex, assign
a switch sequence from $X$ to $Y$ as follows. The matchings decompose $\nabla$
into alternating circuits (Lemma~\ref{th:decomp} in Section~\ref{sec:step1}).
Refine each of the alternating circuits into primitive alternating circuits in a
canonical way (Section~\ref{sec:order}). By concatenating the switch sequences
given by Algorithm~\ref{alg:sweep} for the primitive alternating circuits, a
switch sequence from $X$ to $Y$ is obtained.

Here we reached a very important  point: Algorithm~\ref{alg:sweep} does not
require an order for processing the primitive circuits; in principle it can be
done arbitrarily. One of the novelties of this paper leading to the unified
proof is constructing a very delicate order of processing the primitive
circuits, which ultimately enables the unification of the proofs. We will return
to this point in Section~\ref{sec:order}. We design the order of processing the
primitive circuits in Section~\ref{sec:design}.

\medskip

The \textsc{Sweep} procedure in Algorithm~\ref{alg:sweep} will be used to
construct the switch sequence between two realizations whose symmetric
difference is a primitive alternating circuit $C=x_1x_2\cdots x_{2\ell}$.
Suppose that some of the alternating primitive circuits from $\nabla$ have
already had their edges and non-edges exchanged compared to $X$ along the switch
sequence; let the current realization be $G$. Let $C$ be the next alternating
circuit in $X$ whose edges and non-edges we exchange using \textsc{Sweep}. We
assume that $x_1x_2\not\in E(G)$.  Most steps of the algorithm involve a single
switch (line~\ref{state:switch}) which removes an edge $x_1x_{2t+2}$, fixes two
edges $x_{2t}x_{2t+1}$, $x_{2t+1}x_{2t+2}$ of the symmetric difference, and
inserts an edge $x_1x_{2t}$ (which has to be a chord). If $x_1x_{2t}$ is
non-chord then a different function, \textsc{Double step} (described in
Algorithm~\ref{alg:switches}), is called twice in succession, to perform two
switches. \textsc{Double step} is designed to avoid the use of $x_1x_{2t}$. When
$x_1=x_{2t}$ and $x_{2t-1}=x_{2t+2}$ or $x_{2t-2}=x_{2t+1}$, however, calling
\textsc{Double step} is not feasible, therefore a special switch avoiding
$x_1x_{2t}$ is performed (lines~\ref{state:switch1}~and~\ref{state:switch2}).
Addition and subtraction of edges naturally means that we add the edge to, or
remove the edge from, the edge set of the first operand.

\begin{algorithm}[H]
	\caption{Sweeping a primitive alternating circuit
	$C=(x_1,x_2,\ldots,x_{2\ell})$ in $G$. The algorithm assumes that
	$x_1x_2\notin E(G)$.}\label{alg:sweep}
	\begin{algorithmic}[1]
		\Procedure{Sweep}{$G,{[x_{1},x_{2},\ldots,x_{2\ell}]}$} $\to
		[Z_{0},Z_{1},Z_{2},\ldots, Z_{\ell-2},\langle Z_{\ell-1}\rangle ]$
		\State $Z_0\gets G$
		\State $q\gets 1$
		\State $\mathit{end}\gets 2$
		\If{$\exists r\in\mathbb{N}$ $x_1=x_{2r}$}
		\State $\mathcal{L}\gets\left\{ 2i\in2\mathbb{N} :\ 4\le 2i\le
		2\ell,\,\, x_1x_{2i}\in E(G)\, \text{ and }\, x_{2i}\neq
		x_{2r+1}  \right\}$
		\Else
		\State $\mathcal{L}\gets\left\{ 2i\in2\mathbb{N} :\ 4\le 2i\le 2\ell,\,\, x_1x_{2i}\in E(G)\right\}$
		\EndIf
		\Statex
		\While{$\mathit{end}<2\ell$}
		\State $\mathit{start}\gets \min\big\{ 2i\in\mathcal{L}\ :\ 2i>\mathit{end} \big\}$
		\State $2t\gets \mathit{start}-2$
		\Statex
		\While{$2t\ge\mathit{end}$}\label{state:innerloop}
		\If{$x_1=x_{2t}$ and $x_{2t+2}= x_{2t-1}$}
		\State $Z_q\gets Z_{q-1}-\{x_{2t}x_{2t+1},x_{2t-2}x_{2t-1}\}+\{x_{2t+1}x_{2t+2},x_{1}x_{2t-2}\}$\label{state:switch1}
		\State $2t\gets 2t-4$
		\ElsIf{$x_1=x_{2t}$ and $x_{2t+1}= x_{2t-2}$}
		\State $Z_q\gets Z_{q-1}-\{x_1x_{2t+2},x_{2t-2}x_{2t-1}\}+\{x_{2t+1}x_{2t+2},x_{2t-1}x_{2t}\}$\label{state:switch2}
		\State $2t\gets 2t-4$
		\ElsIf{$x_1=x_{2t}$ or $x_1x_{2t}\in E(Z_{q-1})$}
		\State $Z_{q}\gets \Call{Double step}{Z_{q-1},x_1,{[x_{2t-2},\ldots,x_{2t+2}]}}$\label{state:doublestep1}
		\State $q\gets q+1$
		\State $Z_{q}\gets \Call{Double step}{Z_{q-1},x_1,{[x_{2t-2},\ldots,x_{2t+2}]}}$\label{state:doublestep2}
		\State $2t\gets 2t-4$
		\Else
		\State $Z_q\gets Z_{q-1}-\big\{x_1x_{2t+2},x_{2t}x_{2t+1}\big\}+\big\{x_1x_{2t},x_{2t+1}x_{2t+2}\}$\label{state:switch}
		\State $2t\gets 2t-2$
		\EndIf
		\State $q\gets q+1$
		\EndWhile
		\Statex
		\State $\mathit{end}\gets \mathit{start}$
		\EndWhile
		\Statex
		\EndProcedure
	\end{algorithmic}
\end{algorithm}

\begin{lemma}\label{th:switch-sequence}
	Suppose that $G$ contains a primitive alternating circuit $C$ of length
	$2\ell$. Let the vertex sequence of $C$ be $(x_1,x_2,\ldots,x_{2\ell})$,
	and suppose that $x_1x_2\notin E(G)$.

	\medskip

	Then Algorithm~\ref{alg:sweep} provides a valid switch sequence
	between $G$ and $G\triangle C$ in the case of unconstrained and
	bipartite graph models. The length of the switch sequence is $\ell-1$ in
	the bipartite case; in the unconstrained case, the length is either
	$\ell-1$ or $\ell-2$.
\end{lemma}
\begin{proof}
	The processing done by Algorithm~\ref{alg:sweep} is governed by two
	nested loops. The outer loop iterates the variable $\mathit{start}$
	through
	\begin{equation*}
		\mathcal{L}=\Big\{ 2i\in2\mathbb{N} :\ 4\le 2i\le 2\ell,\,\,
		x_1x_{2i}\in E(G)\Big\} \cap \Big\{ 2i \in 2\mathbb{N} \mid
	\text{$x_{2i}\neq x_{2r+1}$ whenever $x_{2r}=x_1$}\Big\}
	\end{equation*}
	in increasing order.  The set $\mathcal{L}$ is not empty since
	$x_1x_{2\ell}\in E(G)$ and, if $x_1 = x_{2r}$ for some $r$ then
	$x_{2\ell}\neq x_{2r+1}$. In the first iteration, $\mathit{end}=2$, and
	in the successive iterations $\mathit{end}$ takes the value taken by
	$\mathit{start}$ in the previous iteration.

	\medskip

	The inner loop performs a series of switches that changes the status of
	the edges and non-edges induced by consecutive vertices in the interval
	of vertices between $x_\mathit{start},\ldots,x_\mathit{end}$. As a side
	effect, it also changes chords between $x_1$ and
	$\{x_4,x_6,\ldots,x_{2\ell-4},x_{2\ell-2}\}$, and rarely other
	chords induced by two vertices of $C$.

	\medskip

	We have to check that each time a new graph is obtained from $Z_{q-1}$
	(lines \ref{state:switch1}, \ref{state:switch2},
	\ref{state:doublestep1}, \ref{state:doublestep2}, and~\ref{state:switch}
	in Algorithm~\ref{alg:sweep}), the changes correspond to valid switches.

	\begin{figure}[H]
	\centering
	\begin{tikzpicture}[scale=1.2]
	\def\n{16};
	\def\r{2.5cm};

	\foreach \i in {1,...,11,12}
	{
		\pgfmathsetmacro{\deg}{\i*360/\n-180}
		\node[thick,draw,circle,minimum size=8pt, inner sep=0,label={\deg:$x_{\i}$}] (v\i)
			at (\deg:\r) {};
	}

	\pgfmathsetmacro{\deg}{14*360/\n-180}
	\node[thick,draw,circle,minimum size=8pt, inner
		sep=0,label={\deg:$x_{2\ell-2}$}] (v14)
		at (\deg:\r) {};
	\pgfmathsetmacro{\deg}{15*360/\n-180}
	\node[thick,draw,circle,minimum size=8pt, inner
		sep=0,label={\deg:$x_{2\ell-1}$}] (v15)
		at (\deg:\r) {};
	\pgfmathsetmacro{\deg}{16*360/\n-180}
	\node[thick,draw,circle,minimum size=8pt, inner
		sep=0,label={\deg:$x_{2\ell}$}] (v16)
		at (\deg:\r) {};

	\foreach \i in {1,...,11,14,15,16}
	{
		\pgfmathtruncatemacro{\j}{mod(\i,16)+1}
		\ifthenelse{\isodd{\i}}{\draw[very thick,dotted,red] (v\i) -- (v\j);}
			{\draw[very thick,black] (v\i) -- (v\j);}
	}
	\draw [very thick,decorate,decoration=snake,dotted] (v12) -- (v14);
	\foreach \i in {4,8,10,14} \draw[very thick, dotted,red] (v1) -- (v\i);
	\foreach \i in {6,12}\draw[very thick, black] (v1) -- (v\i);

	\draw[very thick,dotted,red] (15:5cm) -- ++(0.5,0) node[black,anchor=west]
		{chord, non-edge};
	\draw[very thick] ($ 0.5*(15:5cm) + 0.5*(-15:5cm) $) -- ++(0.5,0) node[anchor=west] {edge};
	\draw[very thick,decorate,decoration=snake,dotted] (-15:5cm) --
		++(0.5,0) node[black,anchor=west]
		{not shown};
	\end{tikzpicture}

	\caption{An alternating cycle $C$ in $G$ where every $x_1x_{2t}$ is a
	chord (an edge or a non-edge)}\label{fig:cyclesweep}
\end{figure}

	Let us first check the case when $G$ is bipartite; for an initial state,
	see Figure~\ref{fig:cyclesweep}. As we remarked
	earlier, $C$ must be a cycle: if $x_i=x_j$ and $i\neq j$, then
	Definition~\ref{def:primitive} implies that $i\neq j\pmod {2}$, but
	because $C$ is alternating, $x_i$ and $x_j$ have to be in different
	vertex classes. Therefore lines \ref{state:switch1} and
	\ref{state:switch2} in
	Algorithm~\ref{alg:sweep} are never reached. We will use induction to prove that
	whenever a new iteration of the inner loop (line \ref{state:innerloop})
	starts, we have
	\begin{align}
	\begin{split}
		Z_{q-1}=G&\triangle \{x_ix_{i+1}\ :\ i\in
		[1,\mathit{end})\cup [2t+2,\mathit{start}) \}\triangle \\
		&\triangle \{x_1x_\mathit{end},x_1x_\mathit{start}\}
		\triangle \{x_1x_{2t+2}\}.
	\end{split}\label{eq:switch}
	\end{align}
	If Equation~\eqref{eq:switch} holds when line~\ref{state:innerloop} is
	reached, then $x_1x_{2t}\notin E(Z_{q-1})$, which means that
	line~\ref{state:switch} alone produces the switch sequence. For $q=1$,
	equation~\eqref{eq:switch} holds. Now assume that $q\geq 2$ and
	\eqref{eq:switch} holds. The set $\{x_1,x_{2t},x_{2t+1},x_{2t+2}\}$ is
	guaranteed to be a set of four vertices, and Equation~\eqref{eq:switch}
	guarantees that $x_1x_{2t}\notin E(Z_{q-1})$. It also follows that the
	switch on line~\ref{state:switch} neither creates a multi-edge nor
	deletes an edge which is not present in $Z_{q-1}$. Hence
	\begin{equation*}
		Z_q = Z_{q-1}\triangle \{ x_1x_{2t}, \, x_{2t}x_{2t+1}, \,
		x_{2t+1}x_{2t+2}, \, x_1x_{2t+2}\}
	\end{equation*}
	which implies that \eqref{eq:switch} holds for $Z_q$ (as $2t+2$ becomes
	$2t$). When $2t=\mathit{end}$, the graph assigned to $Z_q$ on
	line~\ref{state:switch} is $G\triangle \{x_ix_{i+1}\ :\ i\in
	[1,\mathit{start})\}\triangle \{x_1x_\mathit{start}\}$. In particular,
	$Z_{\ell-1}=G\triangle C$, which is what we wanted.

	\medskip

	Next, let us verify the algorithm for unconstrained degree sequences.
	Observe the following:
	\begin{enumerate}[label=\textbf{Obs.\ (\arabic*)},ref={Obs.\ (\arabic*)},leftmargin=4em]
		\item If $x_i=x_j$, then either $i=j$ or $i\not\equiv j
			\pmod{2}$, in accordance with
			Lemma~\ref{th:list}\ref{item:listShorter}. In
			particular, $\{x_1,x_3,\ldots,x_{2\ell-1}\}$  and
			$\{x_2,x_4,\ldots,x_{2\ell}\}$ are both sets of size
			$\ell$, and
			\begin{equation*}
				\{x_{2i}x_{2j+1}\ :\ x_{2i}=x_{2j+1}\}
			\end{equation*}
			is a (partial) pairing between the two sets.
			\label{observation1}

		\item Suppose $x_1x_{2i}$ is a chord. (It may or may not be an
			edge in $G$.) If $x_1x_{2i}\in C$, then $\exists
			r\in\mathbb{N}$ such that $x_1=x_{2r}$ and
			$x_{2i}\in\{x_{2r-1},x_{2r+1}\}$.  \label{observation2}
	\end{enumerate}

	\medskip

	Suppose first, that $x_1$ is visited exactly once by $C$. Notice, that
	the arguments of the bipartite case go through seamlessly for the
	unconstrained case because of~\ref{observation1} and~\ref{observation2}.
	The two observations guarantee that~\eqref{eq:switch} holds by
	induction, since $x_1x_{2t}$ does not appear in $\{x_{i}x_{i+1}\ :\ i\in
	[2,2\ell-1] \}$ for any $2\le 2t \le 2\ell$.
	Figure~\ref{fig:cyclesweep} is still a faithful picture, although
	some of the vertices may be identical, edges and non-edges are not
	repeated.

	\medskip

\begin{figure}[H]
	\centering
	\begin{tikzpicture}[scale=1.2]
	\def\n{16};
	\def\r{2.5cm};

	\foreach \i in {1,...,11,12}
	{
		\pgfmathsetmacro{\deg}{\i*360/\n-180}
		\node[thick,draw,circle,minimum size=8pt, inner sep=0,label={\deg:$x_{\i}$}] (v\i)
			at (\deg:\r) {};
	}

	\pgfmathsetmacro{\deg}{14*360/\n-180}
	\node[thick,draw,circle,minimum size=8pt, inner
		sep=0,label={\deg:$x_{2\ell-2}$}] (v14)
		at (\deg:\r) {};
	\pgfmathsetmacro{\deg}{15*360/\n-180}
	\node[thick,draw,circle,minimum size=8pt, inner
		sep=0,label={\deg:$x_{2\ell-1}$}] (v15)
		at (\deg:\r) {};
	\pgfmathsetmacro{\deg}{16*360/\n-180}
	\node[thick,draw,circle,minimum size=8pt, inner
		sep=0,label={\deg:$x_{2\ell}$}] (v16)
		at (\deg:\r) {};

	\foreach \i in {1,...,11,14,15,16}
	{
		\pgfmathtruncatemacro{\j}{mod(\i,16)+1}
		\ifthenelse{\isodd{\i}}{\draw[very thick,dotted,red] (v\i) -- (v\j);}
			{\draw[very thick,black] (v\i) -- (v\j);}
	}
	\draw [very thick,decorate,decoration=snake,dotted] (v12) -- (v14);
	\draw[very thick,dashed,gray] (v1) -- (v8);
	\draw[very thick,dashed,gray] (v7) -- (v14);
	\draw[very thick,dashed,gray] (v6) -- (v9);
	\foreach \i in {4,10,14} \draw[very thick, dotted,red] (v1) -- (v\i);
	\foreach \i in {6,12}\draw[very thick, black] (v1) -- (v\i);

	\draw[very thick,dotted,red] (25:5cm) -- ++(0.5,0) node[black,anchor=west]
		{chord, non-edge};
	\draw[very thick] ($ 0.66*(25:5cm) + 0.34*(-25:5cm) $) -- ++(0.5,0) node[anchor=west] {edge};
	\draw[very thick,dashed,gray] ($ 0.34*(25:5cm) + 0.66*(-25:5cm) $) -- ++(0.5,0) node[black,anchor=west]
		{non-chord};
	\draw[very thick,decorate,decoration=snake,dotted] (-25:5cm) --
		++(0.5,0) node[black,anchor=west]
		{not shown};
	\end{tikzpicture}

	\caption{A primitive circuit with 3 non-chords shown. Each non-chord
	joins two copies of a vertex, i.e., $x_1=x_8$, etc. The end-points of
	non-chords are pairwise disjoint.}\label{fig:sweep}
\end{figure}

	Suppose now, that $x_1$ is visited twice by $C$, that is, $x_1=x_{2r}$
	for some integer $r$, where $2<2r<2\ell$. See Figure~\ref{fig:sweep},
	where $x_1=x_8$. There are two ways the induction argument may break when
	$x_1=x_{2r}$. Since $x_1x_{2r}$ is a non-chord, it cannot participate
	in a switch. Furthermore, it is possible that $x_{1}x_{2t}\in
	\{x_{i}x_{i+1}\,:\,i\in [2,2\ell-1] \}$ for some $2t$.  The only way for
	this to happen is if $x_1=x_{2r}$ and $x_{2t}\in \{x_{2r+1},x_{2r-1}\}$.

	\medskip

	We will verify that Equation~\eqref{eq:switch} holds whenever
	Algorithm~\ref{alg:sweep} reaches line~\ref{state:innerloop}.  Notice,
	that if $2t\in \mathcal{L}$, then by definition $x_{2t}\neq x_{2r+1}$,
	and $x_{2t}\neq x_{2r-1}$, because $x_1x_{2t}=x_{2r-1}x_{2r}\notin E(G)$
	($C$ alternates in $G$). When
	$x_{2t},x_{2t+2}\notin\{x_{2r-1},x_{2r},x_{2r+1}\}$, the arguments of
	the previous case are still applicable, because $x_1x_{2t},x_1x_{2t+2}$
	are chords and
	\begin{equation*}
		x_1x_\mathit{start},x_1x_\mathit{end},x_1x_{2t},x_1x_{2t+2}\notin
		\{x_{i}x_{i+1}\ :\ i\in [2,2\ell-1] \}.
	\end{equation*}

	\medskip

	Let us check that Equation~\eqref{eq:switch} is preserved by the
	iteration of the inner loop that starts with $2t=2r$.  Suppose $2t=2r$,
	$x_{2r-1}\neq x_{2r+2}$ and $x_{2r-2}\neq x_{2r+1}$.  We choose not to
	perform the standard switch on line~\ref{state:switch}, instead, we
	perform a \textsc{Double step} (Algorithm~\ref{alg:switches}) to avoid
	using $x_1x_{2r}$ in a switch. \textsc{Double step} is called on
	line~\ref{state:doublestep1} to apply a switch to $Z_{q-1}$. We have
	$\mathit{end}\le 2r-2$, $\mathit{start}\ge 2r+2$, and by
	induction,~\eqref{eq:switch} holds for $Z_{q-1}$, thus
	$x_1x_{2r-2}\notin E(Z_{q-1})$ and $x_1x_{2r+2}\in E(Z_{q-1})$.

	\begin{algorithm}[H]
		\caption{The \textsc{Double step} avoids inserting or deleting
		$x_1x_{2t}$. The algorithm assumes that at most one of
		$x_{2t-2}x_{2t+1}$, $x_{2t-1}x_{2t+2}$ is a
		non-chord.}\label{alg:switches}
		\begin{algorithmic}
			\Function{Double step}{${Z},x_1,{[x_{2t-2},x_{2t-1},x_{2t},x_{2t+1},x_{2t+2}]}$}
			\If{$x_{2t-2}x_{2t+1}$ is a chord}
			\If{$x_{2t-2}x_{2t+1}\in E({Z})$}
			\State\Return ${Z}+\{x_1x_{2t-2},x_{2t+1}x_{2t+2}\}-\{x_{2t-2}x_{2t+1},x_1x_{2t+2}\}$
			\ElsIf{$x_{2t-2}x_{2t+1}\notin E({Z})$}
			\State\Return ${Z}+\{x_{2t-2}x_{2t+1},x_{2t-1}x_{2t}\}-\{x_{2t-2}x_{2t-1},x_{2t}x_{2t+1}\}$
			\EndIf
			\ElsIf{$x_{2t-1}x_{2t+2}$ is a chord}
			\If{{$x_{2t-1}x_{2t+2}\in E(Z)$}}
			\State\Return ${Z+\{x_{2t+1}x_{2t+2},x_{2t-1}x_{2t}\}-\{x_{2t-1}x_{2t+2},x_{2t}x_{2t+1}\}}$
			\ElsIf{$x_{2t-1}x_{2t+2}\notin E({Z})$}
			\State\Return ${Z+\{x_1x_{2t-2},x_{2t-1}x_{2t+2}\}-\{x_{2t-2}x_{2t-1},x_{1}x_{2t+2}\}}$
			\EndIf
			\EndIf
			\EndFunction
		\end{algorithmic}
	\end{algorithm}

	All in all, if $2t=2r$ and $x_{2r-1}x_{2r+2},x_{2r-2}x_{2r+1}$ are
	chords, then the vertices
	\begin{equation*}
		(x_1,x_{2t-2},x_{2t-1},x_{2t},x_{2t+1},x_{2t+2},x_1)
	\end{equation*}
	in this order trace out a bow-tie. Thus the first call to \textsc{Double
	step} on line~\ref{state:doublestep1} performs a valid switch. The next
	call to \textsc{Double step} on line~\ref{state:doublestep2} likewise
	performs a valid switch, and it restores alternation between edges along
	$(x_1,x_{2t-2},x_{2t-1},x_{2t},x_{2t+1},x_{2t+2},x_1)$; furthermore,
	$Z_{q+1}$ satisfies~\eqref{eq:switch} (for $q+2$).

	\medskip

	Suppose $2t=2r$. If $x_{2r-2}x_{2r+1}$ or $x_{2t-1}x_{2t+2}$ is a
	non-chord, we perform a special switch, either on
	line~\ref{state:switch1} or on line~\ref{state:switch2}. It cannot
	happen that both $x_{2r-2}x_{2r+1}$ and $x_{2t-1}x_{2t+2}$ are
	non-chords, because $x_{2r+1}x_{2r+2}\notin E(G)$ and
	$x_{2r-2}x_{2r-1}\in E(G)$. It is simple to check that in either case, a
	valid switch is performed, and~\eqref{eq:switch} holds in the next
	iteration of the inner loop.

	In any case, we will never have $x_{2t+2}=x_{2r}$ when the algorithm
	reaches line~\ref{state:innerloop}.

	\medskip

	If $x_{2t}=x_{2r+1}$, $x_{2t+2}\neq x_{2r-1}$, and $x_{1}x_{2t}\notin
	E(Z_{q-1})$, then line~\ref{state:switch} performs a valid switch on
	$Z_q$, and~\eqref{eq:switch} holds in the next iteration. Suppose, that
	$x_{2t}=x_{2r+1}$, $x_{2t+2}\neq x_{2r-1}$, and $x_1x_{2t}\in E(G)$. In
	this iteration of the inner loop, two calls to \textsc{Double step} are
	made on lines~\ref{state:doublestep1} and~\ref{state:doublestep2}. Note,
	that $2t\notin\mathcal{L}$, which means that
	$\mathit{end}<2t<\mathit{start}$. By Equation~\eqref{eq:switch}, we have
	$x_1x_{2t+2}\in E(Z_{q-1})$. By the assumption, $x_{2r}x_{2r+1}\notin
	\{x_{i}x_{i+1}\ :\ i\in [1,\mathit{end})\cup [2t+2,\mathit{start}) \}$,
	which implies that $x_{2r-1}x_{2r}$ is also not contained in this set,
	so $x_{2r-1}x_{2r}\notin E(Z_{q-1})$. Therefore, even if
	$x_{2t-2}=x_{2r-1}$, the vertices
	$(x_1,x_{2t-2},x_{2t-1},x_{2t},x_{2t+1},x_{2t+2},x_1)$ form an
	alternating cycle of 6 vertices if $x_{2t-1}\neq x_{2t+2}$ and $x_{2t-2}\neq
	x_{2t+1}$. The alternating circuit
	$(x_1,x_{2t-2},x_{2t-1},x_{2t},x_{2t+1},x_{2t+2},x_1)$ is a bow-tie if
	$x_{2t-1}=x_{2t+2}$ or $x_{2t-2}=x_{2t+1}$. By Lemma~\ref{th:C6}, we
	cannot have $x_{2t-1}=x_{2t+2}$ and $x_{2t-2}=x_{2t+1}$ simultaneouly,
	thus the two calls to \textsc{Double step} perform valid switches, and
	Equation~\eqref{eq:switch} holds at the start of the next iteration.

	\medskip

	The case of $x_{2t}=x_{2r-1}$, $x_{2t+2}\neq x_{2r+1}$, and $x_1x_{2t}\notin
	E(Z_{q-1})$ is similar to the previous case. We can deduce that
	\begin{equation*}
		x_{2r-1}x_{2r},x_{2r}x_{2r+1}\in \{x_{i}x_{i+1}\ :\ i\in
		[1,\mathit{end})\cup [2t+2,\mathit{start}) \},
	\end{equation*}
	i.e., we have $x_1x_{2t-2}\notin E(Z_{q-1})$, even if
	$x_{2t-2}=x_{2r+1}$.

	\medskip

	Lastly, suppose that $x_{2t+2}\in \{x_{2r-1},x_{2r+1}\}$ and
	$x_{2t}\notin \{x_{2r-1},x_{2r},x_{2r+1}\}$. From the analysis of the
	previous cases it follows that $x_1x_{2t+2}$ is a chord and
	$x_1x_{2t+2}\in E(Z_{q-1})$. We get $x_1x_{2t}\notin E(Z_{q-1})$ from
	Equation~\eqref{eq:switch}. As before, line~\ref{state:switch} performs
	a valid switch and~\eqref{eq:switch} holds at the start of the next
	iteration.

	\medskip

	We have thus shown that Algorithm~\ref{alg:sweep}
	produces valid switches. If line~\ref{state:switch1}
	or~\ref{state:switch2} is reached during \textsc{Sweep}, then
	$Z_{\ell-2}=G\triangle C$, otherwise $Z_{\ell-1}=G\triangle C$.
	This concludes the proof of Lemma~\ref{th:switch-sequence}.
\end{proof}

We will demonstrate the algorithm on Figure~\ref{fig:sweep}. In the first
iteration of the outer loop, $\mathit{start}$ takes index $12$ as its
initial value.  We call $x_1x_{12}$ the \textbf{start-chord} and $x_1x_2$ the
\textbf{end-chord}. The algorithm \textbf{sweeps} the alternating chords along
the circuit between $x_{2}$ and $x_{12}$, and vertex $x_1$ will be the
\textbf{cornerstone} of this procedure.

\medskip

The inner loop works from the start-chord $x_1x_{12}$ (edge) towards the
end-chord $x_1x_2$ (non-edge). The first value taken by $2t$ is $10$. Since
$x_{1}x_{10}$ is a chord, $Z_1$ is obtained by switching along
$x_1,x_{12},x_{11},x_{10}$. In the next step, $2t=8$. However, $x_1x_{8}$ is a
non-chord, therefore \textsc{Sweep} branches into line~\ref{state:switch2}
instead of line~\ref{state:switch}. $Z_2$ is obtained by switching along
$x_1,x_{10},x_9=x_6,x_7,x_8=x_1$. Subsequently $Z_3$ is obtained by switching
along $x_1,x_6,x_5,x_4$; note that $x_8x_9=x_1x_6\notin E(Z_3)$. The last
iteration of the inner loop (for the first iteration of the outer loop) switches
along $x_1,x_4,x_3,x_2$ and produces $Z_4$. Notice, that all of the chords on
the circuit from $x_2$ to $x_{10}$ changed their status and $x_1x_{12}$ is no
longer an edge in $Z_4$, but the rest of the chords (except $x_1x_6$ and
$x_1x_{2\ell-2}$) have the same status in $Z_4$ as they had in $G$. Furthermore,
the edges and non-edges of the circuit have all been exchanged on the segment
from $x_1$ to $x_{12}$.

\medskip

For the second iteration of the outer loop, $\mathit{end}=12$ and
$\mathit{start}$ is assigned a new value too. Eventually, $\mathit{start}=2\ell$
is set, which marks the last iteration of the outer loop. Until $t=2\ell-2$ (and
$q<\ell-3$), the switches are produced by line~\ref{state:switch}. Note,
that $x_1x_{2\ell-4}\notin E(Z_{\ell-4})$, even when $2\ell-4\in\mathcal{L}$.
However, $x_1x_{2\ell-2}\in E(Z_{\ell-4})$, thus \textsc{Double step} is called
on lines~\ref{state:doublestep1}. The next call to \textsc{Double step} on
line~\ref{state:doublestep2} returns $G\triangle C$ and finishes the switch
sequence.

\medskip

The demonstration shows that some chords that are not necessarily traversed by
$C$ change from being an edge to a non-edge and vice versa during this
procedure. However, there are strict patterns that these irregularities must
abide. We will collect in $R$ the set of pairs of vertices $xy$ of $Z_q$ whose
edge vs.\ non-edge status is only temporarily modified by the switch sequence
from $G$ to $G\triangle C$. In other words, $R$ contains those chords $xy$ of $Z_q$,
that will be changed by a switch before the sequence reaches $G\triangle
C$ and $xy\in E(G)\Leftrightarrow xy\notin E(Z_q)$.

\begin{lemma}\label{lemma:primitivecircuit}
	Suppose $Z_q$ is an intermediate realization produced by
	Algorithm~\ref{alg:sweep} on the switch sequence between $G$ and
	$G\triangle C$, where $C=(x_1,x_2,\ldots,x_{2\ell})$ is an alternating
	primitive circuit. Let
	\begin{equation}\label{eq:Rdef}
		R:=\left((Z_q\triangle G)\setminus {E(C)}\right)\cup Q,
	\end{equation}
	where $Q:=\emptyset$, except if $Z_q$ is the return value of a call to
	{\normalfont\textsc{Double step}} on line~\ref{state:doublestep1}, in
	which case
	\begin{align*}
		Q:=\left\{
			\begin{array}{ll}
				\{x_{2t+1}x_{2t+2}\}, & \text{ if } x_{2t-2}\neq
				x_{2t+1}\text{ and } x_{2t-2}x_{2t+1}\notin
				E(Z_{q-1}), \\
				\{x_{2t-2}x_{2t-1}\}, & \text{ if } x_{2t-2}=
				x_{2t+1}\text{ and } {x_{2t-1}x_{2t+2}\in
				E(Z_{q-1})}, \\
				\emptyset, &\text{ otherwise.}
			\end{array}
		\right.
	\end{align*}
	The following statements hold at the moment when $Z_q$ is assigned a
	value in {\normalfont\textsc{Sweep}}.
	\begin{enumerate}[label=(\alph*)]
		\item $R$ is a set of chords with both end-vertices in $\nabla$,
		\item $R=\emptyset$, if $Z_q\in \{G,G\triangle C\}$,\label{item:Rempty}
		\item $R\subseteq\{x_1x_\mathit{start},x_1x_\mathit{end},x_1x_{2t}\}$,
			if $Z_q$ is produced on
			line~\ref{state:switch},\label{item:Rswitch}
		\item $R\subseteq\{x_1x_\mathit{start},x_1x_\mathit{end},x_1x_{2t-2}\}$,
			if $Z_q$ is produced on line~\ref{state:switch1},
			line~\ref{state:switch2}, or
			{line~\ref{state:doublestep2}} ({the second} call to
			{\normalfont\textsc{Double step}}),\label{item:Rdoublestepsecondgraph}
		\item $R\subseteq \{x_1x_\mathit{start},x_1x_\mathit{end}\}\cup
			H$, where $H$ is defined in Equation~\eqref{eq:ZqH}, if
			$Z_q$ is the return value of a call to
			{\normalfont\textsc{Double step}} on
			{line~\ref{state:doublestep1}}.\label{item:Rdoublestep}
		\item ${(Z_q\triangle R)\triangle G}$ is a set of at most $2$
			edge-disjoint subtrails of $C$, starting and
			ending at endpoints of chords in $R$.\label{item:Rtrails}
		\item The edges in $R$ cover at most 5 vertices besides $x_1$.
			In the bipartite case, $R$ covers at most 3 vertices
			other than $x_1$.\label{item:Rsize}
	\end{enumerate}
\end{lemma}
\begin{proof}
	Most statements follow from the proof of Lemma~\ref{th:switch-sequence}.
	To verify~\ref{item:Rdoublestep} and~\ref{item:Rtrails} when $Z_q$ is
	assigned a value on line~\ref{state:doublestep1} in
	Algorithm~\ref{alg:sweep}, we shall provide a formula for the $Z_q$
	returned on line~\ref{state:doublestep1}. Let us define:
	\begin{equation*}
	F:=\left\{x_ix_{i+1} : i\in\left\{
	        \begin{array}{ll}
			{[1,\mathit{end})}\cup [2t+1,\mathit{start}),
			 & \text{if }x_{2t-2}\neq x_{2t+1}\wedge x_{2t-2}x_{2t+1}\in
			 E(Z_{q-1}) \\
			{[1,\mathit{end})}\cup [2t-2,\mathit{start}),
			 & \text{if }x_{2t-2}\neq x_{2t+1}\wedge x_{2t-2}x_{2t+1}\notin
			 E(Z_{q-1}) \\
			{[1,\mathit{end})}\cup [2t-1,\mathit{start}),
			 & \text{if }x_{2t-2}= x_{2t+1}\wedge x_{2t+2}x_{2t-1}\in
			 E(Z_{q-1}) \\
			{[1,\mathit{end})}\cup [2t+2,\mathit{start}),
			 & \text{if }x_{2t-2}= x_{2t+1}\wedge x_{2t+2}x_{2t-1}\notin
			 E(Z_{q-1})
	        \end{array}
	\right. \right\}
	\end{equation*}
	Moreover, set:
	\begin{equation}
	        H:=
	        \left\{
	                \begin{array}{ll}
	                        x_1x_{2t-2},x_{2t-2}x_{2t+1} & \text{if }x_{2t-2}\neq
				x_{2t+1}\wedge x_{2t-2}x_{2t+1}\in E(Z_{q-1})\\
	                        x_1x_{2t+2},x_{2t-2}x_{2t+1},x_{2t+1}x_{2t+2} & \text{if }
				x_{2t-2}\neq x_{2t+1}\wedge x_{2t-2}x_{2t+1}\notin
				E(Z_{q-1})\\
	                        x_1x_{2t+2},x_{2t+2}x_{2t-1} & \text{if }x_{2t-2}=
				x_{2t+1}\wedge x_{2t+2}x_{2t-1}\in E(Z_{q-1})\\
	                        x_1x_{2t-2},x_{2t+2}x_{2t-1},x_{2t-2}x_{2t-1} & \text{if }
				x_{2t-2}= x_{2t+1}\wedge x_{2t+2}x_{2t-1}\notin
				E(Z_{q-1})
	                \end{array}
		\right\}\label{eq:ZqH}
	\end{equation}
	The following equation holds when \textsc{Double step}
	is called on line~\ref{state:doublestep1}:
	\begin{equation}
		Z_{q}=G\triangle F \triangle
		\{x_1x_\mathit{end},x_1x_\mathit{start}\}\triangle H,
		\label{eq:doublestep1}
	\end{equation}
	which verifies~\ref{item:Rdoublestep}. If $Q$ is non-empty, the single
	chord in $Q$ is not involved in the \textsc{Double step} which created
	$Z_q$, therefore $Q\cap E(Z_q\Delta G)=\emptyset$, which completely
	verifies~\ref{item:Rtrails}.
\end{proof}

Let $X$ and $Y$ be two realizations. Assume that we can decompose the symmetric
difference $\nabla=X\triangle Y$ into $k$ primitive circuits. Let \textsc{Sweep}
process primitive circuits of the decomposition one by one. The concatenation of
the switch sequences returned by \textsc{Sweep} is a switch sequence from $X$ to
$Y$ with at most $\frac{|\nabla|}{2} - k$ switch operations. The process only
changes the statuses of chords induced by vertices of the circuits (which
includes the edges of the circuit).

\section{Sinclair's multicommodity flow method}\label{sec:sinclair}

For \UC\ sequences  we define our Markov chain $(\G_{\bd}, P_{\bd})$ as follows:
in the Markov graph $\G_{\bd} (\G(\bd), \mathbb{E}_{\bd})$ the pair $(G,G')$ is
an edge if these two realizations differ in exactly one switch. To make a move,
choose an unordered pair $F$ of unordered pairs of distinct vertices, uniformly
at random from $G$, say $F=\{(x,y), (z,w)\}$ and choose a perfect matching $F'$
from the other two perfect matchings on the same four vertices. If $F\subseteq
E(G)$ and $F'\cap E(G)=\emptyset$, then perform the switch (so $E(G')=(E(G)\cup
F')\setminus F$). Assuming that $P(G,G')\ne 0$ and $G\ne G'$, we have
\begin{equation}\label{eq:proa}
	\Prob(G\rightarrow G')=P(G,G'):=\frac{1}{2\binom{n}{2}\binom{n-2}{2}}.
\end{equation}
Equation~\eqref{eq:proa} immediately gives that the Markov chain is symmetric.
Notice, that if $(F,F')$ corresponds to a feasible switch, then $(F',F)$ does
not, therefore $P(G,G)\ge \frac{1}{2}$ for any realization $G$. A Markov chain
possessing this property of staying in the current state with probability at
least $\frac12$ is called \textbf{lazy}. Laziness implies that the eigenvalues
of the transition matrix of the Markov chain are non-negative, and it also
implies that the chain is aperiodic.

\medskip

For bipartite degree sequences  we define our Markov chain $(\G_{\bD}, P_{\bD})$
as follows: in the Markov graph $\G_{\bD} (\mathbf{V}_{\bD}, \mathbb{E}_{\bD})$
the pair  $(G,G')$ is an edge, if these two realizations differ in exactly one
switch. The transition matrix $P$ is defined as follows: we choose uniformly an
unordered pair of distinct vertices from $U$, and an unordered pair of distinct
vertices from $V$, then uniformly randomly choose one of the two matchings
between these two pairs.  If it preserves the degree sequence, we remove the
chosen matching, and add the other. The switch moving from $G$ to $G'$ is
unique, therefore the probability of this transformation (the \emph{jumping
probability} from $G$ to $G'\ne G$) is:
\begin{equation}\label{eq:prob}
	\Prob(G\rightarrow G')= P(G,G') :=
	\frac{1}{2\binom{n_1}{2}\binom{n_2}{2}}.
\end{equation}
The transition probabilities are time- and edge-independent, and symmetric.
The probability of staying in the same state is at least $\frac12$.

\bigskip

To begin with, we recall some definitions and notations from the literature.
Since the uniform distribution is the desired stationary distribution of the
switch Markov chains, we conveniently present cited theorems for this specific
case only. Let $P^t$ denote the $t$\textsuperscript{th} power of the transition
probability matrix and let $N:=|V(\G)|$ be the size of the state space of the
Markov chain.  For any element of the state space $X\in V(\G)$, define
\begin{displaymath}
	\Delta_X(t):=\frac 12\sum_{Y\in V(\G)} \left| P^t(X,Y)- 1/N \right|,
\end{displaymath}
We define the \emph{mixing time} of a Markov chain $\mathcal{M}$ as
\begin{displaymath}
	\tau_{\varepsilon}(\mathcal{M}):=\max_{X\in V(\G)}\min_{t}
	\big\{\Delta_X(t')\le {\varepsilon}\text{\ for all $t'\ge t$}\big\}.
\end{displaymath}
The Markov chain is said to be \emph{rapidly mixing} if and only if
\begin{displaymath}
	\tau_{\varepsilon}(\mathcal{M})\le 
	\mathcal{O}\Big (\poly\!\big(\!\log( N/{\varepsilon})\big)\Big)
\end{displaymath}
In this case the switch Markov chain method provides a \emph{fully polynomial
almost uniform sampler} (FPAUS) of the realizations of the given degree
sequences. Using a different Markov chain, which changes at most 2 edges per step (but
which may not preserve the degree sequence), Jerrum and Sinclair have shown that
realizations of $P$-stable degree sequences have a fully polynomial almost
uniform sampler~\cite{JS90} (Jerrum and Sinclair, 1990).

\medskip

Consider the different eigenvalues of the transition matrix of $\mathcal{M}$ in
non-increasing order: if the Markov chain is lazy, we have
\begin{displaymath}
	1={\lambda}_1>{\lambda}_2\ge \dots \ge{\lambda}_N \ge 0.
\end{displaymath}
The \emph{relaxation time} ${\tau}_\mathrm{rel}(\mathcal{M})$ is defined as
\begin{displaymath}
	{\tau}_\mathrm{rel}(\mathcal{M})=\frac 1 {1-{\lambda}^*}
\end{displaymath}
where $\lambda^*$ is the {\em second largest eigenvalue modulus}. So
$\tau_\mathrm{rel}(\mathcal{M})={(1-\lambda_2)}^{-1}$ for lazy chains. The
following result was proved implicitly by Diaconis and Strook in 1991, and
explicitly stated by Sinclair~\cite[Proposition 1]{S92}:
\begin{theorem}[Sinclair]\label{tm:DS}
	\qquad $\displaystyle \tau_{\varepsilon}(\mathcal{M})\le
	\tau_\mathrm{rel}(\mathcal{M})\cdot \log(N/{\varepsilon})$. \qed{}
\end{theorem}
Note, that the cardinality of the state space of the switch Markov chain
trivially satisfies $\log|V(\G)|\le n^2$. By Theorem~\ref{tm:DS}, it is
sufficient to find a $\poly(\log N)$ upper bound on
${\tau}_\mathrm{rel}(\mathcal{M})$ to show that the Markov chain is rapidly
mixing. Rapid convergence to the uniform distribution is highly desirable in any
practical application.

\bigskip

There are different methods to prove rapid convergence, here we use {--}
similarly to~\cite{KTV97} (Kannan, Tetali and Vempala, 1997) {--} Sinclair's
\emph{multicommodity flow method}~\cite[Theorem~5']{S92} (Sinclair, 1992).

\begin{theorem}[Sinclair]\label{th:sinc}
	Let $\mathbb{H}$ be a graph whose vertices represent the possible states
	of a time reversible finite state Markov chain $\mathcal{M}$, and where
	$(U,V) \in E(\mathbb{H})$ if and only if the transition probabilities of
	$\mathcal{M}$ satisfy $P(U,V) P(V,U) \ne 0$. For all $X\ne Y \in
	V(\mathbb{H})$ let $\Gamma_{X,Y}$ be a set of paths in $\mathbb{H}$
	connecting $X$ and $Y$ and let $\pi_{\raisebox{-1pt}{\scriptsize X,Y}}$
	be a probability distribution  on $\Gamma_{X,Y}$. Furthermore let
	\begin{displaymath}
		\Gamma := \bigcup_{ X\ne Y \in V(\mathbb{H})} \Gamma_{X,Y}.
	\end{displaymath}
	We also assume that there is a stationary distribution $\pi$ on the
	vertices $V(\mathbb{H})$. We define the capacity of an edge $e=(W,Z)$ as
	\begin{displaymath}
		Q(e) := \pi(W) P(W,Z)
	\end{displaymath}
	and we denote the length of a path $\gamma$ by $|\gamma|$. Finally let
	\begin{equation}\label{eq:multiflow}
		\kappa_{\Gamma} := \max_{e \in E(\mathbb{H})}
		\frac{1}{Q(e)}
		\sum_{\genfrac{}{}{0pt}{}{\scriptstyle  X,Y \in
				V(\mathbb{H})}{\scriptstyle  \gamma \in \Gamma_{X,Y}\ : \ e \in
		\gamma}} \pi(X) \pi(Y)\pi_{\raisebox{-1pt}{\scriptsize X,Y}}(\gamma) |\gamma|.
	\end{equation}
	Then
	\begin{equation}\label{eq:relax1}
		\tau_{\mathrm{rel}}(\mathcal{M})  \le \kappa_\Gamma
	\end{equation}
	holds. \qed{}
\end{theorem}

\bigskip

We are going to apply Theorem~\ref{th:sinc} for $(\G, P)$, which is either the
unconstrained $(\G, P)=(\G(\bd), P_\bd)$ or the bipartite $(\G, P)=(\G(\bD),
P_\bD)$ switch Markov chain. Using the notation $N:=|V(\G)|$, the (uniform)
stationary distribution is given by $\pi(X) = N^{-1}$ for all $X\in V(\G)$.  So
if we can design a multicommodity flow which is composed of paths shorter than
an appropriate $\poly^\ell(n)$ function, then inequality~\eqref{eq:relax1}
becomes 
\begin{equation}\label{eq:multiflow_simple}
	\tau_{\mathrm{rel}}(\G,P) \leq \frac{\poly^\ell(n)}{N}
	\left (\max_{e \in E(\mathbb{H})}\frac{1}{P(e)}\cdot
	\sum_{\genfrac{}{}{0pt}{}{\scriptstyle X,Y \in V(\mathbb{H})}{\scriptstyle \gamma \in \Gamma_{X,Y}\ : \ e \in     \gamma}}\pi_{\raisebox{-1pt}{\scriptsize X,Y}}(\gamma)\right ).
\end{equation}
If $Z\in e$, then
\begin{equation}\label{eq:multiflow_simplex1}
	\sum_{\genfrac{}{}{0pt}{}{\scriptstyle X,Y \in     V(\mathbb{H})}{\scriptstyle
	\gamma \in \Gamma_{X,Y}\ : \ e \in     \gamma}}\pi_{\raisebox{-1pt}{\scriptsize X,Y}}(\gamma) \leq
	\sum_{\genfrac{}{}{0pt}{}{\scriptstyle X,Y \in     V(\mathbb{H})}{\scriptstyle
	\gamma \in \Gamma_{X,Y}\ : \ Z \in     \gamma}}\pi_{\raisebox{-1pt}{\scriptsize X,Y}}(\gamma),
\end{equation}
so we have
\begin{equation}\label{eq:multiflow_simple1.5}
	\tau_{\mathrm{rel}}(\G,P) \leq 
	\left(\max_{e\in E(\mathbb{H})}\frac{1}{P(e)}\right)\cdot
	\frac{\poly^\ell(n)}{N}\cdot
	\left (\max_{Z \in
		V(\mathbb{H})}
		\sum_{\genfrac{}{}{0pt}{}{\scriptstyle X,Y \in     V(\mathbb{H})}{\scriptstyle
	\gamma \in \Gamma_{X,Y}\ : \ Z \in     \gamma}}\pi_{\raisebox{-1pt}{\scriptsize X,Y}}(\gamma)\right ).
\end{equation}

We make one more assumption. Namely, that for each pair of realizations $X,Y\in
V(\mbb G)$, there is a non-empty finite set $S_{X,Y}$ (which draws its elements
from a pool of symbols), and for each $s\in S_{X,Y}$ there is a path
$\varUpsilon(X,Y,s)$ from $X$ to $Y$. Let
\begin{equation}\label{eq:param}
	\Gamma_{X,Y} := \{\varUpsilon(X,Y,s) : s\in S_{X,Y}\}.
\end{equation}
It can happen that $\varUpsilon(X,Y,s) = \varUpsilon(X,Y,s')$ for $s\ne s'$,
so we consider $\Gamma_{X,Y}$ a ``multiset''. For any ${\gamma}\in \Gamma_{X,Y}$, we
have
\begin{displaymath}
	\pi_{X,Y}(\gamma)=\frac{\left | \big \{s\in S_{X,Y}: \gamma=\varUpsilon(X,Y,s) \big \} \right |
	}{\left |S_{X,Y}\right |}.
\end{displaymath}

\bigskip

Putting together the observations and simplifications above we obtain a slightly
weaker version of Theorem~\ref{th:sinc}. Theorem~\ref{th:simplified_sinclair} is
simpler to use in our applications.

\begin{theorem}\label{th:simplified_sinclair}
Let us assume that there exists two polynomials
$\poly^{\ell}_{\mathcal{D}},\poly^{\rho}_{\mathcal{D}}\in
\mathbb{R}[x]$ (which only depend on the degree sequence class $\mathcal{D}$),
a non-empty finite set $S_{X,Y}$ for each $X\ne Y\in V(\mbb G)$, and a path
$\varUpsilon(X,Y,s)$ from $X$ to $Y$ in $\mbb G$ for each $s\in S_{X,Y}$. If
\begin{itemize}
	\item each path $\varUpsilon(X,Y,s)$ is shorter than
		$\poly^{\ell}_{\mathcal{D}}(n)$, and
	\item for each $Z\in V(\mbb G)$
		\begin{equation}\label{eq:multiflow_simple2}
			\sum_{\genfrac{}{}{0pt}{}{\scriptstyle X,Y \in     V(\mathbb{G})}{ }}
			\frac{\left | \big \{s\in S_{X,Y}: Z\in \varUpsilon(X,Y,s)\big \} \right | } {|S_{X,Y}|}
			\\\le \poly^{\rho}_{\mathcal{D}}(n)\cdot |V(\G)|,
		\end{equation}
\end{itemize}
then the Markov chain $(\mbb G, P)$ is rapidly mixing. Specifically, the mixing
time is at most
\begin{equation}
\tau_{\varepsilon}(\G,P)\le{\left(\min_{P(X,Y)\neq 0}P(X,Y)\right)}^{-1}
	\cdot\poly^{\ell}_{\mathcal{D}}(n)\cdot\poly^{\rho}_{\mathcal{D}}(n)\cdot
	\left(n^2-\log\varepsilon\right)
\end{equation}
In the simple and bipartite cases, each transition $(X,Y)\in E(\G)$ 
satisfies $P(X,Y) \ge n^{-4}$.
\end{theorem}

\section{Multicommodity flow - general considerations}\label{sec:order}

Typically, a successful application of Sinclair's method requires decomposing
$\nabla=X\triangle Y$ into alternating circuits in many ways, each decomposition
yielding a different path in $\G$. The decompositions are
parameterized by $S_{X,Y}$ (see~\eqref{eq:param}). This parametrization
(described in detail in Lemma~\ref{th:decomp}) and its application to \ref{step1}
was introduced in~\cite{KTV97} (Kannan, Tetali and Vempala, 1997).

\medskip

Let $X$ and $Y$ be two realizations of the same (unconstrained or bipartite)
degree sequence; they both belong to $\G$.  For each $s\in S_{X,Y}$ we will
construct a path $\varUpsilon(X,Y,s)$ in $\G$ which connects $X$ to $Y$. A high
level description of the construction of these paths follows.
\begin{enumerate}[label=\textbf{(Step \arabic*)},ref=\mbox{(Step~\arabic*)},leftmargin=6em]
	\item Guided by $s\in S_{X,Y}$, decompose the symmetric difference
		$\nabla = E(X) \triangle E(Y)$ into $X,Y$-alternating circuits:
		$W^s_1,W^s_2\dots, {W^s_{p_s}}$.\label{step1}
	\item Decompose every alternating circuit $W^s_k$ into primitive
		alternating circuits $C^k_1,C^k_2\dots, C^k_{\ell_k}$ (the
		parameter $s$ is omitted from the labels).\label{step2}
	\item Process the primitive circuits iteratively via
		Algorithm~\ref{alg:sweep}. The returned switch sequences are
		concatenated. The resulting path in $\G$ is labeled
		$\varUpsilon(X,Y,s)$.\label{step3}
\end{enumerate}

\noindent We will describe how to perform~\ref{step1} in Section~\ref{sec:step1}
and how to perform~\ref{step2} in Section~\ref{sec:reconst}.  Now suppose that
the decomposition processes in~\ref{step1} and~\ref{step2} are already complete.

\begin{definition}\label{def:milestone}
For any $1\le k\le p_s$ and $1\le r\le \ell_k+1$ let
\begin{equation}\label{eq:milestone}
	G^k_r=X\triangle \left(\cup_{i=1}^{k-1} W^s_i\right)\triangle
	\left(\cup_{j=1}^{r-1}C^k_j\right).
\end{equation}
Each graph $G^k_r$ is a called a \textbf{milestone} of the path $\varUpsilon(X,Y,s)$.
\end{definition}

Clearly, $X=G^1_1$ and $Y=G^{p_s}_{\ell_{p_s}+1}$. Also,
$G^k_{\ell_k+1}=G^{k+1}_1$ for $1\le k< p_s$. We are now equipped to provide a
more technical description of~\ref{step3}.

\textsc{Sweep} is defined in Algorithm~\ref{alg:sweep}. It takes a realization
$G$ and a primitive circuit $C$ which is alternating in $G$ as its input.
\textsc{Sweep} returns a sequences of graphs in which successive elements are
joined by a switch.
\begin{definition}\label{def:path}
	To any $s\in S_{X,Y}$ we associate a
	path $\varUpsilon(X,Y,s)$ in the Markov graph $\G$ defined as the following
	sequence:
	\begin{equation}\label{eq:path}
		\varUpsilon(X,Y,s):=
		\left(\left(\textsc{Sweep}(G^k_r,C^k_r)\right)_{r=1}^{\ell_{k}+1}
		\right)_{k=1}^{p_s}.
	\end{equation}
\end{definition}
\begin{lemma}\label{lemma:pathLength}
	The length of $\varUpsilon(X,Y,s)$ is at most $\frac12|E(X)\triangle
	E(Y)|$.
\end{lemma}
\begin{proof}
	Follows from Lemma~\ref{th:switch-sequence}.
\end{proof}

The most sensitive part of the construction is~\ref{step2}. We will need to ensure
that the following property holds:

\vspace{1em}

\begin{center}
\noindent
\fbox{
	\begin{minipage}{\dimexpr\textwidth-4em}

	\mbox{\ }

	\textbf{(Reconstructability)}
	\customlabel{reconstructability}{\textbf{(Reconstructability)}}
	
	\vspace{4pt}
	
	\hfill\begin{minipage}{\dimexpr\textwidth-3em}
	Let $Z\in \G$ denote an arbitrary vertex along a path $\varUpsilon(X,Y,s)$.
	To apply Sinclair's method we will need that $s\in S_{X,Y}$ can be
	reconstructed from an element of $S_{\nabla\cap E(Z'),\nabla\setminus
	E(Z')}$, $\nabla$, and another small parameter set $B$; here $Z'$ is a
	slight perturbation of $Z$ described by $B$, such that $|E(Z\triangle
	Z')|$ is at most a small constant.

	\mbox{\ }
	\end{minipage}
	\end{minipage}
}
\end{center}

\vspace{1em}

In case of unconstrained degree sequences, Cooper, Dyer and Greenhill
(\cite{CDG07}, 2007) and Greenhill and Sfragara (2018~\cite{G18}) decompose
$W^s_k$ into ``simple'' circuits which have the following property: in each
``simple'' circuit $C$ there is one predefined vertex (actually, the smallest
vertex in a predefined vertex order), which occurs at most twice in $C$. This
made the reconstruction above relatively simple, but made processing such
``simple'' circuits relatively complicated.

As mentioned earlier, primitive circuits are cycles in bipartite graphs. For
bipartite degree sequences such a cycle decomposition is available, which is
provided by the $T$-operator defined in Section~5.2 of~\cite{MES} (Miklós, Erdős
and Soukup, 2013). To adapt this method to the unconstrained degree sequences we
cannot expect to be able to decompose $W^s_k$ into alternating cycles (recall
the bow tie or Figures~\ref{fig:c10} and~\ref{fig:sweep}). In
Section~\ref{sec:design} we generalize the $T$-operator to simple graphs. For
bipartite graphs, the generalized and the original $T$-operator in~\cite{MES}
produce the same decomposition of alternating cycles. The new proof described in
Section~\ref{sec:T-operator} is simpler than that of~\cite{MES} because it is
described on a higher level of abstraction.

\subsection{\ref{step1} {-} parameterizing the circuit decomposition}\label{sec:step1}
Now we describe the parametrization process which was originally introduced by
Kannan, Tetali and Vempala~\cite{KTV97}.  Let $W:=[n]$ in the unconstrained case
and let $W:=([n_1],[n_1+1,n_1+n_2])$ in the bipartite case. Let $K:=(W,F\cup F')$ be a simple
graph where $F\cap F'=\emptyset$ and assume that for each vertex $w \in W$ the
$F$-degree and the $F'$-degree of $w$ are the same: $d_F(w)=d_{F'}(w)$ for all
$w\in W$. An \emph{alternating circuit decomposition} of $(F, F')$ is a circuit
decomposition such that successive edges come alternately from $F$ and $F'$. By
definition, that means that each circuit is of even length. To be more verbose,
we may say that the circuit is $F,F'$-alternating.

The set of all edges in $F$ (in $F'$) which are incident to a vertex $w$ is
denoted by $F(w)$  (by $F'(w)$, respectively).  If $\mathfrak{A}$ and
$\mathfrak{B}$ are sets, denote by $[\mathfrak{A},\mathfrak{B}]$ the complete
bipartite graph with classes $\mathfrak{A}$ and $\mathfrak{B}$.
\begin{definition}\label{def:SFF}
	\begin{equation}
		\begin{split}
			\mbb S(F,F')=\big\{&s:W\to 2^{E(\left [ F(w), F'(w)\right ])}
				\text{\ such that}\\ &\text{$s(w)$ is a maximum
			matching of $\left[F(w),F'(w)\right]$ for all $w\in W$}\big\}
		\end{split}
	\end{equation}
\end{definition}
Naturally, $s(w)$ is a perfect matching of $[F(w),F(w')]$ if $d_F(w)=d_{F'}(w)$.
(The definition is meaningful even when $d_F(w)=d_{F'}(w)$ does not hold for
every $w$; this will be the case in Lemma~\ref{lemma:numberofmatchings}.)
\medskip
\begin{lemma}\label{th:decomp}
	There is a natural one-to-one correspondence between the family of all
	alternating circuit decompositions of $(F,F')$ and the elements of
	$\S(F,F' )$.
\end{lemma}
\begin{proof}
	If $\C=\{C_1,C_2,\dots,C_n\}$ is an alternating circuit decomposition of
	$(F,F')$, then define $s_{\mathcal{C}}\in \mbb S(F,F')$ as follows:
	\begin{equation}
	\begin{split}
		s_{\mathcal{C}}(w):=\big \{\big(&(w,u),(w,u')\big)\in [F(w), F'(w)]:\\
		   &\text{$(w,u)$ and $(w,u')$ are successive edges in some  $C_i\in \C$} \big \}.
	\end{split}\label{eq:matching}
	\end{equation}
	In turn, to each  $s\in \mbb S(F,F')$ assign  an alternating circuit
	decomposition
	\begin{displaymath}
		\mathcal{C}_s=\{W^s_1,W^s_2\dots, W^s_{p_s}\}
	\end{displaymath}
	of $(F,F')$ as follows: Consider the bipartite graph $\mathcal{F}_s$,
	whose vertex classes are $F$ and $F'$, which are the edges of the simple
	graph $K=(W,F\cup F')$. The edge set of $\mathcal{F}_s$ is
	\begin{equation*}
		E(\mathcal{F}_s)=  \big\{ \big(  (u,w), (u',w)\big ): w\in W \text{\ and }
		\big(  (u,w), (u',w)\big )\in s(w) \big\}.
	\end{equation*}
	In other words, $E(\mathcal{F}_s)$ is the union of $s(w)$ for $w\in W$.
	$\mathcal{F}_{s}$ is a $2$-regular graph, because for each edge
	$(u,v)\in F\cup F'$ there is exactly one $(u,w)\in F\cup F'$ with
	$\big((u,v),(u,w')\big)\in s(u)$, there is exactly one $(t,v)\in
	F\cup F'$ with $\big((u,v),(t,v)\big)\in s(v)$, therefore the
	$\mathcal{F}_{s}$-neighbors of $(u,v)$ are $(u,w)$ and $(t,v)$.

	Since $\mathcal{F}_{s}$ is $2$-regular, it is the union of vertex
	disjoint cycles $\{W^s_i:i\in I\}$. Now $W^s_i$ can also be viewed as a
	sequence of edges in $F\cup F'$, which is a circuit in the graph $K$
	that alternates between $F$ and $F'$. In conclusion, $\{W^s_i:i\in I\}$
	is an alternating circuit decomposition of $(F,F')$. Since
	\begin{displaymath}
		s_{{\mathcal{C}}_{s}}=s,
	\end{displaymath}
	the proof is complete.
\end{proof}

\medskip

If the degree sequence of both $F$ and $F'$ is $(d_1,\dots d_k)$, then write
\begin{equation}\label{eq:tFF}
	t_{F,F'}= \prod_{i=1}^k (d_i!).
\end{equation}
Clearly,
\begin{equation}\label{eq:cardS}
	\left |\S \big (F,F' \big ) \right |=t_{F,F'}.
\end{equation}

\medskip

We are ready to describe $S_{X,Y}$ of Theorem~\ref{th:simplified_sinclair}
($S_{X,Y}$ first appears after Equation~\ref{eq:multiflow_simple1.5}).
\begin{definition}\label{def:SXY}
	For any two graphs $X$ and $Y$ whose degree sequences are identical, let
	\begin{equation*}
		S_{X,Y}:=\S(E(X)-E(Y),E(Y)-E(X)).
	\end{equation*}
\end{definition}
By Lemma~\ref{th:decomp}, every $s\in S_{X,Y}$ corresponds to an
{$(X,Y)$-alternating} circuit decomposition of $\nabla=X\triangle Y$:
\begin{equation}\label{eq:nabla}
 \nabla=W_1^s\uplus W_2^s\uplus \cdots\uplus W_p^s,
\end{equation}
where $\uplus$ is the disjoint union. Since the vertex set of a realization is
$[n]$, the natural ordering on $[n]$ induces a lexicographical order on the
edges, which are unordered pairs of $[n]$. We order the circuits of the
decomposition~\eqref{eq:nabla} in the order of their lexicographically first
edges: for $W^s_i$ and $W^s_j$, we have $i<j$ if and only if the
lexicographically first edge of $W^s_i$ is lexicographically smaller than the
lexicographically first edge of $W^s_j$.

\medskip

In each $W^s_i$, the matching-system $s$ induces an $(X,Y)$-alternating Eulerian
circuit. For readability, we omit $s$ from the superscript of $W_i^s$ in
Section~\ref{sec:design}. The detailed description of~\ref{step1} is now
complete.

\section{Multicommodity flow {-} designing and counting paths in
\texorpdfstring{$\G(\bd)$}{the Markov graph}.}\label{sec:design}

\subsection{Preparatory considerations}\label{sec:cons}

So far, we have described Sinclair's multicommodity flow method in general and
we have learnt how to shape a realization into another one if their difference
is exactly one primitive circuit. We defined our irreducible Markov chain on the
state space $\G$. Finally we described a set of parametrizations of the
difference of any two realizations $X$ and $Y$, and decided that for each such
parametrization we will design one switch sequence which transform $X$ into $Y$.
We have reached the technically most challenging part of any such proof: we have
to design one particular path (switch sequence) for every possible parameter set
$(X,Y,s)$ in such a way that the ensemble does not overload any edge in the
Markov graph. Typically this is not an easy problem. To quote~\cite{DJM17} Dyer,
Jerrum and Müller, ``Achieving low congestion [\ldots] is a delicate matter.''

\medskip

In our case the problems come from two sources. The first one is connected to
the alternating circuit decomposition of the symmetric difference of the two
realizations $X$ and $Y$. It will be relatively simple to keep track of both
$E(X)\cap E(Y)$ and $E(X)\triangle E(Y)$, but extra care is needed to decide
whether an edge of $E(X)\triangle E(Y)$ originates from $X$ or $Y$. At the
beginning of the switch sequence, it is clear which edges come from $X$, and
which come from $Y$, but after processing several primitive circuits, we lose
this information (without the help of further parameters).

\medskip

One can try the following trivial strategy: for each primitive circuit, we
define a Boolean variable, initiated to zero. Whenever we process a primitive
circuit, we change the value of the associated Boolean variable to one. This is
useful when we want to list the edges of $W_k$ in their original order defined
by $s$.

When a primitive circuit $C$ is processed by \textsc{Sweep}, we reverse the
order of its edges compared to the original order $s$. If the original
alternating enumeration is $(a,b,c,d)$, then after \textsc{Sweep} processes $C$,
the order is $(d,c,b,a)$.  The resulting order is an Eulerian circuit of $W_k$
which alternates in the next milestone.  If $X\triangle Y$ is known, then the
Boolean variables tell us which which primitive circuits have already been
processed by \textsc{Sweep}. If the Boolean switch is zero for $C$, then the
current realization restricted to $E(C)$ is identical to $X\cap E(C)$; if the
Boolean switch is one, then the current realization restricted to $E(C)$ is
identical to $Y$. This allows us to restore the original order $s$. Thus the
current values of the variables determine the origin of the edges.

This is a plausible solution if we have a constant number of primitive circuits.
However, if this number is not bounded by a constant, say, it is linear in $n$,
then the cardinality of all possible configurations is exponential in $n$,
therefore the cardinality of possible values taken by the auxiliary parameter
set $B$ is also exponential. This is simply not sufficient to prove rapid
mixing.

\medskip

To produce a successful proof (of the polynomial mixing time of the switch
Markov chain) we shall keep track of both the current alternating edge sequence
and the origin of the edges (whether an edge of $X\triangle Y$ comes from $X$ or
$Y$). Our solution relies on an enumeration algorithm, which, after processing a
certain primitive circuit, will change the trailing order of the edges not only
in the current circuit, but also in a well-defined neighborhood of that circuit.
This is achieved via the $T$-operator (see Section~\ref{sec:T-operator}), which
is the tool that provides the~\ref{reconstructability} property of the
multicommodity-flow we are building. Ultimately, it enables us to reconstruct
the realizations $X$ and $Y$ when at milestones
(Definition~\ref{def:milestone}).

\medskip

The second problem crops up while \textsc{Sweep} processes a primitive circuit.
Algorithms~\ref{alg:sweep} and~\ref{alg:switches} change the status of some
chords which do not belong to the symmetric difference of $X$ and $Y$. We will
use the auxiliary matrix $\widehat M$ introduced in
Definition~\ref{def:widehatM} to track $X\cup Y$ and $X\triangle Y$ using the
current $Z$. Very often, the row- and column-sums of $\widehat M$ correspond to
$\bd$, i.e., $\widehat M$ is the adjacency matrix of a realization of $\bd$.
However, when we switch an edge from outside the symmetric difference, then
$\widehat M$ is not a 0{-}1 matrix anymore.  This does not impact the
reconstructability of $X$ and $Y$ from our auxiliary data, but we have to ensure
that the same $\widehat M$ does not appear too often. Papers~\cite{EKM}
and~\cite{EKMS} used the same coding method as this paper does, but others,
like~\cite{CDG07} or~\cite{G18}, used different parameters. All known proofs
control the number of ``invalid'' non-0{-}1 occurrences, by an appropriate
``critical lemma'' (a term coined in~\cite{G15}). In this paper, this quantity
is controlled directly by $P$-stability.

\medskip

In the next subsection we will introduce and study the $T$-operator to tackle
the first problem mentioned at the beginning of this section.

\subsection{Preparing for~\ref{step2} {-} the
\texorpdfstring{\boldmath$T$}{T}-operator}\label{sec:T-operator}

The $T$-operator is an abstract description of the algorithm we use to decompose
$W_k$ into primitive $X,Y$-alternating circuits, as foretold in
Section~\ref{sec:order}. Suppose the edges of $W_k$ are colored \textit{green}
and $\pi_0$ is a permutation of the edges in which $s\in S_{X,Y}$ traverses
$W_k$, started from a fixed predefined vertex. The $T$-operator will be called
repeatedly, where each iteration outputs a primitive circuit. As its input, the
$T$-operator takes an Eulerian circuit on $W_k$ and a
\textit{red}-\textit{green} coloring of the edges $W_k$. To find a primitive
circuit, $W_k$ is traversed along the Eulerian circuit until the set of visited
\textit{green} edges contains an $X,Y$-alternating primitive circuit $C$. The
output of the $T$-operator modifies the input coloring by recoloring the edges
of $C$ red, moreover, reverses the order of the edges of $C$ and the
\textit{red}-neighborhood of its edges in the Eulerian circuit.

This procedure
tracks the path $\varUpsilon(X,Y,s)$: when the primitive circuit $C$ is
processed by \textsc{Sweep}, the portion of $s$ corresponding to $W_k$ is
modified to encode the Eulerian circuit produced by the $T$-operator. The
modified trail is alternating between edges and non-edges in the current
realization, because the order of the edges of $C$ is reversed by the
$T$-operator.

\medskip

We will show that iteratively applying the $T$-operator to $\pi_0$ and the
identically \textit{green} coloring stabilizes at the reverse (not inverse!) of
$\pi_0$ and the identically \textit{red} coloring. This will show, that the
$T$-operator is indeed capable of producing a primitive circuit decomposition of
$W_k$.

\medskip

Let us list a number of abstract objects. After the definitions, we will roughly
give the correspondence between the abstractions and their application to $W_k$.

Let $[\m ] = \{1,2,\ldots, \m \}$ be a base set, denote by $S_{[\m ]}$ the
symmetric group on $[\m ]$ and let $\Pos$ denote the set of \textbf{position}s,
where $\Pos=\{\posp{1}, \posp{2}, \ldots, \posp{\m -1}\}$.  For convenience, we
consider $\posp{i} = \posm{i+1}$ and allow the alternating naming
$\Pos=\{\posm{2},\ldots, \posm{\m }\}$.  Let $f$ be a two-coloring on $\Pos$
with $f\in {\{\mathit{green}, \mathit{red}\}}^{\Pos}$. We will describe the
state of our system with the pair
\begin{equation*}
	(\pi, f) \quad : \quad \pi\in S_{[\m ]},\ f\in
	{\{\mathit{green},\mathit{red}\}}^{\Pos}.
\end{equation*}
Let $\mathcal{E}\subset \binom{[\m ]}{2}$ be a fixed subset which we call the
set of \textbf{eligible reversals}.  Assume that
\begin{equation}\label{ass:star}
	\text{the connected components of the simple graph }
	([\m ],\mathcal{E})\text{ are cliques.}
\end{equation}
It is important to recognize that {each} eligible reversal {consists of a pair
of} elements of the base set, and they do not depend on the {image of those
elements under $\pi$.} Accordingly, to make the definitions more readable, let
us define
\begin{displaymath}
	\pi^{-1}(\mathcal{E})=\Big\{\{\pi^{-1}(x),\pi^{-1}(y)\}
	: \{x,y\}\in \mathcal{E} \Big\}.
\end{displaymath}
Let us emphasize that $\mathcal{E},\pi,f$ are abstractions, they will gain their
meaning and actual contents later on in Section~\ref{sec:reconst}.

\medskip

In applications, each element in $[\m]$ corresponds to a visit made by an
Eulerian-trail to a vertex of $W_k$. The positions will correspond to the edges
of $W_k$. The set $\mathcal{E}$ pairs the visits which occur at identical
vertices \textbf{with the same parity}. Stated simply, in applications,
$\mathcal{E}$ describes every possible primitive circuit which is contained as a
subsequence in the Eulerian circuit. As mentioned in
Section~\ref{sec:T-operator}, $\pi$ describes an Eulerian circuit on $W_k$, and
$f^{-1}(\mathit{red})$ is the union of the edges of primitive circuits already
processed by \textsc{Sweep}.

\medskip

We now define an operator $T_{\mathcal{E}}$, or $T$ for short, as
$\mathcal{E}$ is fixed anyway. This $T$ is a function mapping $S_{[\m ]}\times
{\{\mathit{green},\mathit{red}\}}^{\Pos}$ into itself. To determine the image of
$(\pi,f)$ under $T$, an interval will be selected first. For that end let
\begin{equation*}
	j_{(\pi,f)}:=\min\left \{j'\in [\m ]\  \Big \vert \ \exists i'<j':\
	f\left (\posp{i'} \right ) =f \left (\posm{j'} \right )=\mathit{green},\
	\{i',j'\}  \in\pi ^{-1}(\mathcal{E}) \right \}
\end{equation*}
then let
\begin{equation*}
	i_{(\pi,f)}:=\max \left \{i'<j_{(\pi,f)}\ \Big \vert \  f\left (\posp{i'}\right ) =\mathit{green},\ \{i',j_{(\pi,f)}\} \in\pi^{-1}(\mathcal{E}) \right \}.
\end{equation*}
We define $\max\emptyset=-\infty$ and $\min\emptyset=+\infty$.  For any integer
$k : 1 \le k \le \m $ we select two positions from $\Pos$. Let
$a_{(\pi,f)}(k):=k$ if $f(\posm{k})=\mathit{green}$, and let
\begin{equation*}
	a_{(\pi,f)}(k) :=\min\Big\{ i'\le k\ \big \vert \  \forall i''
	\text{\ s.t.\ }i'\le i''<k\ : \  f\left (\posp{i''}\right )=\mathit{red}\Big\}
\end{equation*}
otherwise. Furthermore, let $b_{(\pi,f)}(k):=k$ if $f(\posp{k})=\mathit{green}$,
and let
\begin{equation*}
	b_{(\pi,f)}(k) :=\max\Big\{ j'\ge k\ \big \vert \  \forall j''
	\text{\ s.t.\ }k<j''\le j'\ : \  f\left (\posm{j''} \right )=\mathit{red}\Big\}
\end{equation*}
otherwise. By definition, $a_{(\pi,f)}(k) \leq k\leq b_{(\pi,f)}(k)$ for all
$k\in [\m ]$. For a visualization, see the top half of Figure~\ref{fig:demonst}.

\medskip

For any $k\in [\m ]$, let ${\overleftarrow
k^{(\pi,f)}:=a_{(\pi,f)}(k)+b_{(\pi,f)}(k)-k}$. We omit the index $_{(\pi,f)}$
in the following. Informally, $\overleftarrow k$ takes the largest red segment
$[a_{(\pi,f)}(k),b_{(\pi,f)}(k)]$
around $k$, and flips $k$ around the center of the
segment; see the bottom half of Figure~\ref{fig:demonst}.

\begin{figure}[ht!]
	\centering
	\begin{tikzpicture}
		\begin{scope}
			\def\n{11}
			\node at (7,2) {$\mathcal{E}$};
			\node at (-1,0) {$(\pi:x\mapsto x,f)$};
			\foreach \x in {1,...,\n}
			\node[draw,circle,inner sep=0pt,minimum size=12pt] (v\x) at (\x,0) {\tiny \textrm{\x}};
			\foreach \x/\c [count=\y from 2] in {1/red,2/green,3/red,4/green,5/green,6/red,7/red,8/green,9/green,10/green}
			\draw[very thick,\c] (v\x) -- (v\y);
			\foreach \x/\y in {1/2,6/8,2/6,1/6,2/8,1/8,1/11,2/11,6/11,8/11,10/3,10/4,4/3,10/9,9/3,9/4}
			\draw (v\x)  edge[bend left=60] (v\y);
			\node (ai) at (1,-0.5) {$a(i)$};
			\node (i) at (2,-0.5) {$i$};
			\node (j) at (6,-0.5) {$j$};
			\node (bj) at (8,-0.5) {$b(j)$};
		\end{scope}
		\begin{scope}[yshift=-150]
			\def\n{11}
			\node at (7,2) {$\mathcal{E}$};
			\node at (-1,0)
				{$(x\mapsto \pi(\overleftarrow{x}),f)$};
			\foreach \y [count=\x from 1] in
			{2,1,4,3,5,8,7,6,9,10,11}
			\node[draw,circle,inner sep=0pt,minimum size=12pt] (v\x) at (\x,0) {\tiny \textrm{\y}};
			\foreach \x/\c [count=\y from 2] in {1/red,2/green,3/red,4/green,5/green,6/red,7/red,8/green,9/green,10/green}
			\draw[very thick,\c] (v\x) -- (v\y);
			\foreach \x/\y in {1/2,6/8,2/6,1/6,2/8,1/8,1/11,2/11,6/11,8/11,10/3,10/4,4/3,10/9,9/3,9/4}
			\draw (v\x)  edge[bend left=60] (v\y);
		\end{scope}
	\end{tikzpicture}
	\caption{An demonstrative example where $\pi=\mathrm{id}_{11}$. The
		curved arcs represent the pairs in $\mathcal{E}$. The encircled
		numbers are $\pi(x)$ and $\pi(\protect\overleftarrow{x})$,
		respectively. Do note, that the displayed $\mathcal{E}$ cannot
		occur in applications, as it connects vertices of different
		parity.}\label{fig:demonst}
\end{figure}

\begin{definition}[$T$-operator]\label{def:T-operator}
	We define the function
	\begin{equation*}
		T:S_{[\m]}\times {\{\mathit{green},\mathit{red}\}}^\Pos\to
		S_{[\m]}\times {\{\mathit{green},\mathit{red}\}}^\Pos.
	\end{equation*}
	Let $(\pi,f)\in S_{[\m]}\times
	{\{\mathit{green},\mathit{red}\}}^\Pos$. If $j_{(\pi,f)}=+\infty$, then
	let $T(\pi,f):=(\pi,f)$ be a fixed point. If $j_{(\pi,f)}$ is finite,
	define ${T:(\pi,f)\mapsto (\pi',f')}$ as follows:
	\begin{align*}
	\pi'(k)          & =\left\{
		\begin{array}{ll}
			\pi\left(\overleftarrow{a(i)+b(j)-k}\right) & \text{ if }k\in
			[a(i),b(j)],\\
			\pi(k)                           & \text{ if }k\notin [a(i),b(j)],
		\end{array}
	\right.                     \\
		{f'}{(\posp{k})} & =\left\{
			\begin{array}{ll}
				f(\posp{k})  & \text{ if }1\le k < a(i)\text{ or }b(j)\le k <\m , \\
				\mathit{red} & \text{ if }a(i)\le k<b(j).
			\end{array}
		\right.
	\end{align*}
\end{definition}

\begin{figure}[!ht]
	\centering
	\begin{tikzpicture}
		\begin{scope}
			\def\n{11}
			\node at (7,2) {$\mathcal{E}$};
			\node at (-0.5,0) {$(\pi,f)$};
			\foreach \x in {1,...,\n}
			\node[draw,circle,inner sep=0pt,minimum size=12pt] (v\x) at (\x,0) {\tiny \textrm{\x}};
			\foreach \x/\c [count=\y from 2] in {1/red,2/green,3/red,4/green,5/green,6/red,7/red,8/green,9/green,10/green}
			\draw[very thick,\c] (v\x) -- (v\y);
			\foreach \x/\y in {1/2,6/8,2/6,1/6,2/8,1/8,1/11,2/11,6/11,8/11,10/3,10/4,4/3,10/9,9/3,9/4}
			\draw (v\x)  edge[bend left=60] (v\y);
			\node (ai) at (1,-0.5) {$a(i)$};
			\node (i) at (2,-0.5) {$i$};
			\node (j) at (6,-0.5) {$j$};
			\node (bj) at (8,-0.5) {$b(j)$};
		\end{scope}
		\begin{scope}[yshift=-150]
			\def\n{11}
			\node at (7.5,1.5) {$\mathcal{E}$};
			\node at (-0.5,0) {$(\pi',f')$};
			\foreach \y [count=\x from 1] in {6,7,8,5,3,4,1,2,9,10,11}
			\node[draw,circle,inner sep=0pt,minimum size=12pt] (v\y) at (\x,0) {\tiny \textrm{\y}};

			\foreach \x/\y/\c in {6/7/red,7/8/red,8/5/red,5/3/red,3/4/red,4/1/red,1/2/red,2/9/green,9/10/green,10/11/green}
			\draw[very thick,\c] (v\x) -- (v\y);
			\foreach \x/\y in {1/2,6/8,6/2,6/1,8/2,8/1,1/11,2/11,6/11,8/11,10/3,10/4,4/3,10/9,9/3,9/4}
			\draw (v\x)  edge[bend left=60] (v\y);
		\end{scope}
	\end{tikzpicture}
	\caption{An example for $T(\pi,f)=(\pi',f')$. The curved arcs represent
	the pairs in $\mathcal{E}$. The encircled numbers are $\pi(x)$ and
	$\pi'(x)$, respectively, where {$x=1,\ldots, 11$ from left to right ($\pi$ is
	identity).}}\label{fig:Toperator}
\end{figure}

Observe Figure~\ref{fig:Toperator}. The $T$-operator reverses the order of the
\textit{green} sections in $[a(i),b(j)]$, such that the order is not reversed
on the already \textit{red} sections. Then the region $[a(i),b(j)]$ is colored red.

Let us give another, more verbose description of $T(\pi,f)$.
Let
\begin{equation*}
\left\{ \posp{a_{(\pi,f)}(i_{(\pi,f)})}, \ldots,\posm{b_{(\pi,f)}(j_{(\pi,f)})}\right\}
\end{equation*}
be the maximal interval (with integer endpoints) containing
$\left\{ \posp{i_{(\pi,f)}},\ldots, \posm{j_{(\pi,f)}} \right\} $
such that the $f$-image of the new positions of the extended interval are $\mathit{red}$. Take a look at Figure~\ref{fig:Toperator}. To construct $\pi'$ from $\pi$, every maximal $\mathit{red}$ interval $\{x^+,\ldots,y^-\}$ in $\{a^+,\ldots,b^-\}$ is shifted to $\{ {(a+b-y)}^+,\ldots,{(a+b-x)}^- \}$, and the $\mathit{green}$ positions {within $\{ a^+,\ldots, b^-\}$} are taken in reverse order in the remaining positions between the shifted $\mathit{red}$ intervals.

Now that $T$ is defined, we will use $T$ iteratively. For $r\in\mathbb{N}$, let
\begin{equation*}
	T^r=\overbrace{T\circ T\circ \ldots \circ T}^{r}.
\end{equation*}
Given any permutation $\pi_0$ on $[\m ]$, let $(\pi_r,f_r):=T^r(\pi_0,\mathbf{green})$, where $\mathbf{green}$ is the identically $\mathit{green}$ function.
In subscripts, we shorten $(\pi_r,f_r)$ by writing $r$ instead. For example, $j_r=j_{(\pi_r,f_r)}$, etc.
\begin{lemma}\label{lemma:maxred}
	For any $r\geq 0$, the pair of endpoints of a maximal path formed by elements of $f_r^{-1}(\mathit{red})$ is an element of $\pi^{-1}_r(\mathcal{E})$.
\end{lemma}

\begin{proof}
	The statement is vacuously true when $r=0$.
	Use property~\eqref{ass:star} and the fact that $f_r(\posp{i_r})=f_r(\posm{j_r}) =\mathit{green}$. By induction, either $a_r(i_r) =i_r$ or $\{a_r(i_r),i_r\}\in\pi_r^{-1}(\mathcal{E})$. By definition, $\{i_r,j_r\}\in\pi_r^{-1}(\mathcal{E})$, thus we also have $\{a_r(i_r),j_r\}\in\pi_r^{-1}(\mathcal{E})$. The same argument goes through for $j_r$ and $b_r(j_r)$.
\end{proof}

\begin{lemma}\label{lemma:jinc}
	The following statements hold for $r\geq 0$.
	\begin{itemize}
		\item[\emph{(i)}] If $(\pi_r,f_r)$ is not a fixed point of the $T$ operator then $j_r < j_{r+1}$;
		\item[\emph{(ii)}] $b_r(j_r)=j_r$;
		\item[\emph{(iii)}] $f_{r+1}((k)^+)=\mathit{green}$ for $j_{r}\leq k\leq \m $;
	\end{itemize}
\end{lemma}

\begin{proof}
	We proceed by induction on $r$. Statements (ii), (iii) are true when $r=0$.
	Now suppose that $r\geq 0$ is such that statements (ii), (iii) are true for $r$.
	If $(\pi_r,f_r)$ is a fixed point of $T$ then we are done.
	Otherwise, since
	\[
		f_{r+1}^{-1}(\mathit{green})\subsetneq f_r^{-1}(\mathit{green}),
	\]
	we have $f_r(\posp{i_{r+1}})=f_r(\posm{j_{r+1}})=\mathit{green}$.
	Clearly, $\{\pi_{r+1}(i_{r+1}),\pi_{r+1}(j_{r+1})\}\in\mathcal{E}$, so
	\begin{equation}\label{eq:preimage}
		\left\{\pi_r^{-1}(\pi_{r+1}(i_{r+1})),\pi_r^{-1}(\pi_{r+1}(j_{r+1}))\right\}\in\pi_r^{-1}(\mathcal{E}).
	\end{equation}
	Since $f_{r+1}(\posm{j_{r+1}})=\mathit{green}$, we must have $j_{r+1}>b_r(j_r)=j_r$ or $j_{r+1}\le a_{r}(i_r)$. The first case gives $j_r < j_{r+1}$ immediately.

	\medskip

	Next, suppose that $j_{r+1}< a_r(i_r)$. Then $\pi_{r+1}(j_{r+1})=\pi_r(j_{r+1})$ and
	$\pi_{r+1}(i_{r+1})=\pi_r(i_{r+1})$. Plugged into Equation~\eqref{eq:preimage}, the definition of $j_r$ implies that $j_{r+1}\ge j_r > a_r(i_r)$, a contradiction.

	\medskip

	Finally, suppose that $j_{r+1}=a_r(i_r)$.
	Recalling that $\pi_{r+1}(a_r(i_r)) = \pi_r(j_r)$, we have
	\[\left\{\pi_r^{-1}(\pi_{r+1}(i_{r+1})),\pi_r^{-1}(\pi_{r+1}(j_{r+1}))\right\}=
	\left\{i_{r+1}, j_r\right\}. \]
	Now Equation~\eqref{eq:preimage} gives
	\[\left\{i_{r+1},j_r\right\}\in\pi_r^{-1}(\mathcal{E}),\]
	and property~\eqref{ass:star} implies that
	$\left\{ i_{r+1},i_r\right\}\in \pi_r^{-1}(\mathcal{E})$.
	If $f_r((i_r)^-) = \mathit{green}$ then this contradicts the definition of $j_r$.
	Otherwise, $\{a_r(i_r),\ldots, i_r\}$ is a maximal $\mathit{red}$ interval in $f_r$ and
	hence, by Lemma~\ref{lemma:maxred} and property~\eqref{ass:star} we conclude that
	$\left\{ i_{r+1},a_r(i_r)\right\}\in \pi_r^{-1}(\mathcal{E})$.
	Again, this contradicts the choice of $j_r$.

	Hence in all cases we conclude that $j_r < j_{r+1}$, which implies that
	$b_{r+1}(j_{r+1})=j_{r+1}$ and $f_{r+2}((k)^+) = \mathit{green}$ for $j_{r+1}\leq k\leq \m $.
	This completes the proof.
\end{proof}

The next lemma is very important for understanding how the $T$-operator works.
In Definition~\ref{def:T-operator}, we could have taken
$\pi(\overleftarrow{i+j-k})$ for $k\in [i,j]$, and so far, every stated lemma
would still hold. The reason we take $a(i)$ instead of $i$ is to ensure that the
order $\pi_r$ on maximal red regions in $f_r$ are just the reverse of their
original order in $\pi_0$.
\begin{lemma}\label{lemma:pi_r}
	For arbitrary $\pi_0$, $r\ge 0$, and $k\in [\m ]$, we have
	\[\pi_r(k)=\pi_0\Big(a_r(k)+b_r(k)-k\Big).\]
\end{lemma}
\begin{proof}
	If $r=1$, the statement immediately follows from the definition. Suppose
	the statement holds for $r-1$. If $k\notin
	[a_{r-1}(i_{r-1}),b_{r-1}(j_{r-1})]$ then $a_{r-1}(k)=a_r(k)$ and
	$b_{r-1}(k)=b_r(k)$, so
	\[\pi_r(k)=\pi_{r-1}(k)=\pi_0(a_{r-1}(k)+b_{r-1}(k)-k)=\pi_0(a_{r}(k)+b_{r}(k)-k),\]
	as we wished.

	\medskip

	Suppose from now on that $k\in [a_{r-1}(i_{r-1}),b_{r-1}(j_{r-1})]$.
	Let
	\begin{equation*}
		\ell=a_{r-1}(i_{r-1})+b_{r-1}(j_{r-1})-k.
	\end{equation*}
	Since the edges in $[a_{r-1}(i_{r-1}),b_{r-1}(j_{r-1})]$ are all
	$\mathit{red}$ in $f_r$, we have
	\[
		a_{r}(k)=a_{r-1}(i_{r-1}),\quad b_{r}(k)=b_{r-1}(j_{r-1}).
	\]
	Writing $\overleftarrow{\ell}^{(r-1)}$ for $\overleftarrow{\ell}^{(\pi_{r-1},f_{r-1})}$, by induction
	we have
	\begin{equation*}
		\pi_r(k)=\pi_{r-1}\left(\overleftarrow \ell^{(r-1)}\right)=\pi_0\left(a_{r-1}\left(\overleftarrow \ell^{(r-1)}\right)+b_{r-1}
		\left(\overleftarrow \ell^{(r-1)}\right)-\overleftarrow \ell^{(r-1)}\right).
	\end{equation*}
	Since $a_{r-1}\left(\overleftarrow \ell^{(r-1)}\right)=a_{r-1}(\ell)$ and $b_{r-1}\left(\overleftarrow \ell^{(r-1)}\right)=b_{r-1}(\ell)$, the right hand side is equal to $\pi_0(\ell)$. Expanding it, we get
	\begin{equation*}
		\pi_0(\ell)=\pi_0(a_{r-1}(i_{r-1})+b_{r-1}(j_{r-1})-k)=\pi_0(a_{r}(k)+b_{r}(k)-k),
	\end{equation*}
	which is what we intended to prove.
\end{proof}

Repeated application of the $T$-operator turns every position \textit{red}
eventually, if the trivial necessary condition is satisfied:
\begin{lemma}\label{lemma:complete}
	If $\{1,\m \}\in\mathcal{E}$, then $\exists z\in \mathbb{N}$ such that
	$f_z^{-1}(\mathit{red})=\Pos$.
\end{lemma}
\begin{proof}
	Lemma~\ref{lemma:maxred} implies that $\{1,\min\{t: f_r(\posp{t})=
	\mathit{green}\}\}\in\pi^{-1}_r(\mathcal{E})$, therefore 
	\begin{equation*}
	\{\min\{t: f_r(\posp{t})= \mathit{green}\},\m \}\in\pi^{-1}_r(\mathcal{E}),
	\end{equation*}
	except if $f_r^{-1}(\mathrm{red})=\Pos$ already.
\end{proof}

The following lemma shows that those segments whose endpoints form an eligible
reversal cannot simultaneously start on a \textit{green} position and end on a
\textit{red} position.
\begin{lemma}\label{lemma:nogreenred}
	Given $r\ge 0$ and any $\{x,y\}\in\pi^{-1}_r(\mathcal{E})$ such that $x<y$, either $\{x,y\}=\{i_r,j_r\}$, or $f(\posp{x})=\mathit{red}$, or $y\ge j_r+1$.
\end{lemma}
\begin{proof}
	The lemma trivially holds for $r=0$. Suppose now, that $r\ge 1$.

	\medskip

	Suppose first, that $f_r(\posp{x})=f_r(\posm{y})=\mathit{green}$. By
	definition, $y\ge j_r$. If $y\ge j_r+1$, the lemma holds. If $y=j_r$, then definition of $i_r$ implies that $x\le i_r$. If $y=j_r$ and $x<i_r$, then $x<a_r(i_r)$. By property~\eqref{ass:star} {and Lemma~\ref{lemma:maxred}}, $\{x,a_r(i_r)\}\in\pi^{-1}_r(\mathcal{E})$ holds. Since $f_r(\posm{a_r(i_r)})=\mathit{green}$, we have a contradiction with the definition of $j_r$.

	\medskip

	Suppose, that $f_r(\posp{x})=\mathit{green}$, $f_r(\posm{y})=\mathit{red}$, and the lemma does not hold. By Lemma~\ref{lemma:jinc}
	we have $y\le j_{r-1}$. Then we must also have $x<a_{r-1}(i_{r-1})$
	(otherwise $f_r(\posp{x})=\mathit{red}$, a contradiction). Thus
	\begin{equation*}
		\pi^{-1}_{r-1}(\pi_r(x))=x\text{ and }
		x<\pi^{-1}_{r-1}(\pi_r(y))\le j_{r-1}
	\end{equation*}
	so $\{x,\pi^{-1}_{r-1}(\pi_r(y))\}\in
	\pi^{-1}_{r-1}(\mathcal{E})$. By induction, we should have
	$f_{r-1}(\posp{x})=\mathit{red}$, which implies
	$f_{r}(\posp{x})=\mathit{red}$, a contradiction.

	\medskip

	We have checked and eliminated every possible case where the statement of the lemma is not satisfied.
\end{proof}

The next lemma provides our desired~\ref{reconstructability}
property. By including a natural number $w$ between 0 and $n^2$ in the parameter
set $B$ in Equation~\eqref{eq:additionalParams}, we need not know the current
coloring $f_r$ to reconstruct $\pi_0$ from $\pi_r$:
\begin{theorem}\label{thm:Toperator}
	$\forall r\in\mathbb{N}\ \exists w\in \mathbb{N}$ and
	$\exists g\in{\{\mathit{green},\mathit{red}\}}^{\Pos}$ such that
	\begin{equation*}
		T^w(\pi_r,\mathbf{green})=(\pi_0,g).
	\end{equation*}
\end{theorem}
\begin{proof}
	If $f_r(\Pos)\equiv \mathit{red}$, then Lemma~\ref{lemma:maxred} implies
	that $\{1,\m \}\in\pi^{-1}_r(\mathcal{E})$.  By Lemma~\ref{lemma:pi_r},
	$\pi_r(k)=\pi_0(1+\m-k)$.  By Lemma~\ref{lemma:complete}, there exists a
	$\pi\in S_{[\m]}$ and a $z\in\mathbb{N}$ for which
	$T^z(\pi_r,\mathbf{green})=(\pi,\mathbf{red})$. By
	Lemma~\ref{lemma:pi_r}, we have
	\begin{equation*}
		\pi(k)=\pi_r(1+\m-k)=\pi_0(k),
	\end{equation*}	
	which is what we wanted.

	\medskip

	If $f_r^{-1}(\mathit{red})$ is composed of multiple components, then a
	repeated application of the $T$-operator works successively in these
	components. Lemma~\ref{lemma:pi_r} says that the order of elements in
	each of these components have been reversed in $\pi_r$ compared to
	$\pi_0$. Outside these intervals, however, $\pi_r$ is identical to
	$\pi_0$.

	\medskip

	Because of Lemma~\ref{lemma:maxred}, we see that
	Lemma~\ref{lemma:complete} implies that the maximal $\mathit{red}$
	intervals will be completely processed after a certain number of steps.
	Lemmas~\ref{lemma:jinc} and~\ref{lemma:nogreenred} together imply that
	if the $T$-operator starts working inside a component of
	$f_r^{-1}(\mathit{red})$ then the next selected interval $[i,j]$ will
	also be inside until the whole component becomes $\mathit{red}$ again.
	At this point, by Lemma~\ref{lemma:pi_r}, the order of the elements
	inside each $\mathit{red}$ component have been reversed for a second
	time, so the final permutation equals $\pi_0$, as claimed.
\end{proof}

\begin{definition}
	The restriction of a permutation $\pi$ to an interval $[\alpha,\beta]$,
	denoted by $\pi|_{[\alpha,\beta]}$, means that original domain $[\m ]$
	of $\pi$ is replaced with $[\alpha,\beta]$.
\end{definition}
The $T$-operator naturally generalizes to injective maps from $[\alpha,\beta]$
to an arbitrary set.

\begin{lemma}\label{lemma:Trestricted}
	Suppose that $r\in\mathbb{N}$ and $[\alpha,\beta]\subset [i_{r},j_{r}]$
	(where $\alpha<\beta$) is a proper subinterval such that any component
	of $f_r^{-1}(\mathit{red})$ is either disjoint from
	$\Pos[\alpha,\beta]:=\{\posp{\alpha},\ldots,\posm{\beta}\}$ or entirely
	contained by it. If
	$f_r(\posp{\alpha})=f_r(\posm{\beta})=\mathit{green}$, then
	\begin{align*}
		&(\pi_{r}|_{[\alpha,\beta]},f_{r}|_{\Pos[\alpha,\beta]})\text{\
		is a fixed point of the $T$-operator and}\\
		&(\pi_{r}|_{[\alpha,\beta]},f_{r}|_{\Pos[\alpha,\beta]})=
		T^r(\vartheta,\mathbf{green}),
	\end{align*}
	where $\vartheta=\pi_{0}|_{(\alpha,\beta)}\cup
	\{\alpha\mapsto\pi_r(\alpha),\beta\mapsto\pi_r(\beta)\}$.
\end{lemma}
\begin{proof}
	Both $\pi_{0}|_{[\alpha,\beta]},\pi_{r}|_{[\alpha,\beta]}$ map
	$[\alpha,\beta]$ to the same set, because of the assumption on the
	components of $f_r^{-1}(\mathit{red})$.
	Lemma~\ref{lemma:nogreenred} implies that $\{\alpha,\beta\}\notin \mathcal{E}$
	and that $(\pi_r|_{[\alpha,\beta]}|f_r|_{[\alpha,\beta]})$ is a fixed
	point of the $T$-operator. Lemma~\ref{lemma:maxred} guarantees that for
	$k\le r$, either $\alpha<i_{(\pi_k,f_k)}<j_{(\pi_k,f_k)}<\beta$, or
	$[i_{(\pi_k,f_k)},j_{(\pi_k,f_k)}]\cap (\alpha,\beta)=\emptyset$.
	In the latter case, it is possible that $j_{(\pi_k,f_k)} = \alpha$ or
	$i_{(\pi_k,f_k)}=\beta$, which is why we defined
	$\vartheta(\alpha)$ and $\vartheta(\beta)$ separately.
	In the former case, the $T$-operator and the restriction operation
	trivially commute:
	\begin{equation*}
		T(\pi_k|_{[\alpha,\beta]},f_k|_{[\alpha,\beta]})=
		T(\pi_k,f_k)|_{[\alpha,\beta]}.
	\end{equation*}
	This implies that 
	\begin{equation*}
		(\pi_{r}|_{[\alpha,\beta]},f_{r}|_{\Pos[\alpha,\beta]})=
		T^{q}(\vartheta,\mathbf{green})
	\end{equation*}
	for some $q\le r$. Since the left hand side is a fixed point of the
	$T$-operator, applying $T^{r-q}$ acts identically on it. 
\end{proof}

\subsection{\ref{step2} - decomposing a circuit into primitive circuits}\label{sec:reconst}

Given $X,Y\in \mathbb{G}$ (we do not specify whether the degree sequence is
unconstrained or bipartite), and $s\in S_{X,Y}$, we construct a path between $X$
and $Y$ in $\mathbb{G}$ as follows. Recall Equation~\eqref{eq:nabla}: the
matching-system $s$ decomposes $\nabla=X\triangle Y$ into alternating circuits
\begin{equation}\label{eq:circuits}
	\nabla=W_1\uplus\ldots\uplus W_{p}.
\end{equation}
The number $p$ and the circuits $W_k$ depend on $s$; we omit $s$ from super- and
subscripts, as usual. The goal in this section is to apply the $T$-operator to
decompose each circuit $W_k$ into primitive circuits, and show
that~\ref{reconstructability} holds.

Let us fix an arbitrary circuit $W_k$ from the collection~\eqref{eq:circuits}.
The vertices of $W_k$ are a subset of $\mathbb{N}$, so there is a natural
ordering on them.  Let $v_1v_2$ be the lexicographically first edge of $W_k$,
where $v_1$ precedes $v_2$. As stated below Equation~\ref{eq:nabla}, the
matching system $s$ induces an $(X,Y)$-alternating Eulerian circuit on $W_k$.
Let us list this Eulerian trail starting with the edge $v_1v_2$: 
\begin{equation}\label{eq:WkEulerianCircuit}
	W_k:\ (v_1,v_2,v_3,v_4,\ldots,v_{|E(W_k)|},v_{|E(W_k)|+1}),\text{\ where\ }v_{|E(W_k)|+1}=v_1.
\end{equation}
Let $\m =|E(W_k)|+1$,
$\pi_0 = \mathrm{id}_{[\m ]}$, $f_0=\mathbf{green}$, and
\begin{equation}\label{eq:cE}
	\mathcal{E}=\left\{ \{x,y\}\in\binom{[\m ]}{2}
	\ \Big.:\ v_x=v_y\text{\ and\ }x\equiv y\pmod2 \right\}.
\end{equation}
By transitivity, this set possesses property~\eqref{ass:star}, so we can apply
the $T$-operator on $\pi_0$ with $\mathcal{E}$ as the set of eligible reversals.
Let $(\pi_r,f_r)=T^r(\pi_0,f_0)$ for $r\in\mathbb{N}$.
\begin{lemma}\label{lemma:trail}
	Given $r\in\mathbb{N}$, visiting the vertices
	$v_{\pi_r(1)}v_{\pi_r(2)}\dots v_{{\pi_r(|E(W_k)|+1)}}$ is
	an Eulerian circuit of $W_k$.
\end{lemma}
\begin{proof}
	Easily seen by induction on $r$. Lemma~\ref{lemma:maxred} and the
	definition of the $T$ operator implies that we get $\pi_r$ by reversing
	some intervals of the trail defined by $\pi_{r-1}$ whose first and last
	vertices are identical. Consequently, every edge is visited by the new
	trail too.
\end{proof}

It is clear by definition, that $\posp{x}$ in the set $\Pos$ coincides with
$v_{\pi_r(x)}v_{\pi_r(x+1)}$ on the Eulerian circuit determined by $\pi_r$ in
Lemma~\ref{lemma:trail}. Subsequently, $f_r$ defines a corresponding coloring:
the edge $v_{\pi_r(x)}v_{\pi_r(x+1)}$ has color $f_r(\posp{x})$. By
Lemma~\ref{lemma:pi_r}:
\begin{equation}\label{eq:vpirx}
v_{\pi_r(x)}v_{\pi_r(x+1)}=\left\{
	\begin{array}{ll}
		v_{x}v_{x+1} & \text{if }f_r(\posp{x})=\mathit{green}, \\
		v_{a_r(x)+b_r(x)-x-1}v_{a_r(x)+b_r(x)-x} &
		\text{if }f_r(\posp{x})=\mathit{red}.
	\end{array}
	\right.
\end{equation}
By definition,
\begin{equation}\label{eq:rededges}
\{v_{a_r(x)+b_r(x)-x-1}v_{a_r(x)+b_r(x)-x} \,:\,f_r(\posp{x})=\mathit{red}\}=
\{v_{x}v_{x+1}\,:\,f_r(\posp{x})=\mathit{red}\}.
\end{equation}

To complete Definition~\ref{def:milestone}, we shall
determine the primitive circuits $C^k_r$.

\begin{definition}\label{def:Ckr}
Let $\ell_k$ be the maximum $\ell$ for which $\pi_{\ell-1}\neq\pi_{\ell}$.
For any $1\le r\le \ell_k$, let
\begin{align}
	\begin{split}
		E(C^k_r) :=&\left \{ v_{x}v_{x+1}\, :\,  x\in [i_{r-1},
		j_{r-1}-1]\text{\ and\ } \posp{x} \in f_{r-1}^{-1}
	(\mathit{green})  \right \}=\\
			=&\left \{ v_{x}v_{x+1}\, :\,  \posp{x} \in f_{r}^{-1}(\mathit{red})
	\setminus f_{r-1}^{-1}(\mathit{red})\right \}.
	\end{split}\label{eq:Ckr} 
\end{align}
The equality of the two right hand sides in~\eqref{eq:Ckr} follows
from~\eqref{eq:rededges}. Take the list of edges of $W_k$ starting with $v_1v_2$
in the order defined by the Eulerian circuit $\pi_0=\mathrm{id}_{[\m ]}$. This
order can be restricted to the edges of $C^k_r$, so there is a natural Eulerian
circuit on $C^k_r$, too.
\end{definition}

Equations~\eqref{eq:vpirx} and~\eqref{eq:Ckr} show that $C^k_r$ is a
subsequence of edges of the Eulerian circuit defined by $\pi_r$ on $W_k$:
\begin{equation}\label{eq:Ckr2}
	E(C^k_r)=\left \{ v_{\pi_r(x)}v_{\pi_r(x+1)}\, :\,  x\in [i_{r-1},
		j_{r-1}-1]\text{\ and\ } \posp{x} \in f_{r-1}^{-1}
	(\mathit{green})  \right \}
\end{equation}
We want to show that $C^k_r$ is $(X,Y)$-alternating. We need the following
lemma.
\begin{lemma}\label{lemma:cEpreserved}
	For any $r\in\mathbb{N}$, we can describe $\pi^{-1}_r(\mathcal{E})$ as
	the set of endpoints of even circuits formed by subintervals of the
	Eulerian circuit defined by $\pi_r$ in $W_k$:
	\begin{equation*}
	\pi^{-1}_r(\mathcal{E})=\left\{ \{x,y\}\in\binom{[\m ]}{2}\ \Big.:
	\ v_{\pi_r(x)}=v_{\pi_r(y)}\text{ and }x\equiv y\pmod2 \right\}.
	\end{equation*}
\end{lemma}
\begin{proof}
	From~\eqref{eq:cE}, we have
	\begin{align*}
		\pi^{-1}_r(\mathcal{E}) & =\left\{ \{\pi^{-1}_r(x),\pi^{-1}_r(y)\}\ \Big.:\ \{x,y\}\in\binom{[\m ]}{2},\ v_x=v_y\text{ and }x\equiv y\pmod2 \right\} \\
					& =\left\{ \{x,y\}\in\binom{[\m ]}{2}\ \Big.:\ v_{\pi_r(x)}=v_{\pi_r(y)}\text{ and }\pi_r(x)\equiv \pi_r(y)\pmod2 \right\},
	\end{align*}
	because $\pi_r$ is a permutation. It is enough to show that $\pi_r$ preserves parity, i.e.,\ $\pi_r(x)\equiv x\pmod 2$ for any $x$. For $r=0$ this is trivial. Suppose $\pi_{r-1}$ preserves parity.
	For $x\notin [a_{r-1}(i_{r-1}),b_{r-1}(j_{r-1})]$, we have $\pi_r(x)\equiv\pi_{r-1}(x)\pmod 2$, so parity is preserved.

	\medskip

	Now suppose that $x\in [a_{r-1}(i_{r-1}),b_{r-1}(j_{r-1})]$.
	First observe that Lemma~\ref{lemma:maxred} implies that $a_{r-1}(y)\equiv b_{r-1}(y)\pmod 2$ for arbitrary $y$.
	Hence $a_{r-1}(i_{r-1})\equiv i_{r-1}\pmod 2$, and
	since $\{ i_{r-1},j_{r-1}\}\in \pi_{r-1}(\mathcal{E})$
	we know that $i_{r-1}\equiv j_{r-1}$.
	As $j_{r-1} =b_{r-1}(j_{r-1})$ by Lemma~\ref{lemma:jinc}, it follows that $a_{r-1}(i_{r-1})\equiv
	b_{r-1}(j_{r-1})\pmod 2$.
	Now
	\[ \pi_r(x) = \pi_{r-1}\big(a_{r-1}(z) + b_{r-1}(z) - z\big)\]
	where
	\[ z = a_{r-1}(i_{r-1}) + b_{r-1}(j_{r-1}) - x \, \equiv \, x\pmod{2}.\]
	Hence, by induction,
	\[ \pi_r(x)\equiv a_{r-1}(z) + b_{r-1}(z) - z \, \equiv \, z\, \equiv \, x\pmod{2},\]
	completing the proof.
\end{proof}

\begin{lemma}\label{lemma:primitivecircuits}
	$C^k_r$ is a primitive $(X,Y)$-alternating circuit for every $1\le r\le
	\ell_k$. Moreover,
	\begin{equation*}
		W_k=C^k_1\uplus C^k_2\uplus\ldots\uplus C^k_{\ell_k}.
	\end{equation*}
\end{lemma}
\begin{proof}
	Since $W_k$ is an alternating circuit, $|E(W_k)|$ is divisible by two,
	so $\{1,|E(W_k)|+1\}\in\mathcal{E}$. According to
	Lemma~\ref{lemma:complete} there exists a smallest $\ell_k\in
	\mathbb{N}$ such that $f_{\ell_k}(v_{x}v_{x+1})=\mathit{red}$ for every
	$x=1,\ldots,|E(W_k)|$.  As $r$ increases, $\mathit{red}$ edges stay
	$\mathit{red}$, so it follows from Equation~\ref{eq:Ckr} that
	$\bigcup_{r=1}^{\ell_k}C^k_r$ is an edge-disjoint partition of $W_k$.

	\medskip

	Recall \eqref{eq:Ckr2}. By Lemma~\ref{lemma:maxred} and
	Lemma~\ref{lemma:cEpreserved}, $C^k_r$ is an even circuit.  The proof of
	Lemma~\ref{lemma:cEpreserved} also shows that $\pi_r$ preserves parity.
	By \eqref{eq:Ckr}, it is such a subsequence of the Eulerian circuit
	$\pi_{0}$, which starts with the odd indexed $v_{i_{r-1}}$ and ends with
	the odd indexed $v_{j_{r-1}}$; moreover, by \eqref{eq:rededges} and
	Lemma~\ref{lemma:cEpreserved}, the \textit{red} edges between
	$v_{i_{r-1}}$ and $v_{j_{r-1}}$ start and end on vertices of the same
	parity. This shows that $C^k_r$ is $(X,Y)$-alternating.

	\medskip

	Suppose $C^k_r$ visits some vertex three times, that is
	\begin{equation*}
		\exists\, x<y<z\text{\ such that\ } {i_{r-1}\le x,y,z <
		j_{r-1}}\text{\ and\ }v_{\pi_r(x)}=v_{\pi_r(y)}=v_{\pi_r(z)}.
	\end{equation*}
	If $x\equiv y\pmod 2$ then $y\ge j_{r-1}$, a contradiction.
	Similarly, we must have $y\not\equiv z\pmod 2$ and $x\not\equiv z\pmod
	2$, which is a contradiction.

	If an even number of steps lead from one copy of a vertex to
	another copy of it on $C^k_r$, then $j_{r-1}$ is not minimal, a
	contradiction. This proves that $C^k_r$ is a primitive circuit.
\end{proof}

Recall Definition~\ref{def:milestone}. Let us repeat
Equation~\eqref{eq:milestone}:
\begin{equation*}
	G^k_r=X\triangle \left(\cup_{i=1}^{k-1} W^s_i\right)\triangle
	\left(\cup_{j=1}^{r-1}C^k_j\right).
\end{equation*}
In words, we obtain the milestone $G^k_r$ from $X$ by exchanging edges with
non-edges (and vice versa) in the following subsets: each $W_i$ for $1\le i\le
k-1$ and each edge in $W_k$ which is $\mathit{red}$ in $f_r$. Milestones are
special realizations: both $G^k_r\triangle X$ and $G^k_r\triangle Y$ are
subgraphs of $X\triangle Y$. Milestones are uniquely determined by $(X,Y,s)$
and the fixed lexicographical order.

\begin{lemma}\label{lemma:color}
	For any $0\le r\le \ell_k$,
	\begin{equation}\label{eq:color}
		\cup_{j=1}^{r}C^k_j=
		\left \{v_{x}v_{x+1}\,:\,\posp{x}\in
		f^{-1}_r(\mathit{red})\right\}.
	\end{equation}
	Furthermore, the Eulerian circuit described by $\pi_r$ in $W_k$ is
	alternating in $G^k_r$. For any $1\le r\le\ell_k$, the circuit $C^k_{r}$
	is alternating in $G^k_r$.
\end{lemma}
\begin{proof}
	Equation~\eqref{eq:color} follows from~\eqref{eq:rededges} and
	\eqref{eq:Ckr}. The Eulerian circuit determined by $\pi_0$ is by
	definition $(X,Y)$-alternating in $W_k$, and thus it is alternating in
	$G^k_0$. Consequently, by \eqref{eq:Ckr}, $C^k_r$ is alternating in
	$W_k$. As described by \eqref{eq:vpirx}, $\pi_r$ reverses the order of
	maximal \textit{red} segments (with respect to the coloring $f_r$). From
	the definition of $G^k_r$, it follows that the Eulerian circuit
	determined by $\pi_r$ in $W_k$ is alternating in $G^k_r$.
\end{proof}

\subsection{\ref{step3} {-} Describing the switch sequence
along a primitive circuit}\label{sec:step3}

Recall Equation~\ref{eq:path}, which describes the switch sequence from $X$
to $Y$ determined by $s$.
\begin{equation*}
	\varUpsilon(X,Y,s):=
	\left(\left(\textsc{Sweep}(G^k_r,C^k_r)\right)_{r=1}^{\ell_{k}+1}
	\right)_{k=1}^{p_s}.
\end{equation*}
In this section we provide one last missing detail of the construction, and
check that $\varUpsilon(X,Y,s)$ is indeed a path in $\G$. By
Equation~\ref{eq:milestone},
\begin{equation*}
	G^k_{r+1}=G^k_r\triangle C^k_r.
\end{equation*}
By Lemma~\ref{lemma:color}, the primitive circuit $C^k_r$ alternates in $G^k_r$.
By Lemma~\ref{th:switch-sequence}, \textsc{Sweep} returns a valid switch
sequence between $G^k_r$ and $G^k_{r+1}$. We have one degree of freedom left
in~\eqref{eq:path}: the order in which $\textsc{Sweep}$ enumerates the edges
(and vertices) of $C^k_r$ is not specified yet.
\begin{align}
	\begin{split}
	&\text{Let }x_1:=v_y\in V(C^k_r),\text{\ where $y$ minimizes}\\
	&\qquad\deg_{X[V(C^k_r)]}(v_y)+\deg_{Y[V(C^k_r)]}(v_y)-\deg_{G^k_r[V(C^k_r)]}(v_y),\\
	&\quad\text{and }y\text{\ is minimal with respect to this condition.}\\
	\end{split}\label{ass:u0}
\end{align}

Label the vertices of $C^k_r$ by $x_1,x_2,\ldots$ following the natural Eulerian
circuit~\eqref{eq:WkEulerianCircuit} on $W^k_r$ (either forwards or in reverse),
starting with $x_1x_2$, where $x_1x_2\notin G^k_r$. By
Lemma~\ref{th:switch-sequence}, $\textsc{Sweep}(G^k_r,C^k_r)$ returns a switch
sequence between $G^k_{r}$ and $G^k_{r+1}$.

\subsection{Reconstructing the endpoints \texorpdfstring{\boldmath $X,Y$}{} of the
switch sequence}\label{sec:switch-sequence}

Suppose that $Z\in\textsc{Sweep}(G^k_r,C^k_r)$ is a realization which lies on
the path $\varUpsilon(X,Y,s)$ from $X$ to $Y$ with respect to some $s\in
S_{X,Y}$. Recall Lemma~\ref{lemma:primitivecircuit} and
Equation~\eqref{eq:Rdef}, which in our setting translates to
\begin{equation}\label{eq:Rnew}
	R= \big((Z\triangle G^k_r)\setminus E(C^k_r)\big)\cup Q,
\end{equation}
where either $Q=\emptyset$ or $Q$ contains exactly one edge of $C^k_r$.
Unfortunately, $Z$ may differ from $X$ on $W_i$ even when $i\neq k$. These are
exactly the differences that are tracked by $R$, therefore it is much more
comfortable to continue working in this section with 
\begin{equation}\label{eq:R}
	Z'=Z\triangle R.
\end{equation}
Simply stated, the edges which are temporarily modified by
$\textsc{Sweep}(G^k_r,C^k_r)$ in $Z$ are returned to their appropriate state in
$Z'$. Generally $Z'$ is not a realization of $\bd$, but it is well-behaved with
respect to $W_i$ for $i\neq k$. From Definition~\ref{def:milestone} and
Equation~\eqref{eq:Rnew} it follows that
\begin{align}\label{eq:Zprime}
	\begin{split}
		Z'\triangle X,\,\, Z'\triangle Y\,\, &\subseteq E(C^k_r)\subseteq X\triangle Y, \\
		E(Z')\cap E(W_i)&=E(X)\cap E(W_i) \,\, \text{ for }\,\, i>k, \\
		E(Z')\cap E(W_i)&=E(Y)\cap E(W_i) \,\, \text{ for }\,\, i<k, \\
		E(Z')\cap E(C^k_j)&=E(X)\cap E(C^k_j) \,\, \text{ for }\,\, j>r, \\
		E(Z')\cap E(C^k_j)&=E(Y)\cap E(C^k_j) \,\, \text{ for }\,\, j<r.
	\end{split}
\end{align}
By Lemma~\ref{lemma:primitivecircuit}\ref{item:Rempty}, if $Z$ is a
\textbf{milestone} (that is, $Z=G^k_r$ for some $k,r$), then $R=\emptyset$ and
we have $Z'\cap C^k_r=Y\cap C^k_r$. The Eulerian circuit associated
to $\pi_r$ is alternating in $Z=G^k_r$.

\medskip

If $Z$ is not a milestone, then it is an \textbf{intermediate realization}
strictly between $G^k_{r}$ and $G^k_{r+1}$ on the switch sequence. There
generally does not exist an Eulerian circuit on $W_k$ which is alternating in
$Z'=Z\Delta R$. For $W_i$ ($i\neq k$), the trail induced by $s\in S_{X,Y}$ on
$W_i$ is alternating in $Z'$ (but generally it is not alternating in $Z$). Our
next goal is to slightly modify $\pi_{r-1}$, such that it induces an Eulerian
circuit in $W_k$ which is alternating in $Z'$, except at a constant number of
vertices.

\medskip

According to Lemma~\ref{lemma:primitivecircuit}\ref{item:Rtrails}, the symmetric
difference of $Z'$ and $G^k_{r}$ are two subtrails of $C^k_r$.
Without loss of generality,
\begin{align}\label{eq:twointervals}
	Z'\triangle G^k_r=\{v_x v_{x+1}\ :\ x\in\mathcal{I}
	\text{ and }f_{r-1}(\posp{x})=\mathit{green}\}
\end{align}
where $\mathcal{I}=[\alpha,\beta)\cup[\gamma,\delta)$ or
$\mathcal{I}=[i_{r-1},\alpha)\cup[\beta,\gamma)\cup [\delta,j_{r-1})$, such that
$i_{r-1}\le \alpha\le \beta\le \gamma\le \delta \leq j_{r-1}$ and
$\alpha,\beta,\gamma,\delta$ are chosen in such a way that $\mathcal{I}$ is
minimal. Let us define $\pi_{Z'}$ as follows ($\pi_{Z'}$ depends on $X$, $Y$,
and $s$, too). If $\mathcal{I}=[\alpha,\beta)\cup[\gamma,\delta)$, then let
\begin{align}\label{eq:piZprime1}
	\pi_{Z'}(x):=\left\{
		\begin{array}{ll}
			\pi_{r-1}(\overleftarrow x^{(r-1)}), & \text{ if }x\in[\alpha,\beta]\cup[\gamma,\delta] \\
			\pi_{r-1}(x), & \text{ otherwise.}
		\end{array}
	\right.
\end{align}
In other words, $\pi_{r-1}$ is reversed on the maximal $f_{r-1}$-$\mathit{red}$
intervals of $[\alpha,\beta]\cup[\gamma,\delta]$.

\medskip

If $\mathcal{I}$ does not match the previous form, then we have
$\mathcal{I}=[i_{r-1},\alpha)\cup[\beta,\gamma)\cup [\delta,j_{r-1})$, such that
$i_{r-1}<\alpha$ and $\delta<j_{r-1}$. In this case, let
\begin{align}\label{eq:piZprime2}
	\pi_{Z'}(x):=\left\{
		\begin{array}{ll}
			\pi_{r-1}\left(\overleftarrow{({i_{r-1}+j_{r-1}- x})}^{(r-1)}\right), & \text{ if }x\in[i_{r-1},\alpha]\cup[\beta,\gamma]\cup [\delta,j_{r-1}] \\
			\pi_{r-1}(i_{r-1}+j_{r-1}-x), & \text{ otherwise.}
		\end{array}
	\right.
\end{align}

\begin{lemma}\label{lemma:WkSplit}
	The Eulerian circuit defined by $\pi_{Z'}$ on $W_k$ alternates in $Z'$
	with the exception of at most 4 pairs of chords.
\end{lemma}
\begin{proof}
	Since $\pi_{r-1}$ alternates on $G^k_r$, it follows that $\pi_{r-1}$
	also alternates at $v_{\pi_{r-1}(x)}$ in $Z'$ for {$x\notin
	\mathcal{I}$}.  For any {$x\in\mathcal{I}\setminus
	\{\alpha,\beta,\gamma,\delta\}$}, $\pi_{r-1}$ alternates at $v_{\pi_{r-1}(x)}$
	in $Z'$ if $f_{r-1}(\posm{x})=f_{r-1}(\posp{x})$.

	To ensure alternation at $v_{\pi_{r-1}(x)}$ when $f_{r-1}(\posm{x})\neq
	f_{r-1}(\posp{x})$ and ${x\in\mathcal{I}\setminus
	\{\alpha,\beta,\gamma,\delta\}}$ the two subtrails are reversed on the
	edges which also belong to $f^{-1}_{r-1}(\mathit{red})$. In addition, if
	$\mathcal{I}=[i_{r-1},\alpha)\cup[\beta,\gamma)\cup [\delta,j_{r-1})$
	such that $i_{r-1}<\alpha$ and $\delta<j_{r-1}$, then we reverse the
	trail on the whole circuit $C^k_r$, so that the trail defined by
	$\pi_{Z'}$ alternates at $v_{\pi_{Z'}(i_{r-1})}$ and
	$v_{\pi_{Z'}(j_{r-1})}$ in $Z'$. Overall there are at most 4
	non-alternations of the trail associated to $\pi_{Z'}$ in $Z'$, which
	occur at a subset of $\{v_{\pi_{Z'}(\alpha)},v_{\pi_{Z'}(\beta)},
	v_{\pi_{Z'}(\gamma)},v_{\pi_{Z'}(\delta)}\}$.
\end{proof}

\begin{definition}\label{def:sigma}
	Denote by $\sigma$ the assembly structure which describes how to put
	together the Eulerian circuit defined by $\pi_{Z'}$ on $W_k$ from the at
	most three alternating subtrails of it that run between the at most four
	sites $\{v_{\pi_{Z'}(\alpha)},v_{\pi_{Z'}(\beta)},
	v_{\pi_{Z'}(\gamma)},v_{\pi_{Z'}(\delta)}\}$ of possible
	non-alternation.
\end{definition}

Ideally we want to choose a matching-system from $S_{\nabla\cap E(Z'),
\nabla\setminus E(Z')}$, but unfortunately this set is not well-defined: the
degree sequences of $\nabla\cap E(Z')$ and $\nabla\setminus E(Z')$ can be
slightly different. Recall, that Definition~\ref{def:SFF} can be applied to
graphs with different degree sequences. Instead of $\nabla\cap E(Z')$ and
$\nabla\setminus E(Z')$, we will often write $\nabla\cap Z'$ and
$\nabla\setminus Z'$.
\begin{definition}\label{def:sXYZ}
	Let $s(X,Y,Z)\in \S(\nabla\cap Z',\nabla\setminus Z')$ be the system of
	matchings we obtain by modifying the original $s\in S_{X,Y}$ such that
	it traces $\pi_{Z'}$ on $W_k$ in $Z'$ (see Lemma~\ref{lemma:trail}). The
	non-alternating pairs are not stored in $s(X,Y,Z)$.
\end{definition}

\begin{lemma}\label{lemma:numberofmatchings}
	$|\S(\nabla\cap Z', \nabla\setminus Z')|\le n^4 \cdot |S_{X,Y}|.$
\end{lemma}
\begin{proof}
	According to Lemma~\ref{lemma:WkSplit}, the sites of non-alternation of
	$\pi_{Z'}$ in $Z'$ are vertices of certain edges in $R$ ($x_1$ and the
	other ends of the start-, end- and current-chord). At each of these
	sites, $\deg_{\nabla\cap Z'}(v)=\frac12\deg_{\nabla}(v)+1$ and
	$\deg_{\nabla\setminus Z'}(v)=\frac12\deg_{\nabla}(v)-1$, or the other
	way around. Recall Equation~\eqref{eq:tFF} and~\eqref{eq:cardS}. The
	number of maximum
	matchings between edges incident to $v$ in  $\nabla\cap Z'$ and in
	$\nabla\setminus Z'$ is $\left[ (\frac12\deg_{\nabla}(v)+1)!\right]$.
	Thus there is an extra factor of $\frac12\deg_{\nabla}(v)+1\le n$
	compared to the respective factor in the enumeration of
	$S_{X,Y}=\S(X\setminus Y,Y\setminus X)$.
\end{proof}

Lemma~\ref{thm:Toperator} proves that an appropriate $w$ exists in the next
definition.
\begin{definition}\label{def:additionalParams}
Let us define the list of additional parameters:
\begin{equation}\label{eq:additionalParams}
	B(X,Y,Z,s):=(x_1,\sigma,R,w),
\end{equation}
where $w\in\mathbb{N}$ satisfies $T^w(\pi_r,\mathbf{green})=(\pi_0,g)$ for some
$g$, and $w$ is minimal.
\end{definition}

\begin{definition}
	Let $\mathbb{B}$ be the set which contains all possible tuples of
	additional parameters, defined by 
	\begin{equation}
		\mathbb{B}:=\Big\{B(X,Y,Z,s)\ :\ Z\in\varUpsilon(X,Y,s),
			\ X,Y\in\G,\ s\in S_{X,Y}\Big\}.
	\end{equation}	
\end{definition}	

\begin{lemma}\label{lemma:Bsize}
	The cardinality of the parameter set $\mathbb{B}$ is
	$\mathcal{O}(n^{8})$ in the unconstrained case, and
	$|\mathbb{B}|=\mathcal{O}(n^{6})$ in the bipartite case.
\end{lemma}
\begin{proof}
	Clearly $0\le w\le n^2$, because the generalized $T$-operator decreases
	the number of $\mathit{green}$ positions in each iteration.
	Lemma~\ref{lemma:primitivecircuit}\ref{item:Rsize} claims that given
	$x_1$, $R$ has only $\mathcal{O}(n^5)$ possible values
	($\mathcal{O}(n^3)$ in the bipartite case). Since $\sigma$ has a
	constant size description and there are at most $n$ choices for
	$x_1$, we have $|\mathbb{B}|\le \mathcal{O}(n^{8})$ (and
	$\mathcal{O}(n^{6})$ in the bipartite case).
\end{proof}
\begin{lemma}\label{lemma:sigma}
	The quadruplet composed of the graphs $Z$, $\nabla$, the matching-system $s(X,Y,Z)$ and
	$B(X,Y,Z,s)$ uniquely determines the triplet $(X,Y,s)$.
\end{lemma}
\begin{proof}
	The matching $s(X,Y,Z)$ assembles $W_i$ and the trails given by $s$ on
	them for $i\neq k$ (because the lexicographical order is used to number
	the $W_i$). The edges $E(W_k)$ are assembled into at most three trails.

	\medskip

	If $E(W_k)$ is assembled into one trail by $s(X,Y,Z)$, then
	$Z'=Z\triangle R=G^k_r$ and thus $\pi_{Z'}=\pi_{r-1}$. A single bit
	in $\sigma$ is dedicated to indicate which of the two directions the
	trail of $\pi_{r-1}$ starts from {$v_1$}. We can reconstruct $\pi_0$
	from $\pi_{r-1}$ and $w$ via Theorem~\ref{thm:Toperator}. The
	$T$-operator reproduces the primitive alternating circuits $C^k_j$ for
	$1\le j\le \ell_k$ (Lemma~\ref{lemma:primitivecircuits}). Clearly,
	$X=G^k_r\triangle\cup_{i=1}^{k-1}W_i \triangle
	\cup_{j=1}^{r-1}C^k_j$ and $Y=X\triangle \nabla$.

	\medskip

	If $E(W_k)$ is not assembled into one trail by $s(X,Y,Z)$, the
	reconstruction is trickier. Let $\sigma$ contain the list of assembly
	structures of the at most three trails (the original trail $s$ is a
	closed trail). The non-alternations in $s(X,Y,Z)$ occur precisely at the
	boundaries of intervals of positions corresponding to $\mathcal{I}$. A
	flag in $\sigma$ is dedicated to signaling whether
	Equation~\eqref{eq:piZprime1} or~\eqref{eq:piZprime2} holds for
	$\pi_{Z'}$.

	\medskip

	Suppose first, that Equation~\eqref{eq:piZprime1} holds. Then
	$\mathcal{I}=[\alpha,\beta)\cup [\gamma,\delta)$, and the values of
	$\alpha,\beta,\gamma,\delta$ are known (these are the sites of
	non-alternation). Apply the $T$-operator repeatedly to the positions in
	the interval $[\alpha,\beta)$ (i.e., restrict the $T$-operator to this
	interval), starting with an identically $\mathit{green}$ coloring until
	a fixed point of the $T$-operator is reached. According to
	Lemma~\ref{lemma:Trestricted}, this transforms
	$\pi_{Z'}|_{[\alpha,\beta]}$ back to $\pi_{r-1}|_{[\alpha,\beta]}$.
	Repeat the procedure for the $[\gamma,\delta)$ interval, too. As in the
	case when $Z$ is a milestone, $\pi_{r-1}$ and $w$ determine $\pi_0$,
	which in turn determines $G^k_r$ and then $X$ and $Y$.

	\medskip

	Lastly, if Equation~\eqref{eq:piZprime2} holds, then
	$\mathcal{I}=[i_{r-1},\alpha)\cup[\beta,\gamma)\cup [\delta,j_{r-1})$,
	such that $i_{r-1}<\alpha$ and $\delta<j_{r-1}$. Compared to the
	previous case, the extra complexity is determining $i_{r-1}$ and
	$j_{r-1}$ before we could determine $\pi_{r-1}$. Because
	$f_{r-1}((\alpha)^+)=f_{r-1}((\delta)^-)=\mathit{green}$, we have
	\begin{align*}
		j_{r-1}&=\min\{ x>\delta\ :\,\, \exists y<\alpha\,\,\, \text{ such that }\,\,\,
		(v_x=v_y)\wedge \big(x\equiv y\!\!\!\pmod{2}\big) \}, \\
			i_{r-1}&=\max\{ y<\alpha\ :\,\, (v_{j_{r-1}}=v_y)\wedge \big(j_{r-1}\equiv y\!\!\!\pmod{2}\big) \}.
		\end{align*}
\end{proof}

\begin{definition}\label{def:widehatM}
Let us define the auxiliary structure:
\begin{equation}\label{eq:widehatM}
	\begin{split}
		&\widehat M (X,Y,Z):=A_X+A_Y-A_Z, \\
		\mathcal{M}:=\Big\{&\widehat M(X,Y,Z)\ :\ Z\in
			\varUpsilon(X,Y,s),
			\ X,Y\in\G,\ s\in S_{X,Y}\Big\}.
	\end{split}
\end{equation}
where $A_X,A_Y,A_Z$ are the adjacency matrices of $X,Y,Z$, respectively.
\end{definition}
(The rows and columns of adjacency matrices are enumerated in such a way that
the matrix is symmetric and it has an identically zero main diagonal.) In this
section we only focus on the role of $\widehat M$ in the reconstruction process.
We will further discuss its properties in Section~\ref{sec:aux}. We have arrived
at the finale of the reconstruction method.
\begin{lemma}\label{lemma:Psi} The function
	\[\Phi_Z(X,Y,s):=(B(X,Y,Z,s),\widehat M(X,Y,Z),s(X,Y,Z))\]
	is injective on the tuples satisfying $Z\in \varUpsilon(X,Y,s)$. In
	other words, $\Phi_Z$ has an inverse function
	\begin{equation*}
	\Psi_Z: \mathcal{M}\times \mathbb{B}\times \bigcup_{Z',\nabla}
	\S(\nabla\cap Z',\nabla\setminus Z') \to
	\Big\{(X,Y,s):Z\in\varUpsilon(X,Y,s),
			\ X,Y\in\G,\ s\in S_{X,Y}\Big\}
	\end{equation*}
	(naturally, $\Psi_Z$ is only defined on the image set of $\Phi_Z$).
\end{lemma}
\begin{proof}
	The graph $\nabla=X\triangle Y$ is determined by $Z$ and $\widehat M$,
	but this information alone does not separate the $X$-edges and
	$Y$-edges. The graph $Z'$ is determined by $B(X,Y,Z,s)$ and $Z$.
	Lemma~\ref{lemma:sigma} claims that $(X,Y,s)$ can be reconstructed from
	these objects.
\end{proof}

\section{Directed degree sequences}\label{sec:directed}

Now it is time to extend our management for directed degree sequences. This
short description goes more or less in parallel with~\cite{linear} (Erdős,
Mezei, Miklós and Soltész, 2018).

Let $\vec G$ be a simple directed graph (parallel edges and loops are forbidden,
but oppositely directed edges between two vertices are allowed) with vertex set
$X(\vec G) = \{x_1,x_2,  \ldots,x_n \}$ and edge set $E( \vec G)$. For every
vertex $x_i\in X$ we associate two numbers: the {\em out-degree\/} and the {\em
in-degree\/} of $x_i$. These numbers form the directed degree bi-sequence $\vec
\bd =(\vec \bd_\textrm{out},\vec \bd_\textrm{in})$.

We introduce the following {\em bipartite representation\/} of $\vec G$: let
$B({\vec G})=(U,V;E)$ be a bipartite graph where each class consists of one copy
of every vertex from $X(\vec G)$. The edges adjacent to a vertex $u_x$ in class
$U$ represent the out-edges from $x$, while the edges adjacent to a vertex $v_x$
in class $V$ represent the in-edges to $x$ (so a directed edge $xy$ corresponds
to the edge $u_x v_y$). If a vertex has zero in- (respectively out-) degree in
$\vec G$, then we delete the corresponding vertex from $B({\vec G})$. (This
representation was used by Gale~\cite{G57}, but one can find it already
in~\cite{pet} (Petersen, 1981).) The directed degree bi-sequence $\vec \bd$
gives rise to a  bipartite degree sequence $\vec \bD$:
\begin{equation*}
	\vec{\bD}=\big((\vec\bd_\textrm{out}(x_1),\ldots,\vec\bd_\textrm{out}(x_n)),
	(\vec\bd_\textrm{in}(x_1),\ldots,\vec\bd_\textrm{in}(x_n))\big).
\end{equation*}

Since there are no loops in our directed graph, there cannot be any $(u_x,v_x)$
edge in its bipartite representation: these vertex pairs are {\bf
non-chords}. It is easy to see that these forbidden edges form a forbidden
(partial) matching $\F$ in the bipartite graph $B(\vec G)$, or in more general
terms, in $B(\vec \bD)$, and we call this a \textbf{restricted} bipartite degree
sequence.
\begin{definition}
	For restricted bipartite degree sequences, the set of chords is the
	vertex pairs of the form ${u_x}{v_y}$ where $x\neq y$.
\end{definition}
By definition, $\G(\vec\bD)$ is the set of all bipartite realizations of
$\vec\bD$ which avoid the non-chords from $\F$. Now it is easy to see that the
bipartite graphs in $\G(\vec\bD)$ are in one-to-one correspondence with the
possible realizations of the directed degree bi-sequence.

Consider two oppositely oriented triangles, $\overrightarrow{C_3}$ and
$\overleftarrow {C_3}$. In order to move between these two realizations of the
same degree sequence, Kleitman and Wang~\cite{KW73} observed that a new
operation is needed, and they proved that with this extra ``switching''
operation the space of realizations becomes irreducible. Take the symmetric
difference $\nabla$ of the bipartite representations $B(\overrightarrow{C_3})$
and $B(\overleftarrow{C_3})$. It contains exactly one alternating cycle (the
edges come alternately from $B(\overrightarrow{C_3})$ and
$B(\overleftarrow{C_3})$), s.t.\ each vertex pair of distance 3 along the cycle
in $\nabla$ is a non-chord. In this alternating cycle no ``classical'' switch
can be performed. To address this issue, we need an extra switch operation,
which is the bipartite analogue of the operation introduced by Kleitman and
Wang: we exchange all edges coming from $B(\overrightarrow{C_3})$ with all edges
coming from $B(\overleftarrow{C_3})$ in one operation.

In general, a \textbf{triple-switch} is defined as follows: take a length-6
alternating cycle $C$ in $\nabla$, and if one of the three vertex pairs of
distance 3 in $C$ forms a non-chord, we exchange all edges of $C$ to non-edges
and vice versa.  It is a well-known fact (\cite{EKM} (Erdős, Király and Miklós,
2013),~\cite{EKMS} (Erdős, Kiss, Miklós and Soukup, 2015)) that the set
$\G(B(\vec \bD))$ of all realizations is irreducible under switches and
triple-switches that avoid the $\F$-edges.

The example of $\overrightarrow{C_3}$ and $\overleftarrow {C_3}$ demonstrates
why the triple-switch operation is necessary. However, as long as some steps of
the Markov-chain require choosing 6 vertices, it seems wasteful to not perform
the triple-switch simply because some of the vertex pairs of distance 3 are
chords.

\medskip

In this paper, we relax the restrictions on triple-switches: given a length-6
alternating cycle $C$ in $\nabla$, a triple switch is valid if and only if at
least one of the three vertex pairs of distance 3 in $C$ is a non-chord. This
relaxation allows us to shave off a factor of $n^4$ from the mixing time of the
Markov chain. To see this, compare the proofs of Theorem~\ref{thm:unconstrained}
and Theorem~\ref{thm:directed}.

\medskip

The inner loop of Algorithm~\ref{alg:sweep} has to be modified, because the
conclusion of Lemma~\ref{th:C6} does not necessarily hold in the directed case.
The adaptation of \textsc{Sweep} in Algorithm~\ref{alg:directedsweep} works on
the bipartite representation $B(\vec G)$ instead of the directed graph $\vec G$.
\begin{algorithm}
	\caption{Sweeping a primitive circuit in the bipartite representation.
	The \textsc{Directed Sweep} assumes that $x_1x_2\notin E(G)$ and that
	$C=(x_1,x_2,\ldots,x_{2\ell})$ is a primitive alternating
	circuit.}\label{alg:directedsweep}
	\begin{algorithmic}
		\Function{Triple-switch}{$G,x_1,{[x_{2t-2},x_{2t-1},x_{2t},x_{2t+1},x_{2t+2}]}$}
		\State\Return $G+\{x_{1}x_{2t-2},x_{2t-1}x_{2t},x_{2t+1}x_{2t+2}\}-\{x_{2t-2}x_{2t-1},x_{2t}x_{2t+1},x_{2t+2}x_1\}$
		\EndFunction
		\Statex
		\Procedure{Directed Sweep}{$G,{[x_{1},x_{2},\ldots,x_{2\ell}]}$} $\to [Z_{1},Z_{2},\ldots, Z_{\ell-2}, \langle Z_{\ell-1}\rangle ]$
		\State $Z_0\gets G$
		\State $q\gets 1$
		\State $\mathit{end}\gets 2$
		\State $\mathcal{L}\gets\left\{ 2i\in2\mathbb{N} :\ 4\le 2i\le 2\ell,\,\,  x_1x_{2i}\text{ is a chord and }x_1x_{2i}\in E(G)\right\}$
		\Statex
		\While{$\mathit{end}<2\ell$}
		\State $\mathit{start}\gets \min\big\{ 2i\in\mathcal{L}\ :\ 2i>\mathit{end} \big\}$
		\State $2t\gets \mathit{start}-2$
		\Statex
		\While{$2t\ge\mathit{end}$}
		\If{$x_1x_{2t}$ is a non-chord}
		\State $Z_{q}\gets \Call{Triple switch}{Z_{q-1},x_1,{[x_{2t-2},\ldots,x_{2t+2}]}}$\label{state:tripleswitch}
		\State $2t\gets 2t-2$
		\ElsIf{$x_1x_{2t}$ is a chord}
		\State $Z_q\gets Z_{q-1}-\big\{x_1x_{2t+2},x_{2t}x_{2t+1}\big\}+\big\{x_1x_{2t},x_{2t+1}x_{2t+2}\}$
		\EndIf
		\State $q\gets q+1$
		\State $2t\gets 2t-2$
		\EndWhile
		\Statex
		\State $\mathit{end}\gets \mathit{start}$
		\EndWhile
		\Statex
		\EndProcedure
	\end{algorithmic}
\end{algorithm}
If $Z_q$ gets its value from \textsc{Triple-switch}, then
Lemma~\ref{lemma:primitivecircuit}\ref{item:Rdoublestepsecondgraph} applies to
it, otherwise Lemma~\ref{lemma:primitivecircuit}\ref{item:Rswitch} holds for
$Z_q$. Because of this, the statements of
Lemma~\ref{lemma:primitivecircuit}\ref{item:Rtrails} and~\ref{item:Rsize}, and
Lemma~\ref{lemma:Psi} about the bipartite case apply to the directed case as
well.

\medskip

We are ready to define our switch Markov chain on $(\G(\vec\bD), P)$ for the
restricted bipartite degree sequence $\vec\bD$.  The transition (probability)
matrix  $P$ of the  Markov chain is defined as follows: let the current
realization be $G$. Then
\begin{enumerate}[label=(\alph*)]
	\item with probability $1/2$ we uniformly choose a set of two vertices
	$\{u,u'\}$ from $U$ and a set of two vertices $\{v,v'\}$ from $V$. There
	are two matchings, $\{uv,u'v'\}$ and $\{uv',u'v\}\}$, between these sets. Let
	$F$ be one of these matchings, chosen randomly, and let $F'$ be the other
	matching. If both $F$ and $F'$ consist of chords only and $F\subseteq E(G)$ and
	$F'\cap E(G)=\emptyset$, then perform the switch (so $E(G')=(E(G)\cup
	F')\setminus F$), otherwise $G'=G$.

	\item With probability $1/2$ we choose a set of three vertices from $U$
		and a set of three vertices from $V$. Let $F$ and $F'$ be a
		uniformly randomly selected pair of disjoint perfect matchings
		between these sets. If both $F$ and $F'$ consist of chords only,
		and the remaining matching between the two chosen sets contains
		a non-chord, and $F\subseteq E(G)$ and $F'\cap E(G)=\emptyset$,
		then perform the triple-switch (so $E(G')=E(G)\cup F'\setminus
		F$), otherwise $G'=G$.
\end{enumerate}
The (triple-)switch moving from $G$ to $G'$ is unique, therefore the probability
of this transformation (the {\em jumping probability\/} from $G$ to $G'\ne G$)
is:
\begin{equation}
	\mathrm{Prob}(G \rightarrow_{(a)}   G'):= P(G,G') = \frac{1}{4} \cdot
	\frac{1}{\binom{|U|}{2} \binom{|V|}{2}},
\end{equation}
and
\begin{equation}
	\mathrm{Prob}(G\rightarrow_{(b)} G'):= P(G,G') = \frac{1}{24} \cdot
	\frac{1}{\binom{|U|}{3} \binom{|V|}{3} }.
\end{equation}
The probability of transforming $G$ to $G'$ (or vice versa) is time-independent
and symmetric. Therefore, $P$ is a symmetric matrix, where the entries in the
main diagonal are non-zero, but (possibly) distinct values. Again, $P(G,G)\ge
\frac12$, because if $(F,F')$ corresponds to a feasible (triple-)switch, then
$(F',F)$ does not. Therefore the chain is aperiodic and the eigenvalues of its
transition matrix are non-negative. Each transition $(X,Y)\in E(\G(\vec\bD))$
satisfies $P(X,Y) \ge n^{-4}$.

However it is important to recognize that in papers~\cite{G11} (Greenhill, 2011)
and~\cite{G18} (Greenhill and Sfragara, 2018) a slightly different
Markov chain is studied, where it is assumed that the degree sequences under study
are irreducible using switches only. This is the case, for example, for regular
directed degree sequence. Papers~\cite{BM10} (Berger and Müller{-}Hannemann,
2010) and~\cite{lamar} (LaMar, 2011) provide a full characterization of directed
degree sequences with this property.

\section{The auxiliary matrix \texorpdfstring{$\widehat M$}{hat M}}\label{sec:aux}

The auxiliary matrix $\widehat M=A_X+A_Y-A_Z$ defined in~\eqref{eq:widehatM} is
a linear combination of three adjacency matrices. The row and columns sums are
equal to the corresponding degrees prescribed by $\bd$. If $Z=G^k_r$, then
$G^k_r\triangle X\subseteq X\triangle Y$ implies that $\widehat M$ is a 0--1
matrix. If $Z$ is an intermediate realization, $\widehat M$ is still a 0--1
matrix except on the entries associated to edges in $R$, since $(Z\triangle
R)\triangle X\subseteq X\triangle Y$. These $+2$ and $-1$ entries will be called
\textbf{bad entries}, and the chords to which they correspond to are called
type-$(2)$ and type-$(-1)$ chords, respectively.

\medskip

\begin{lemma}\label{lemma:plustwominusone}
	If $R$ falls under
	case~{\rm\ref{item:Rswitch}}~or~{\rm\ref{item:Rdoublestepsecondgraph}}
	of Lemma~\ref{lemma:primitivecircuit}, then $R$ contains at most two
	type-$(2)$ and at most one type-$(-1)$ chords.
\end{lemma}
\begin{proof}
	Lemma~\ref{lemma:primitivecircuit}\ref{item:Rswitch}~or~\ref{item:Rdoublestepsecondgraph}
	claims that $R$ has at most three elements. Of these,
	$x_1x_\mathit{start}$ and $x_1x_\mathit{end}$ are edges in $X$, so the
	entries associated to them in $\widehat M$ are {symmetric pairs of} $+2$
	or $+1$ entries. In case~\ref{item:Rswitch}, if $R$ contains the third
	chord, $x_1x_{2t}$, and it is an edge in $X$, then we must have
	$\mathit{end}=2t$, so $R$ actually does not contain $x_1x_{2t}$. Thus
	$x_1x_{2t}\in R\implies x_1x_{2t}\notin E(X)$, so the entries associated
	to $x_1x_{2t}$ in $\widehat M$ are $-1$'s or $0$'s.
	Case~\ref{item:Rdoublestepsecondgraph} is similar to
	case~\ref{item:Rswitch}.
\end{proof}

\medskip

The \textbf{switch} operation is extended to symmetric matrices as follows.
Suppose $M\in \mathbb{Z}^{[k]\times [k]}$. For any $x,y\in [k]$ we define the
one-edge graph $G^{x,y}=([k];\{xy\})$ with the adjacency matrix $A_{xy}$.
Clearly, $A_{xy}$ is a symmetric matrix with two $1$'s. Let $(x,y;z,w)$ be a
list of four \textbf{pairwise distinct} elements of $[k]$. Switching along these
four vertices produces the symmetric matrix
\begin{equation}\label{eq:mswitch}
	M'=M+A_{xz}-A_{zy}+A_{yw}-A_{wx}.
\end{equation}
Clearly, the row and column sums of $M'$ are identical to that of $M$. Notice,
that a switch in $Z$ translates into a switch on $\widehat M$.

\medskip

Notice, that for bipartite degree sequences, the ``top-right'' submatrix of this
$\widehat M$ is equal to the auxiliary matrix used in~\cite{MES} (Miklós, Erdős
and Soukup, 2013) (the bipartite adjacency matrix).

\begin{lemma}\label{lemma:removingplustwos}
	Let $M\in \mathbb{Z}^{[k]\times [k]}$ be a symmetric matrix with 0's in
	the diagonal, such that each row and column sum is in the interval
	$[1,k-2]$. Also, suppose that the row sum of the first row is minimal.
	If the entries of $M$ are 0 and 1, except for at most two {symmetric}
	pairs of entries of $+2$ in the first row and in the first column, and
	at most {one symmetric pair} or $-1$ entries anywhere in the matrix,
	then there exist at most $2$ switches that transform $M$ into a 0--1
	matrix except for at most {one pair of symmetric $-1$ entries}.
\end{lemma}
\begin{proof}
	Suppose $M_{1,j}=2$. We must have $j\neq 1$, which means that the
	maximum of an entry in the rest of the column of $j$ is 1. Because the
	column sum is at most $k-2$ and there is at most one $-1$ in the column,
	there exist $i\in [k]\setminus\{1,j\}$ such that $M_{i,j}\in \{-1,0\}$.
	We have two cases.
	\begin{enumerate}
		\item[{\bf (I)}] Suppose that there exists $\ell\in
			[k]\setminus\{1,i\}$ such that $M_{i,\ell}> M_{1,\ell}$.
			Since $1\neq i,\ell$, we assume that $M_{i,\ell}<2$,
			therefore $M_{1,\ell}\in \{0,-1\}$.  Switch along
			$(1,i;\ell,j)$ in $M$. The operation decreases $M_{1,j}$
			to $1$. If $M_{i,\ell}=0$ then $M_{1,\ell}=-1$, so when
			the switch creates a symmetric pair of $-1$'s, it also
			eliminates another pair. The matrix resulting from the
			switch operation satisfies the assumptions of this lemma
			and contains two fewer $+2$ entries.
		\item[{\bf (II)}] Otherwise, for all $\ell\in
			[k]\setminus\{1,i\}$ we have $M_{i,\ell}\le M_{1,\ell}$.
			Since the row sum of the first row is minimal, we have
			\[
			\sum_{\ell=1}^k M_{1,\ell}  \le \sum_{\ell=1}^k M_{i,\ell}                       \]
			and hence
			\[
				0  \le \sum_{\ell\in
				\{1,i,j\}}\left(M_{i,\ell}-M_{1,\ell}\right)=
				M_{i,1}-M_{1,i}+M_{i,j}-M_{1,j}=M_{i,j}-2,
			\]
			because $M$ is symmetric with 0 diagonal. The inequality
			implies that $M_{i,j}=2$, so either $i=1$ or $j=1$,
			which contradicts our choice of $i,j$.
	\end{enumerate}
	By recursion a second pair of entries which equal $+2$ can also be eliminated.
\end{proof}

\section{Applications of the unified method}\label{sec:appl}

In this Section we harvest some fruits of our unified machinery,
proving a rather general result for all typical degree sequence types.

In 1990 Jerrum and Sinclair published a very influential paper~\cite{JS90}
(Jerrum and Sinclair, 1990) about fast uniform generation of regular graphs and
about realizations of degree sequences where no degree exceeds $\sqrt{n/2}$.  To
achieve this goal, they applied the Markov chain they have developed
in~\cite{JS89} (Jerrum and Sinclair, 1989). Informally it is known as {\bf JS
chain}, and it is sampling the perfect and near-perfect 1-factors on the
corresponding Tutte gadget. The rapid mixing nature of the JS chain depends on
the ratio of the number of perfect and the number near-perfect 1-factors. As
they proved it is applicable if and only if the degree sequence $\bd$ belongs to
a \textbf{\boldmath{$P$}-stable} class.

\medskip

Recall the definition of $P$-stability (introduced in
Definition~\ref{def:Pstable}).  Careful examination of the known results about
rapidly mixing switch Markov chains revealed the fact that all known ``good''
degree sequence classes (for \UC\, bipartite or directed degree sequences) are
$P$-stable.  It raises the conjecture that the switch Markov chains on
$P$-stable degree classes are rapidly mixing. We resolve this conjecture
affirmatively in this section.

\medskip

For a fixed $B\in \mathbb{B}$, let the set of compatible auxiliary structures be
\begin{align*}
	\mathcal{M}_B= & \left\{ \widehat M\ :\ \exists X,Y,Z\in \G(\bd),\ s\in S_{X,Y} \right . \\
		       & \left . \text{ s.t. }\widehat M=A_X+A_Y-A_Z,\ B=B(X,Y,Z,s) \right\}.
\end{align*}

\medskip

To apply the simplified Sinclair's method, it is sufficient to estimate from above the value of
\begin{equation}\label{eq:UCcount1}
	\sum_{X,Y \in V(\G)}
	\frac{\left| \big \{s\in S_{X,Y}: Z\in \varUpsilon(X,Y,s)\big \} \right| } {|S_{X,Y}|}
\end{equation}
for any realization $Z$.
\begin{lemma}\label{lemma:computation}
	The following bound holds for any unconstrained, bipartite, and directed degree sequence:
	\begin{equation*}
		\sum_{X,Y \in V(\G)}
		\frac{\left| \big \{s\in S_{X,Y}: Z\in \varUpsilon(X,Y,s)\big \} \right| } {|S_{X,Y}|}\le n^4\cdot \sum_{B\in \B}|\mathcal{M}_B|.
	\end{equation*}
\end{lemma}
\begin{proof}
	According to Lemma~\ref{lemma:Psi}, expression~\eqref{eq:UCcount1} can be rewritten as follows:
	\begin{equation}\label{eq:UCcount2}
		\sum_{X,Y \in V(\G)}
		\frac{\left| \big \{\Psi_Z(\widehat M(X,Y,Z),B(X,Y,Z,s),s(X,Y,Z)) : Z\in\varUpsilon(X,Y,s)\big \} \right| } {|S_{X,Y}|}
	\end{equation}
	Observe, that $|S_{X,Y}|$ is already determined by $\nabla=X\triangle Y$, which in turn is determined by $Z$ and $\widehat M$. Let $t_\nabla:=|S_{X,Y}|$. Furthermore, $Z'$ is determined by $B$ and $Z$ (see Equations~\eqref{eq:R}~and~\eqref{eq:additionalParams}). Let
	\begin{equation*}
		\B_{\widehat M}=\left\{ B(X,Y,Z,s)\ :\exists X,Y\text{ such that } \widehat M=A_X+A_Y-A_Z,\ s\in S_{X,Y} \right\}.
	\end{equation*}
	Continue writing~\eqref{eq:UCcount2} as follows and apply Lemma~\ref{lemma:numberofmatchings}.
	\begin{align*}
& =\sum_{\widehat M\in \mathcal{M}}
\frac{\left| \big \{\Psi_Z(\widehat M,B,s^*)\ :\ \exists B\in\B,\ s^*\in \S(\nabla\cap Z',\nabla\setminus Z')\big \} \right| }{t_\nabla}\le                                                      \\
& \le \sum_{\widehat M\in \mathcal{M}}\frac{|\B_{\widehat M}|\cdot n^4\cdot t_\nabla}{t_\nabla}\le n^4\cdot\sum_{\widehat M\in \mathcal{M}}|\B_{\widehat M}| \le n^4\cdot \sum_{B\in \B}|\mathcal{M}_B|.
	\end{align*}
\end{proof}

At this point, the proofs for unconstrained, bipartite, and directed degree sequences slightly diverge. The most general of these is the case of unconstrained degree sequences. First, we discuss this case. Having understood the argument, it is relatively simple to fit it to the cases of the bipartite and directed degree sequence cases. Moreover, the tools required for proving our results on the latter two classes have already been published in~\cite{MES} (Miklós, Erdős and Soukup, 2013), so their proofs will be less verbose than the next section on unconstrained degree sequences.

\subsection{Unconstrained degree sequences}

First, let us bound the number of auxiliary structures compatible with a given
parameter set. Recall Definition~\ref{def:Pstable}. The $x^\text{th}$ unit
vector is denoted by $\mathds{1}_x$.

\begin{lemma}\label{lemma:cardcalM}
	If the stability of an unconstrained degree sequence $\bd$ is bounded by the polynomial $p(n)$, then
	\begin{equation*}
		\left| \mathcal{M}_B \right|\le n^{6}\cdot p(n)\cdot |\G(\bd)|
	\end{equation*}
	holds for any $B\in \mathcal{B}$.
\end{lemma}
\begin{proof}
	Equation~\eqref{eq:param} defines $B=(x_1,\sigma,R,w)$. Recall, that $(Z\triangle R)\triangle X\subseteq X\triangle Y$, so the bad entries ($+2$ and $-1$ values) in $\widehat M$ correspond to positions assigned to chords in $R$. Let $M$ be the symmetric submatrix of $\widehat M$ induced by the vertices of $C^k_r$ as rows and columns. We have two cases.

	\subsection*{Case 1: \boldmath $\centernot{\exists} f\in R$ where $x_1\notin f$
	and $f\notin E(C^k_r)$}

	All of the non 0--1 entries of $\widehat M$ are contained in $M$, in the
	rows and columns associated to $x_1$. Since $Z$ is an intermediate
	realization of the switch sequence from $G^k_r$ to $G^k_{r+1}$, the
	degree sequence of $Z[V(C^k_r)]$ is equal to the degree sequence of
	$G^k_r[V(C^k_r)]$. Therefore, Lemma~\ref{lemma:plustwominusone} and
	Assumption~\eqref{ass:u0} implies that we can use
	Lemma~\ref{lemma:removingplustwos} to remove the $+2$'s from $M$ with at
	most two switches. Hence the same applies to $\widehat M$. For each
	switch, the type-$(2)$ chord determines two vertices of the switch, thus
	there are $n^4$ ways to choose the at most two switches that eliminate
	the $+2$ entries.

	\medskip

	Let $\widehat M'$ be the matrix we get after applying the switches
	defined by Lemma~\ref{lemma:removingplustwos}. Either $\widehat M'$ is
	an adjacency matrix of a realization of $\bd$, or $\widehat M'$ contains
	$-1$ entries at positions associated to the chord $xy$. In the former
	case $\widehat M'\in \G(\bd)$, and in the latter $\widehat
	M'+A_{xy}+A_{yx}\in \G(\bd+\mathds{1}_x+\mathds{1}_y)$.

	\subsection*{Case 2: \boldmath $\exists! f\in R$ such that $x_1\notin f$
	and $f\notin E(C^k_r)$}

	This is only possible if $R$ falls under case~\ref{item:Rdoublestep} of
	Lemma~\ref{lemma:primitivecircuit}. The auxiliary structure belonging to
	the intermediate realization before or after $Z$ on the switch sequence
	is one switch away from $\widehat M$, moreover this switch touches $f$.
	There are at most $n^2$ switches satisfying these conditions, because
	$f$ already determines two vertices. After performing the appropriate
	switch on $\widehat M$, the enumeration of the previous case applies.
\end{proof}

We are ready to prove one of the main results of this paper.

\begin{theorem}\label{thm:unconstrained}
	The switch Markov chain is rapidly mixing on $P$-stable unconstrained
	degree sequence classes.
\end{theorem}
\begin{proof}
	From Lemma~\ref{lemma:cardcalM} and Lemma~\ref{lemma:Bsize} we get:
	\begin{align*}
		n^4\cdot \sum_{B\in \B}|\mathcal{M}_B|\le
		n^4\cdot |\mathbb{B}|\cdot  n^{6}\cdot p(n)\cdot|\G(\bd)|\le
		\mathcal{O}(n^{18})\cdot p(n)\cdot|\G(\bd)|.
	\end{align*}
	Equation~\eqref{eq:multiflow_simple2} now follows from
	Lemma~\ref{lemma:computation}. Every condition of the simplified
	Sinclair's method is satisfied, so the switch Markov chain on $\G(\bd)$
	is rapidly mixing. 
	Theorem~\ref{th:simplified_sinclair} states that the mixing time is
	\begin{equation*}
		\tau_{\varepsilon}(\G(\bd),P)\le n^4\cdot m\cdot
		\mathcal{O}(n^{18})\cdot p(n)\cdot
		\left(n^2-\log\varepsilon\right)\le
		\mathcal{O}(n^{22})\cdot p(n)\cdot m\cdot (n^2-\log\varepsilon).
	\end{equation*}
	where the length of $\varUpsilon(X,Y,s)$ is bounded by Lemma~\ref{lemma:pathLength}.
\end{proof}

\subsection{Bipartite degree sequences}

Let $\bD$ denote a bipartite degree sequence on $n=n_1+n_2$ vertices.

\begin{definition}\label{def:bipPstable}
	Let $\mathcal{D}$ be an infinite set of bipartite degree sequences. We say that $\mathcal{D}$ is \textbf{\boldmath$P$-stable}, if there exists a polynomial $p\in \mathbb{R}\left[x\right]$ such that for any $n_1, n_2\in \mathbb{N}, n_1\ge n_2,$ and any  degree sequence $\bD \in \mathcal{D}$ on $n_1 \mbox{ and } n_2$ vertices we have
	\begin{equation*}
		\left|\G(\bD)\cup \left( \bigcup_{x\in [n_1],\ y\in [n_2]}
		\G(\bD+\mathds{1}_x+\mathds{1}_{(n_1+y)}) \right)\right|\le p(n)\cdot |\G(\bD)|,
	\end{equation*}
	where $\mathds{1}_x$ is the $x^\text{th}$ unit vector.
\end{definition}

Recall that a primitive alternating circuit on a bipartite graph is a cycle. Also, for any $X,Y,Z\in \G(\bD)$, the auxiliary structure $\widehat M=A_X+A_Y-A_Z$ is determined by the submatrix spanned by $U\times V\subset {(U\uplus V)}^2$. This is the ``top-right'' submatrix, often called the bipartite adjacency matrix. The ``top-left'' and the ``bottom-right'' submatrices are zero.

\begin{lemma}\label{lemma:cardcalMbip}
	If the stability of a bipartite degree sequence $\bD$ is bounded by the polynomial $p(n)$, then
	\begin{equation*}
		\left| \mathcal{M}_B \right|\le n^4\cdot p(n)\cdot |\G(\bD)|
	\end{equation*}
	holds for any $B\in \mathcal{B}$.
\end{lemma}
\begin{proof}
	The proof is simpler and slightly different than that of Lemma~\ref{lemma:cardcalM}. \textsc{Double step} is never called in the bipartite case (Lemma~\ref{th:switch-sequence}), so either $R$ is empty or Lemma~\ref{lemma:plustwominusone} applies to it. Hence, there are ${(n_1n_2)}^2$ possibilities to choose the {remaining vertices of the (at most two)} switches that eliminate the type-$(2)$ bad chords.

	\medskip

	Secondly, we have to make sure that the switches produced by Lemma~\ref{lemma:removingplustwos} respect the bipartition.
	As before, let $M$ be the submatrix of $\widehat M$ induced by the vertices of $C^k_r$. Let $H=K_{U(C^k_r)}\uplus K_{V(C^k_r)}$ be the disjoint union of the two cliques within the classes. Instead of applying Lemma~\ref{lemma:removingplustwos} on $M$, apply it on $M+A_H$. Each row and column sum increased by the same number, therefore assumptions of the lemma are still satisfied. Any switch which eliminates a $+2$ from this matrix which is valid in the unconstrained sense also respects the bipartition.
\end{proof}
\begin{theorem}\label{thm:bipartite}
	The switch Markov chain is rapidly mixing on $P$-stable bipartite degree sequence classes.
\end{theorem}
\begin{proof}
	Instead of Lemma~\ref{lemma:cardcalM} we use Lemma~\ref{lemma:cardcalMbip}. The bound on the size of $\B$ is $\mathcal{O}(n^6)$ according to Lemma~\ref{lemma:Psi}. Other than these differences, the proof is identical to that of Theorem~\ref{thm:unconstrained}:
	\begin{equation*}
		n^4\cdot \sum_{B\in \B}|\mathcal{M}_B| \le  n^4\cdot |\mathbb{B}|\cdot  n^4\cdot p(n)\cdot|\G(\bD)|\le \mathcal{O}(n^{14})\cdot p(n)\cdot|\G(\bD)|.
	\end{equation*}
	Theorem~\ref{th:simplified_sinclair} states that the mixing time is
	\begin{equation*}
		\tau_{\varepsilon}(\G(\bD),P)\le n^4\cdot m\cdot
		\mathcal{O}(n^{14})\cdot p(n)\cdot
		\left(n^2-\log\varepsilon\right)\le
		\mathcal{O}(n^{18})\cdot p(n)\cdot m\cdot (n^2-\log\varepsilon).
	\end{equation*}
	where the length of $\varUpsilon(X,Y,s)$ is bounded by Lemma~\ref{lemma:pathLength}.
\end{proof}

\subsection{Directed degree sequences}

Recall from Section~\ref{sec:directed} that instead of directly manipulating directed graphs, we work on their bipartite representations. Formally, the degree sequence of the directed graph is identical to that of its bipartite representation.
Through the bipartite representation, directed degree sequence classes inherit a definition of $P$-stability (thus the realizations of the representing [perturbed] bipartite degree sequences also avoid the non-chords).
Let $\vec \bD$ denote the bipartite representation of a directed degree sequence $\vec \bd$ on $n$ vertices (so the length of the vector $\vec\bD$ is $2n$).

\begin{lemma}\label{lemma:cardcalMdir}
	If the stability of a directed degree sequence $\vec\bd$ is bounded by the polynomial $p(n)$, then
	\begin{equation*}
		\left| \mathcal{M}_B \right|\le n^4\cdot p(n) \cdot |\G(\vec\bd)|=n^4\cdot p(n)\cdot |\G(\vec \bD)|
	\end{equation*}
	holds for any $B\in \mathcal{B}$.
\end{lemma}
\begin{proof}
	The proof of Lemma~\ref{lemma:cardcalMbip} applies to the bipartite representation, but we have to check that applying Lemma~\ref{lemma:removingplustwos} on $M+A_H$ produces switches that avoid the non-chords. Indeed, this is the case, because the non-chords of the form ${u_x}{v_x}$ correspond to the main diagonal in $M$, which the switches chosen by the lemma avoid.
\end{proof}
\begin{theorem}\label{thm:directed}
	The switch Markov chain is rapidly mixing on $P$-stable directed degree sequence classes.
\end{theorem}
\begin{proof}
	By Lemma~\ref{lemma:Psi}, the bound on the size of $\B$ is $\mathcal{O}(n^6)$, as in the proof of Theorem~\ref{thm:bipartite}. From the previous lemma, we get:
	\begin{align*}
		n^4\cdot \sum_{B\in \B}|\mathcal{M}_B| \le n^4\cdot |\mathbb{B}|\cdot  n^{4}\cdot p(n)\cdot|\G(\vec\bD)|\le \mathcal{O}(n^{14})\cdot p(n)\cdot|\G(\vec\bD)|.
	\end{align*}
	Theorem~\ref{th:simplified_sinclair} states that the mixing time is
	\begin{equation*}
		\tau_{\varepsilon}(\G(\bD),P)\le n^6\cdot m\cdot
		\mathcal{O}(n^{14})\cdot p(n)\cdot
		\left(n^2-\log\varepsilon\right)\le
		\mathcal{O}(n^{20})\cdot p(n)\cdot m\cdot (n^2-\log\varepsilon).
	\end{equation*}
	where the length of $\varUpsilon(X,Y,s)$ is bounded by Lemma~\ref{lemma:pathLength}.
\end{proof}

\section{\texorpdfstring{$P$}{P}-stable degree sequence classes}\label{sec:regions}

In the proof of almost every previous result on rapid mixing of the switch
Markov chain, it turns out there is a short hidden proof that the degree
sequences under study are $P$-stable. The unified proof contains most of the
technical difficulty of proving rapid mixing of the switch Markov chain.

There have already been successful attempts at unifying some of the proofs, most
notably by Amanatidis and Kleer~\cite{AK19,AK20}, who study the notion of strong
stability:

\begin{definition}[adapted from~\cite{AK19}]\label{def:stronglystable}
	Let $\mathcal{D}$ be a set of degree sequences. Let
	\begin{equation*}
		\mathbb{G}'(\bd)=\bigcup_{x,y\in [n]}
		\mathbb{G}(\bd-\mathds{1}_x-\mathds{1}_y).
	\end{equation*}
	We say that $\mathcal{D}$ is strongly stable if there exists a constant
	$\ell$ such that for any $\bd\in \mathcal{D}$ and any $G'\in \G'(\bd)$
	there exists $G\in \G(\bd)$ (which depends on $G'$) such that
	$|E(G')\Delta E(G)|\le 2\ell$.
\end{definition}
For bipartite graphs, the definition is analogous.  The above definition is
easily seen to be equivalent with the one given in~\cite{AK19}, and it has the
advantage that it does not rely on the definition of the Jerrum-Sinclair chain.
Having formally defined strong stability, we restate the relevant theorem of
Amanatidis and Kleer.
\begin{theorem}[\cite{AK19,AK20}]\label{thm:AK20}
	The switch Markov chain is rapidly mixing on strongly stable
	unconstrained and bipartite degree sequence classes.
\end{theorem}

In the following subsections of this section we discuss all known $P$-stable
degree sequence regions. It is an intriguing problem to discover other
$P$-stable regions.

\subsection{Unconstrained degree sequences}

For the sake of having more readable and compact formulas, let $\Delta=\max\bd$, $\delta=\min\bd$, and $m=\frac12 \sum_{v\in V}\bd(v)$ be functions of $\bd$.

\medskip

Recently, Greenhill~and~Sfragara~\cite{G18} published a breakthrough result on the rapid mixing of the switch Markov chain.

\begin{theorem}[\cite{G18}]\label{thm:GS18unconstrained}
	The switch Markov chain is rapidly mixing on the following family of unconstrained degree sequences:
	\begin{equation}\label{eq:GS}
		\mathcal{D}_\text{GS}:= \left\{\bd\in \mathbb{Z}^+\ :\ \delta\ge 1,\ 3\le \max\bd\le\frac13\sqrt{2m} \right\}
	\end{equation}
\end{theorem}

It turns out that the authors implicitly prove on page~10 of~\cite{G18} that $\mathcal{D}_\text{GS}$ is a $P$-stable class. However, this implicit result is actually not new: Jerrum, McKay, and Sinclair extensively studied the notion of $P$-stability in their seminal work~\cite{JMS92}.

\begin{theorem}[Jerrum, McKay, Sinclair -- Theorem 8.1 in~\cite{JMS92}]\label{thm:JMS92}
	The family of unconstrained degree sequences
	\begin{equation*}
		\mathcal{D}_\text{JMS}:=\left\{\bd\in \mathbb{N}^n\ :\ {\left(\max\bd-\min\bd+1\right)}^2\le 4\cdot\min\bd\cdot(n-\max\bd-1) \right\}
	\end{equation*}
	is $P$-stable.
\end{theorem}

\begin{theorem}[Jerrum, McKay, Sinclair -- Theorem 8.3 in~\cite{JMS92}]\label{thm:JMS92average}
	The family of unconstrained degree sequences
	\begin{align*}
		\mathcal{D}_\text{JMS+}:=\Big\{\bd\in \mathbb{N}^n\ :\  & (2m-n\delta)(n\Delta-2m)\le \\
									& \le (\Delta-\delta)\big( (2m-n\delta)(n-\Delta-1)+(n\Delta-2m)\delta \big) \Big\},
	\end{align*}
	is $P$-stable.
\end{theorem}

Theorem~\ref{thm:unconstrained} implies that the switch Markov chain is rapidly mixing on elements of $\mathcal{D}_\text{JMS}$ and $\mathcal{D}_\text{JMS+}$. Moreover, it is easy to see that $\mathcal{D}_\text{GS}\subset\mathcal{D}_\text{JMS+}$. However, the proofs of Theorems~\ref{thm:JMS92}~and~\ref{thm:JMS92average} actually prove a bit more than just $P$-stability.
In~\cite{JMS92} it is also shown that $\mathcal{D}_\text{JMS}$ and $\mathcal{D}_\text{JMS+}$ are \textbf{strongly stable} regions with $\ell\le 10$, so Theorem~\ref{thm:AK20} already applies to them.

\medskip

The following corollary is a consequence of the fact that the degrees in an Erdős-Rényi random graph are tightly concentrated around their expected value.

\begin{corollary}\label{cor:ER}
	Let $G(n,p)$ be an Erdős-Rényi random graph of order $n\ge 100$ with edge probability $p$,
	where $p$ is bounded away from 0 and 1 by at least $\frac{5\log n}{n-1}$.
	Then $\Pr(\bd(G(n,p))\in\mathcal{D}_\text{JMS+})\ge 1-\frac{3}{n}$.
\end{corollary}
\begin{proof}
	We may suppose that $p\le \frac12$ by taking the complement of $G$ if necessary.
	Let $p=p(n)$, $\varepsilon_1=\sqrt{\frac{\log {n} }{n-1}}$ and $m=\frac12\sum_{v\in V}\bd(v)=\binom{n}{2}(p+\varepsilon_2)$.
	By Hoeffding's inequality, we have
	\begin{align*}
	& \Pr\big(\Delta(G)> (p+\varepsilon_1)\cdot (n-1)\big)\le                                                                       \\
	& \qquad\le\sum_{v\in V(G)}\Pr\Big(\bd(v)< (p+\varepsilon_1)\cdot(n-1)\Big)\le n\cdot e^{-2\varepsilon_1^2(n-1)}\le \frac{1}{n},
	\end{align*}
	and similarly
	\[ \Pr\big(\delta(G)< (p-\min(\varepsilon_1,p))\cdot (n-1)\big)\le \frac{1}{n}.
	\]
	The degree sequence $\bd(G(n,p))$ is in $\mathcal{D}_\text{JMS+}$ if it satisfies
	\begin{equation}\label{eq:UC-ER-condition}
		(2m-n\delta)(n\Delta-2m)\le (\Delta-\delta)\left( (2m-n\delta)(n-\Delta-1)+(n\Delta-2m)\delta \right).
	\end{equation}
	First suppose that $p\ge \varepsilon_1$.  Because increasing $\Delta$ or decreasing $\delta$ makes the inequality stricter, without loss of generality, we may substitute $\Delta=(p+\varepsilon_1)\cdot(n-1)$ and $\delta=(p-\varepsilon_1)\cdot(n-1)$ into \eqref{eq:UC-ER-condition}. We calculate
	\begin{align*}
		(2m-n\delta)             & =\left(2\binom{n}{2}(p+\varepsilon_2)-n(p+\varepsilon_1)(n-1)\right)=n(n-1)(\varepsilon_1+\varepsilon_2), \\
		(n\Delta-2m)    & =n(n-1)(\varepsilon_1+\varepsilon_2),   \\
		(2m-n\delta)(n-\Delta-1) & = n(n-1)(\varepsilon_1+\varepsilon_2)\cdot (1-p-\varepsilon_1)(n-1),  \\
		(n\Delta-2m)\delta       & = n(n-1)(\varepsilon_1+\varepsilon_2)\cdot (p-\varepsilon_1)(n-1).
	\end{align*}
	Therefore~\eqref{eq:UC-ER-condition} holds if
	\[
	{\left(n(n-1)(\varepsilon_1+\varepsilon_2)\right)}^2  \le 2\varepsilon_1(n-1)\cdot n(n-1)(\varepsilon_1+\varepsilon_2)\cdot (1-2\varepsilon_1)(n-1), \]
	or after simplification,
	\[
		n\varepsilon_2    \le (n-2)\varepsilon_1-4\varepsilon_1^2 (n-1).
	\]
	If $\varepsilon_2\le\frac{\sqrt{\log n}}{n-1}$ then the last inequality is satisfied whenever
	\begin{align*}
		n\frac{\sqrt{\log n}}{n-1} & \le (n-2)\sqrt{\frac{\log n}{n-1}}-4\log n.
	\end{align*}
	Clearly, the right hand side grows $\varTheta(n^{\frac12})$ faster than the left hand side as $n\to \infty$, and the inequality already holds for $n=100$.
	Now for any $c > 0$,
	\begin{equation*}
		\Pr\left(\frac12\sum_{v\in V} d(v)> (p + c)\binom{n}{2} \right)\le e^{-2 c^2\binom{n}{2}},
	\end{equation*}
	and substituting $c=\frac{\sqrt{\log n}}{n-1}$ shows that
	$\varepsilon_2\le \frac{\sqrt{\log n}}{n-1}$ with probability at least $1-1/n$.
	Overall, $\Pr(\bd(G(n,p))\notin\mathcal{D}_\text{JMS+})\le \frac{1}{n}+\frac{1}{n}+\frac{1}{n}$ if $p\ge \varepsilon_1$.

	\medskip

	Now suppose that $\frac{5\log n}{n-1}\le p<\varepsilon_1$. Substituting $\Delta=(p+\varepsilon_1)\cdot(n-1)$ and $\delta=0$ into \eqref{eq:UC-ER-condition}, we find that $\bd(G(n,p))\in\mathcal{D}_\text{JMS+}$
	if
	\[
		2m(n\Delta-2m) \le \Delta\cdot  2m(n-\Delta-1).
	\]
	Rearranging, this inequality holds if $\Delta(\Delta+1) \le 2m$, and substituting for $\Delta$
	and simplifying gives the sufficient condition
	$4\varepsilon_1^2 \le p-\varepsilon_2.$

	The last inequality is satisfied if $\varepsilon_2=\frac{\sqrt{\log n}}{n}$. Therefore, if $\frac{5\log n}{n-1}\le p<\varepsilon_1$, then $\Pr(\bd(G(n,p))\notin\mathcal{D}_\text{JMS+})\le \frac{2}{n}$.
\end{proof}

Similar results have already been proved for bipartite Erdős-Rényi
graphs~\cite{linear-pre,linear} (Erdős, Mezei, Miklós and Soltész,
2018), with the requirement that $p$ is bounded away from 0 and 1 by at least
$4\sqrt{\frac{2\log n}{n}}$.

\subsection{Unconstrained power-law bounded degree sequences}\label{sec:power}
Let us quote two definitions introduced by Gao~and~Wormald (2016).
\begin{definition}[\cite{GW16}]\label{def:powerlaw}
	Suppose $\bd\in \mathbb{N}^n$ is a degree sequence.
	If $\exists C>0$, then $d$ is
	\begin{itemize}
		\item \textbf{power-law density-bounded} with parameter $\gamma$, if
			 for all $i\in[1,n]$,
			 \begin{equation*}
				|d^{-1}(i)|\le Cni^{-\gamma}
			 \end{equation*}
		\item \textbf{power-law distribution-bounded} with parameter
			$\gamma$, if for all $i\in[1,n]$
			\begin{equation*}
				\sum_{j=i}^n\left|d^{-1}(j)\right|\le
			\sum_{j=i}^\infty Cnj^{-\gamma}.
			\end{equation*}
	\end{itemize}
\end{definition}

Barabási and Albert (1999) \cite{BA99} recognized in their seminal paper that a
lot of real world networks grow via some form of preferential attachment, and
these networks have power-law like degree distributions. The preferential
attachment model and its relatives produce graphs with power-law
\emph{density-bounded} degree sequences (see Definition~\ref{def:powerlaw}).
However, the degree sequences of most real world networks deviate somewhat from
the degree sequences of such synthetic networks.  Instead, as noted in
\cite{GW16}, the degree sequence of a real world network is much more likely to
obey the less restrictive power-law \emph{distribution-bound}. Compared to the
former bound, the latter allows relatively high maximum degrees and longer tails
in the degree distribution.

The switch Markov chain is not the only way to exactly sample the uniform
distribution on the realizations of a degree sequence. Recently, Gao~and~Wormald
(2018) presented in~\cite{GW18} the first ``provably practical'' sampler for
power-law distribution-bounded degree sequence where $\gamma$ is allowed to be
less than 3; in fact they can go as low as 2.8811. For such degree sequences,
they provide a \emph{linear time} approximate sampler and a polynomial time
exact sampler.

In degree distributions of empirical networks following a power-law, the
parameter $\gamma$ is usually between 2 and 3.

\medskip

Gao and Wormald~\cite{GW16} compute the number of realizations for several types
of heavy-tailed degree sequences, and in-turn, those formulas imply
$P$-stability of the respective classes. They conjectured that the degree
sequences obeying a power-law distribution-bound with $\gamma>2$ are $P$-stable,
which was shown by Gao~and~Greenhill (2020) in~\cite{GG20}. This is a
corroborative example for the applicability of Theorem~\ref{thm:unconstrained},
because rapid mixing of such degree sequences is immediately verified,
independently from the rapid mixing result of~\cite{GG20}.

\subsection{Bipartite degree sequences}

Let $\bD$ be a bipartite degree sequence on $U$ and $V$ as vertex classes. We use the following shorthands in this sub-section:
\begin{align*}
	\delta_U & =\min_{u\in U}\bD(u), & \delta_V & =\min_{v\in V}\bD(v), \\
	\Delta_U & =\max_{u\in U}\bD(u), & \Delta_V & =\min_{v\in V}\bD(v),
\end{align*}
and $m=\sum_{u\in U} \bD(u)=\sum_{v\in v} \bD(v)$.

\begin{theorem}[implicitly proved in Theorem 2 in~\cite{linear} Erdős,
	Mezei, Miklós, and Soltész, 2018]\label{th:lin-root}
	The set of bipartite degree sequences $\bD$ that satisfy
	\begin{equation}\label{eq:biproot}
		2\le \Delta\le\sqrt{\frac{m}{2}},
	\end{equation}
	is $P$-stable.
\end{theorem}
Clearly, Theorems~\ref{thm:GS18unconstrained}~and~\ref{th:lin-root} are closely related, the difference in constants is caused by the different structural constraints only.

\begin{theorem}[implicitly proved in Theorem 3 in~\cite{linear}]\label{thm:bipmax}
	The set of bipartite degree sequences $\bD$ that satisfy
	\begin{equation}\label{eq:bipmax}
		(\Delta_U-\delta_U-1)(\Delta_V-\delta_V-1)\le \max\Big(\delta_U(|U|-\Delta_V+1),\delta_V(|V|-\Delta_U+1)\Big)
	\end{equation}
	is $P$-stable.
\end{theorem}

Amanatidis and Kleer presented a bipartite analogue of Theorem~\ref{thm:JMS92}.

\begin{theorem}[Corollary~18 in~\cite{AK19}]\label{thm:bipAK}
	The set of bipartite degree sequences that satisfy both
	\begin{equation}
		\begin{split}
			(\Delta_U-\delta_V)^2 \le 4\delta_V\cdot (|V|-\Delta_U) \\
			(\Delta_V-\delta_U)^2 \le 4\delta_U\cdot (|U|-\Delta_V)
		\end{split}\label{eq:bipAK}
	\end{equation}
	is $P$-stable $($because it is strongly stable$)$.
\end{theorem}

The following Theorem~\ref{thm:bip4min} is a bipartite analogue of
Theorem~\ref{thm:JMS92}. In some sense, it is stronger than either
Theorem~\ref{thm:bipmax} or~\ref{thm:bipAK}. If one side is regular
($\Delta_U=\delta_U$), then Inequality~\eqref{eq:bipmax} and~\eqref{eq:linstab}
are automatically satisfied. Inequality~\eqref{eq:bipmax} trivially holds even
for almost half-regular bipartite degree sequences ($\Delta_U\le \delta_U+1$).
Assuming $|U|=|V|$, $\Delta_U=\Delta_V$, $\delta_U=\delta_V$ are all
satisfied,~\eqref{eq:bipAK} and~\eqref{eq:linstab} are equivalent, and are
loosely speaking 4 times better than~\eqref{eq:bipmax}.
\begin{theorem}\label{thm:bip4min}
	The set of bipartite degree sequences that satisfy
	\begin{equation}\label{eq:linstab}
		(\Delta_U-\delta_U)\cdot (\Delta_V-\delta_V)\le
		4\cdot \min\Big(\delta_U(|U|-\Delta_V),\delta_V(|V|-\Delta_U)\Big)
	\end{equation}
	is $P$-stable (because it is strongly stable).
\end{theorem}
\begin{proof}
	If the degrees of two vertices in the same class are increased by one in
	a graphic bipartite degree sequence, then the resulting degree sequence
	is not graphic. Let $G$ be a realization of
	$\bD+\mathds{1}_{u_1}+\mathds{1}_{v_1}$ on $U$ and $V$ as vertex
	classes. The degree sequence of $G$ is
	\begin{align*}
		\bD(G)=\left\{
			\begin{array}{ll}
				d(x)+1 & \text{ if }x=u_1\text{ or }v_1, \\
				d(x)   & \text{ otherwise.}
			\end{array}
		\right.
	\end{align*}

	We claim that there exists an alternating path $P$ of length at most 7
	between $u_1$ and $v_1$ in $G$, such that the first and last edges of
	$P$ are edges of $G$. Assume that no such path exists. Let $U_1,V_2,U_3$
	be the set of vertices that are reachable from $v_1$ via an alternating
	path (starting with an edge of $v_1$ in $G$) of length exactly $1$, $2$,
	and $3$, respectively. Define $V_1,U_2,V_3$ similarly with respect to
	$u_1$.

	\medskip

	By our assumption, $u_1\notin U_1$, $v_1\notin V_1$, otherwise
	$\{u_1,v_1\}\in E(G)$ is an alternating path of length 1. If
	$G[U_1,V_1]$ is not a complete bipartite graph, then there is an
	alternating path of length 3 which is a good candidate for $P$.
	Similarly, $G[U_2,V_2]$ is an empty graph (no alternating paths of
	length 5), and $G[U_3,V_3]$ is a complete bipartite graph (no
	alternating paths of length 7). These observation also imply that
	$\{u_1\}\uplus U_1\uplus U_2\uplus U_3$ and $\{v_1\}\uplus V_1\uplus
	V_2\uplus V_3$ are subpartitions of $U$ and $V$, respectively. Let $U_4$
	and $V_4$ be the remaining vertices of $U$ and $V$, respectively.

	\medskip

	By definition and the fact that $G[U_3,V_3]$ is a complete bipartite
	graph, $G[U_1\cup U_3,V_1\cup V_3]$ is also a complete bipartite graph.
	The vertices in $U_2$ are only adjacent to elements of $V_1\cup V_3$,
	therefore
	\begin{align*}
		|E(U_2,V_1\cup V_3)|=|E(U_2,V)| & \ge \delta_U|U_2|.
	\end{align*}
	Every vertex which is joined by a non-edge to a vertex of $V_1$ is
	contained in $U_2$. Therefore $G[U_4,V_1]$ is a complete bipartite graph
	and $|V_1|=d(u_1)\ge \delta_U+1$, thus
	\begin{align*}
		|E(U_4,V_1\cup V_3)|\ge|E(U_4,V_1)| & > \delta_U|U_4|.
	\end{align*}
	Since $G[U_1\cup U_3,V_1\cup V_3]$ is a complete bipartite graph, we have
	\begin{align*}
		|E(U_2\cup U_4,V_1\cup V_3)| & \le |V_1\cup V_3|\cdot
		\left(\Delta_V-|U_1\cup U_3|\right).
	\end{align*}
	Combining the previous inequalities, we get
	\begin{align*}
		\delta_U(|U|-|U_1\cup U_3|-1) & =\delta_U\cdot |U_2\cup U_4|<
		|V_1\cup V_3|\cdot \left(\Delta_V-|U_1\cup U_3|\right).
	\end{align*}
	Let us substitute $k_1=|U_1\cup U_3|$ and $k_2=|V_1\cup V_3|$, leading to
	\begin{align*}
		\delta_U(|U|-\Delta_V-1) & < (\Delta_V-k_1)\cdot (k_2-\delta_U), \\
		\delta_V(|V|-\Delta_U-1) & < (\Delta_U-k_2)\cdot (k_1-\delta_V).
	\end{align*}
	The second inequality is obtained by symmetry. Solving for $k_1$, we get
	\begin{align*}
		\frac{\delta_V(|V|-\Delta_U-1)}{\Delta_U-k_2}+\delta_V & < k_1 <
		\Delta_V-\frac{\delta_U(|U|-\Delta_V-1)}{k_2-\delta_U}.
	\end{align*}
	Without loss of generality, we may assume that
	$\delta_V(|V|-\Delta_U-1)\le \delta_U(|U|-\Delta_V-1)$. Omitting $k_1$
	from the middle of the above inequality, we have
	\begin{align*}
		\delta_V(|V|-\Delta_U-1)(\Delta_U-\delta_U) & <
		(\Delta_V-\delta_V)(\Delta_U-k_2)(k_2-\delta_U).
	\end{align*}
	The left hand side in the last inequality is maximal if
	$k_2=\frac12(\Delta_U+\delta_U)$. We get
	\begin{align*}
		4\delta_V(|V|-\Delta_U-1) & < (\Delta_V-\delta_V)(\Delta_U-\delta_U),
	\end{align*}
	which contradicts the assumptions of this theorem. Thus there exists a
	suitable alternating path of length at most 7 starting on an edge of
	$u_1$ and ending on an edge of $v_1$.

	\medskip

	Switching the edges along the alternating path $P$ transforms $G$ into a
	realization of $\bD$. The procedure consists of at most $8$
	non-deterministic choices (vertices of the alternating path), so
	$p(|U|+|V|)={(|U|+|V|)}^8$ is a good witness to stability.
\end{proof}

\subsection{Directed degree sequences}

Let $\vec \bd$ be a directed degree sequence on $X$ as vertices. Let $\vec\bd_\mathrm{out}$ be the out-degree sequence and $\vec\bd_\mathrm{in}$ be the in-degree sequence. We use the following abbreviations in this sub-section:
\begin{align*}
	\delta_\mathrm{out} & =\min_{x\in X}\vec\bd_\mathrm{out}(x), & \delta_\mathrm{in} & =\min_{x\in X}\vec\bd_\mathrm{in}(x), \\
	\Delta_\mathrm{out} & =\min_{x\in X}\vec\bd_\mathrm{out}(x), & \Delta_\mathrm{in} & =\min_{x\in X}\vec\bd_\mathrm{in}(x),
\end{align*}
and $m=\sum_{x \in X} \vec\bd_\mathrm{out}(x)=\sum_{x \in X} \vec\bd_\mathrm{in}(x)$.

\begin{theorem}[implicitly proved in~\cite{G18} Greenhill and Sfragara, 2018]
	The set of bipartite degree sequences $\vec\bd$ that satisfy
	\begin{equation}
		2\le \max(\Delta_\mathrm{out},\Delta_\mathrm{in})\le\frac14\sqrt{m},
	\end{equation}
	is $P$-stable.
\end{theorem}

\begin{theorem}[implicitly proved in Theorem 4 in~\cite{linear} Erdős, Mezei, Miklós and Soltész, 2018]
	The set of directed degree sequences $\vec\bd$ satisfying
	\begin{equation*}
		2\le\max(\Delta_\mathrm{out},\Delta_\mathrm{in})< \frac{1}{\sqrt{2}}\sqrt{m-4},
	\end{equation*}
	is $P$-stable.
\end{theorem}

\begin{theorem}[implicitly proved in Theorem 5 in~\cite{linear}]
	The set of directed degree sequences $\vec\bd$ satisfying
	\begin{align*}
		(\Delta_\mathrm{out}-\delta_\mathrm{out}) & \cdot (\Delta_\mathrm{in}-\delta_\mathrm{in})\le 2-n+                                                                                                                                  \\
		+                                         & \max\Big( \delta_\mathrm{out}(n-\Delta_\mathrm{in}-1)+\delta_\mathrm{in}+\Delta_\mathrm{out}, \delta_\mathrm{in}(n-\Delta_\mathrm{out}-1)+\delta_\mathrm{out}+\Delta_\mathrm{in} \Big)
	\end{align*}
	is $P$-stable.
\end{theorem}

\section{Summary}\label{sec:summary}

To summarize the new results of the paper we present
Table~\ref{table:newresults},  an updated version of Table~\ref{table:results}
which contains both entirely new and improved results. A strongly stable class
is, as the name suggests, naturally $P$-stable, see the papers of Amanatidis and
Kleer~\cite{AK19,AK20}. Their results already provide a unified framework for proving
all previously known bipartite and \UC\ sequence results.

\begin{table}[ht]
	\begin{center}
		\renewcommand*{\arraystretch}{1.5}
		\begin{tabular}{| c | c | c | }
			\hline
			\UC\ sequences                                                                         & bipartite deg.\ seq.                                                       & directed deg.\ seq.                                               \\ \hline \hline
			\color{gray!60!black}
			regular~\cite{CDG07}                                                                   & \color{gray!60!black} {(half-)}regular~\cite{MES}                                   & \color{gray!60!black} regular~\cite{G11}                                   \\ \hline
			\cellcolor{gray!50!white}	& \multicolumn{2}{|c|}{\color{gray!60!black} almost half regular~\cite{EKMS}}                                                                             \\ \hline
			\color{gray!60!black} $\Delta \leq \frac{1}{3}\sqrt{2m}$~\cite{G18}                             & \color{gray!60!black} $\Delta \leq \frac{1}{\sqrt{2}} \sqrt{m}$~\cite{linear}       & \color{gray!60!black} $\Delta< \frac{1}{\sqrt{2}}\sqrt{m-4}$~\cite{linear} \\ \hline
			\color{gray!60!black}Power-law distribution- & \cellcolor{gray!50!white} & \cellcolor{gray!50!white} \\
			\color{gray!60!black}bound,
			$\gamma>2$~\cite{GG20}\textsuperscript{$\dagger$} & \cellcolor{gray!50!white}& \cellcolor{gray!50!white}\\ \hline
			\color{gray!60!black} ${(\Delta-\delta+1)}^2\le$                                                & \small $(\Delta_U-\delta_U)\cdot (\Delta_V-\delta_V)\le$                                      & \color{gray!60!black} \small similar to bipartite case                                         \\
			\color{gray!60!black} $\le 4\cdot \delta(n-\Delta-1)$                               &  \small $4\delta_U(|U|-\Delta_V),4\delta_V(|V|-\Delta_U)$    & \color{gray!60!black} \cite{linear-pre, linear}                                                                   \\

			\color{gray!60!black} proof in~\cite{AK19,AK20}& Theorems~\ref{thm:bipmax}~and~\ref{thm:bipAK} &  \\ \hline
			Erdős-Rényi $G(n,p)$                                                                   & \small \color{gray!60!black} Bipartite Erdős-Rényi~\cite{linear-pre,linear}                                   & \color{gray!60!black} \small similar to bipartite case                           \\
			$p,1-p\ge \frac{5{\log n}}{n-1}$  & \color{gray!60!black} $p,1-p\ge 4\sqrt{\frac{2\log n}{n}}$ & \color{gray!60!black} \cite{linear-pre,linear}                             \\ \hline
			\multicolumn{2}{|c|}{\color{gray!60!black} strongly
		stable degree sequence classes~\cite{AK19,AK20}} & \cellcolor{gray!50!white}\\ \hline
			\multicolumn{3}{|c|}{\textbf{\boldmath$P$-stable degree sequence classes}}                                                                                                                                                              \\ \hline
		\end{tabular}
	\end{center}
	\caption{Updated version of Table~\ref{table:results} with the new
		results in this paper. Here $\Delta$ and $\delta$ denote the
		maximum and minimum degrees, respectively. Half of the sum of
		the degrees is $m$, and $n$ is the number of vertices. The
		notation is similar for bipartite and directed degree sequences.
		Some technical conditions have been omitted. Gray text is used
		for previously known results. \newline
		\textsuperscript{$\dagger$} Gao and Greenhill~\cite{GG20} list two
		sufficient conditions for some kind of stability, both of which
		apply to power-law distribution-bounded degree sequences with
		$\gamma>2$.  Rapid mixing is shown through \texttt{Condition 1}
		in~\cite{GG20}. \texttt{Condition 2} of~\cite{GG20} has slightly better
		constants and implies rapid mixing in combination with
		Theorem~\ref{thm:unconstrained} or Theorem~\ref{thm:AK20}.} \label{table:newresults}
\end{table}

The flexibility of our unified method allowed us to extend the rapid mixing
results of the switch Markov chain in two directions (in the table):
\begin{itemize}
	\item vertically (power of machinery) to $P$-stable degree sequence
		classes, and
	\item horizontally (applicability of machinery) to directed degree
		sequences.
\end{itemize}

We have also shown that the degree sequence of the Erdős-Rényi random graph
$G(n,p)$ is rapidly mixing with high probability as $n\to \infty$, for any edge
probability $p$ satisfying $p,1-p\ge \frac{5{\log n}}{n-1}$.

\medskip

The notion of $P$-stability {arises naturally when studying} the rapid mixing of
the switch Markov chain~\cite{JS90,JMS92} (Jerrum and Sinclair, 1990, Jerrum,
McKay and Sinclair, 1992). It would be really intriguing to find even a small
rapidly mixing degree sequence class which is not $P$-stable.  Finding the
bipartite and directed analogues of Theorem~\ref{thm:JMS92average} seems to be a
relatively easy and moderately rewarding open problem.

\section*{Acknowledgement}
We would like to thank Pieter Kleer and Yorgos Amanatidis for their remarks
regarding stable bipartite degree sequences.
We also thank the referees for their helpful comments.

\bibliographystyle{plain}

\end{document}